\documentclass[12pt,leqno,color,amssymb]{amsart}

\usepackage{amssymb, amsmath,latexsym, amsfonts}

\def\emty{\emptyset}
\def\cbar{\widehat{\C}}

\def\stab{\mbox{\rm Stab}}
\def\diam{\mbox{\rm diam}}

\def\dist{\mbox{\rm dist}}

\def\dis{\displaystyle}
\newcommand\ben{\begin{enumerate}}
	\newcommand\een{\end{enumerate}}
\newcommand\bit{\begin{itemize}}
	\newcommand\eit{\end{itemize}}

\def\AAA{{\mathcal A}}

\def\CCC{{\mathcal C}}

\def\EEE{{\mathcal E}}

\def\GGG{{\mathcal G}}
\def\HHH{{\mathcal H}}

\def\KKK{{\mathcal K}}

\def\OOO{{\mathcal O}}

\def\QQQ{{\mathcal Q}}
\def\RRR{{\mathcal R}}

\def\TTT{{\mathcal T}}

\def\th{theorem }

\def\homeo{homeomorphism }

\def\Proof{{\noindent\sc Proof.} }

\def\al{\alpha}
\def\be{\beta}
\def\te{\theta}

\def\g{\gamma}
\def\G{\Gamma}

\def\vp{\varphi}
\def\ep{\varepsilon}
\def\la{\lambda}
\def\La{\Lambda}
\def\De{\Delta}
\def\de{\delta}

\newcommand{\R}{\mathbb R}
\newcommand{\HH}{\mathbb H}
\newcommand{\C}{\mathbb C}
\newcommand{\TT}{\mathbb T}
\newcommand{\Z}{\mathbb Z}
\newcommand{\Q}{\mathbb Q}
\newcommand{\N}{\mathbb N}
\renewcommand{\SS}{\mathbb S}
\newcommand{\DD}{\mathbb D}
\newcommand{\PP}{\mathbb P}
\newcommand{\BB}{\mathbb B}
\newcommand{\E}{\mathbb E}

\def\endp{\hspace*{\fill}$\rule{.55em}{.55em}$ \smallskip}

\newtheorem{newthm}{Theorem}
\newtheorem{theorem}{Theorem}[section]
\newtheorem{lemma}[theorem]{Lemma}
\newtheorem{proposition}[theorem]{Proposition}
\newtheorem{corollary}[theorem]{Corollary}
\newtheorem{defthm}[theorem]{Definition et \th}
\newtheorem{defn}[theorem]{Definition}
\newcommand{\REFEQN}[1] { \begin{equation}\label{#1} }
	\newcommand{\ENDEQN}{\end{equation}}
\newcommand{\REFTHM}[1] { \begin{theorem}\label{#1} }
	\newcommand{\ENDTHM}{\end{theorem}}
\newcommand{\REFNTH}[1] { \begin{newthm}\label{#1} }
	\newcommand{\ENDNTH}{\end{newthm}}
\newcommand{\REFPROP}[1]{\begin{proposition}\label{#1} }
	\newcommand{\ENDPROP}{\end{proposition} }
\newcommand{\REFLEM}[1]{\begin{lemma}\label{#1} }
	\newcommand{\ENDLEM}{\end{lemma} }
\newcommand{\REFCOR}[1]{\begin{corollary}\label{#1} }
	\newcommand{\ENDCOR}{\end{corollary} }
\newcommand{\REFDEFTHM}[1] { \begin{defthm}\label{#1} }
	\newcommand{\ENDDEFTHM}{\end{defthm}}

\usepackage{enumerate}
\usepackage{color}
\usepackage[dvipsnames]{xcolor}

\numberwithin{equation}{section}

\renewcommand{\hat}{\widehat}
\renewcommand{\tilde}{\widetilde}
\newcommand{\fraction}{\varepsilon(x,y)}

\newtheorem{rmk}[theorem]{Remark}
\newtheorem{fact}[theorem]{Fact}

\newtheorem{cl}[theorem]{Claim}
\newtheorem*{theorem*}{Theorem}

\newcommand{\fdefeq}{\overset{\text{def.}}{=}}

\def\cay{\mathcal{C}\!\mbox{\it ay}\,}
\def\cus{\mathcal{C}\!\mbox{\it us}\,}
\def\bop{\partial_\PP}
\newcommand{\demode}[1]{\noindent{\sc Proof of #1}.}

\title{Quasi-isometric rigidity of three-manifold groups
}
\author{Peter Ha\"{\i}ssinsky}
\address{Aix Marseille Univ, CNRS, I2M, UMR 7373, Marseille, France}
\email{phaissin@math.cnrs.fr}
\author{Cyril Lecuire}
\address{UMPA UMR 5669 CNRS,
ENS de Lyon Site Monod,
46 All\'ee d'Italie,
69364 Lyon Cedex 07}
\email{cyril.lecuire@ens-lyon.fr}

\date{\today}

\begin{document}

\date{}

\maketitle


\begin{abstract} We provide a   proof  that the classes of finitely generated Kleinian groups and of three-manifold groups are 
 quasi-isometrically rigid.
\end{abstract}

\tableofcontents

\section{Introduction} 

It was already known to Dehn that any finitely presented group can be realized as the fundamental group of a closed manifold of any 
dimension at least four, but this is not the case
for three-manifolds, e.g., the group $\Z^4$ is not the fundamental group of any closed three-manifold. This makes the class of three-dimensional manifolds special and we may expect that their fundamental groups enjoy specific properties
which characterize them among finitely generated groups.

From the point of view of geometric group theory, one tries to understand the properties of a group by studying the different
actions it admits on metric spaces. For the action of the group $G$ on the geodesic metric space $X$ to properly reflect the properties of $G$, 
we require that it is {\it geometric}: the group $G$ acts by isometries (the action is distance-preserving), properly discontinuously 
(for any compact subsets $K$ and $L$ of $X$, at most finitely many elements $g$ of $G$ will satisfy $g(K)\cap L \ne\emty$) and cocompactly 
(the orbit space $X/G$ is compact). By identifying $G$ with the orbit $Go$ of a point $o\in X$ and by pulling back the induced metric, 
we obtain a metric on $G$. Changing the orbit or the metric space $X$ gives rise to new metrics. Thus we get a metric structure on $G$ which is coarsely defined in the following sense.

A {\em quasi-isometry} between metric spaces $X$ and $Y$ is a coarsely bi-Lipschitz coarsely surjective map  $\vp:X\to Y$, i.e., there are constants $\la>1$ and $c>0$ such that:
\bit
\item (quasi-isometric embedding) for all $x,x'\in X$, the two inequalities $$\frac{1}{\la}d_X(x,x') - c \le d_Y(\vp(x),\vp(x'))\le \la\, d_X(x,x') + c$$
hold and
\item (quasi-surjectivity) the $c$-neighborhood of the image $f(X)$ covers $Y$.
\eit
This defines in fact an equivalence relation on (separable) metric spaces. The \v{S}varc--Milnor lemma asserts that there is  only one geometric action 
of a group on a proper geodesic metric space up to quasi-isometry \cite[Prop.\,3.19]{ghys:delaharpe:groupes}: we equip $G$ with a reference metric induced by identifying $G$ with an orbit $G o$ under a geometric action on some proper geodesic space $Y$ (usually one takes its left action on one of its locally finite Cayley graphs).

\noindent{\bf \v{S}varc--Milnor Lemma.} 
{\it Let $X$ be a proper geodesic metric space.  Let $G$ act properly discontinuously and cocompactly on $X$ by isometries. 
Then $G$ is finitely generated and, for any $x_0\in X$, the map $g\mapsto g.\ x_0$ is a quasi-isometry.}

Thus we have a coarsely well-defined metric on $G$, i.e., defined up to quasi-isometry, coming from its geometric actions. 
We are then naturally led to ask whether or not a property of a finitely generated group is invariant under quasi-isometries or, equivalently, 
whether or not a class of groups is {\it quasi-isometrically rigid}: a class of groups $\mathcal{C}$ is  quasi-isometrically rigid if any group quasi-isometric to a group in $\mathcal{C}$ has a finite index subgroup in $\mathcal{C}$.

Since the pioneering works of Stallings and Gromov, diverse classes of groups have been proved to be quasi-isometrically rigid. We should mention free groups \cite{stallings:yale}, nilpotent groups \cite{gromov:expanding}, Abelian groups \cite{gromov:expanding, bass:degree, guivarch:croissance, pansu:croissance,cornulier:tessera:valette:gafa07} and word hyperbolic groups. A more thorough overview of these results will be given in \textsection \ref{sec:qi 3mfds}.

In contrast with higher dimension (every finitely presented group is the fundamental group of a compact $4$-manifold), fundamental groups of low-dimensional compact manifolds have many restrictive properties and their quasi-isometric rigidity is a challenging question that has already led to many interesting developments. In this paper we focus on compact three-manifolds but the reader should be aware that the quasi-isometric rigidity of surface groups follows from \cite{gabai:S1}, \cite{casson:jungreis} and results we have already mentioned (see \textsection  \ref{sec:qi 3mfds} for more details). Our main theorem completes the work of many people whose combined results can be fairly accurately summarized in the following statement; see \textsection \ref{sec:qi 3mfds} for more details about these results and \textsection \ref{sec:qi 3mfds} and \textsection \ref{3 mfds} for the definition of irreducible three-manifolds.

\begin{theorem} \label{thm:zero exact seq}
Let $G$ be a group quasi-isometric to the fundamental group of a compact irreducible three-manifold $M$ with zero Euler characteristic. Then there is a short exact sequence $$1\to F\to G\to Q \to 1$$ where  $F<G$ is a finite group and $Q$ has a finite index subgroup isomorphic to the fundamental group of a compact three-manifold with zero Euler characteristic.
\end{theorem}

\newtheorem*{theorem:zero}{Theorem \ref{thm:zero exact seq}}

As we will see in \textsection \ref{splittings finite}, it is relatively easy, using the work of \cite{papasoglu:whyte:splitf}, to remove the assumption that $M$ is irreducible. 
Thus extending Theorem \ref{thm:zero exact seq} to manifolds with negative Euler characteristic completely settles the question. This is the purpose of the present paper which leads to the following statement:

\begin{theorem}[Quasi-isometric rigidity of three-manifold groups]\label{thm:main1} The class of fundamental groups of compact three-manifolds 
is quasi-isometrically rigid. More precisely,
 a finitely generated group quasi-isometric to  the  fundamental group of a compact three-manifold $M$ contains a finite
index subgroup isomorphic to the  fundamental group of a compact three-manifold $N$.
\end{theorem}

\newtheorem*{theorem:main1}{Theorem \ref{thm:main1}}

Notice that $\pi_1(M)$ and $\pi_1(N)$ need not be commensurable, i.e., may have no isomorphic finite index subgroups. 
For example, consider closed quotients of $\HH^2\times \E^1$ and $\widetilde{SL_2(\R)}$: since $\HH^2\times \E^1$ and $\widetilde{SL_2}(\R)$ are quasi-isometric, the fundamental groups of all those quotients are quasi-isometric by \v{S}varc--Milnor lemma, 
but $\HH^2\times \E^1$ and $\widetilde{SL_2(\R)}$ have no isomorphic lattices, see \cite[Lemma 6.3]{waldhausen:irreducible} and \cite[Theorem 5.2]{scott:geo3}. 
Non-geometric examples can also be produced, using \cite{leeb:nonpositive, kapovich:leeb:npc3man, kapovich:leeb:discretenpc}.
Nevertheless, it follows from the present work that, for a non-geometric irreducible manifold or a hyperbolic manifold with non-empty boundary, the fundamental groups of the pieces obtained after cutting the manifold along compressing discs, essential tori and annuli are well-defined up to commensurability, cf. Theorem \ref{thm:commensurability}.

Notice also that Theorem \ref{thm:main1} involves a slight upgrade of Theorem \ref{thm:zero exact seq} to go from a short exact sequence to a finite index subgroup. This comes from the following statement, interesting in its own right:

\begin{theorem}\label{thm:quotientman} Let $G$ be a finitely generated group and $p:G\to Q$ a surjective morphism with finite kernel. 
If $Q$ has a finite index subgroup isomorphic to the fundamental group of a compact $2$- or three-manifold $M$ then $G$ is commensurable to $Q$. \end{theorem}

Two groups $G$ and $Q$ are {\it commensurable} if there are subgroups $G'<G$ and $Q'<Q$ of finite index such that $G'$ is isomorphic to $Q'$.

\newtheorem*{theorem:quotientman}{Theorem \ref{thm:quotientman}}

One of our main inputs deals with finitely generated Kleinian groups, i.e., discrete subgroups of $\PP SL_2(\C)$. We provide a new proof of the quasi-isometric rigidity of  convex-cocompact Kleinian groups \cite{ph:unifplanar,ph:abc-qrig} 
which holds for all (finitely generated) Kleinian groups, leading to:

\begin{theorem}[Quasi-isometric rigidity of Kleinian groups]\label{thm:main2} The class of  Kleinian groups is quasi-isometrically rigid. More precisely,
 a finitely generated group quasi-isometric to a  Kleinian group  contains a finite
index subgroup isomorphic to a  (possibly different) geometrically finite  Kleinian group.
\end{theorem}

\newtheorem*{theorem:main2}{Theorem \ref{thm:main2}}

To add a little perspective to this introduction let us remark that quasi-isometric groups may also be fairly different from one another. For instance, simplicity, linearity and residual finiteness \cite{burger:mozes:finitely:products}, Kazhdan's property (T) (\cite{gersten:combings} and  \cite[Theorem 19.76]{drutu:kapovich:ggt})  and the Haagerup property \cite{carette:haagerup}  are not virtually invariant under quasi-isometries. We refer the interested reader to  \cite{burger:mozes:lattices:product},  \cite[\textsection IV.50]{delaharpe:topics} and \cite[Chap. 25]{drutu:kapovich:ggt}
for more properties that are known to be invariant or not under quasi-isometries.

When considering three-manifold groups from the point of view of geometric group theory it is also natural to wonder about their classification up to quasi-isometry. This question has been the subject of extensive work of Behrstock and Neumann \cite{behrstock:neumann:qigraph, behrstock:neumann:qinongeom} who give a nearly complete answer.

{\noindent\bf Outline of the paper.} \label{outline}
The proofs of Theorems \ref{thm:main1} and \ref{thm:main2} can be decomposed into three steps which are similar for both proofs. Let $Q=\pi_1(M)$ be either the fundamental group of a compact three-manifold (for Theorem \ref{thm:main1}) or a compact hyperbolic three-manifold (for Theorem \ref{thm:main2}) and let $G$ be a finitely generated group quasi-isometric to $Q$. 

In the first step, we simultaneously split $G$ and $Q$ to get graphs of groups $$\GGG=(\G_G,\{G_v\},\{G_e\}, G_e\hookrightarrow G_{t(e)}) \hbox{ and } \QQQ=(\G_Q,\{Q_v\},\{Q_e\}, Q_e\hookrightarrow Q_{t(e)})$$ (see \textsection \ref{scn:graph of group} for definitions and notations) 
with quasi-isometric vertex groups. The splittings of $Q$ come from topological splittings of $M$ and the vertex groups $Q_v$ are fundamental groups 
of hyperbolic three-manifolds and of compact three-manifolds with zero Euler characteristic in the proof of Theorem \ref{thm:main1},
while in the proof of Theorem \ref{thm:main2}, they are fundamental groups of pared acylindrical hyperbolic three-manifolds and pared $I$-bundles. 

The second step consists in showing that the classes to which the vertex groups $Q_v$ belong are quasi-isometrically rigid. For Theorem \ref{thm:main1},  it is given by Theorems \ref{thm:main2} and \ref{thm:zero exact seq}. For Theorem \ref{thm:main2}, it is the quasi-isometric rigidity of pared acylindrical Kleinian groups (Theorem \ref{thm:paredmain}) and pared $I$-bundles (Lemma \ref{lma:ibundle}). Applying these results to the vertex groups $G_v$ we get finite index subgroups $G'_v<G_v$ which are fundamental groups of compact three-manifolds $M_v$. 

In the third and last step we find a finite index subgroup $G'<G$ whose intersection with each vertex group $G_v$ is a finite index subgroup of $G'_v$ and deduce that $G'$ is the fundamental group of a compact three-manifold obtained by gluing together finite covers of the manifolds $M_v$. 

The arguments used in these steps are independent and for a better exposition we will study them in a different order. We will give more insights on these steps while we detail the plan of the paper.

First, in Section \ref{sec:qi 3mfds} we review results on quasi-isometric rigidity that are related to our topic. 
In particular, we recount the results that lead to Theorem \ref{thm:zero exact seq}. In Section \ref{3 mfds} we provide 
background on three-manifolds and Kleinian groups. In particular, we show that a finitely generated
Kleinian group is always isomorphic to a so-called geometrically finite Kleinian group with minimal parabolics, cf. Proposition \ref{prop:minimal}.
In Section \ref{sec:hyperbolicity} we introduce word hyperbolic and relatively hyperbolic groups and establish some facts that will be used later on. 
In Section \ref{sec:cansplit} we establish the first step described above, i.e., we study characteristic splittings of $Q$ and $G$ and their quasi-isometric invariance. 
In Section \ref{splittings finite}, we introduce a maximal splitting along spheres and discs and its analog for groups \cite{dunwoody:accessibility}
which are used in the proofs of both Theorems \ref{thm:main1} and \ref{thm:main2}. 
Its quasi-isometric invariance has been established in \cite{papasoglu:whyte:splitf}. 
In Section \ref{scn:torus decomposition}, we describe the characteristic torus decomposition and the induced splittings which are used in the proof of Theorem \ref{thm:main1}. 
The analog for groups and its quasi-isometric invariance are built on the work of Kapovich--Leeb \cite{kapovich:leeb:qi3man}. 
In Section \ref{scn:annulus decomposition}, we introduce the characteristic annulus decomposition which is used in the proof of Theorem \ref{thm:main2}. 
Its analog for groups is defined using the work of \cite{papasoglu:swenson:continua} on trees associated to the cut points and cut pairs of a continuum. 
The quasi-isometric invariance of the splitting thus produced is proved using the relation between the limit set of a Kleinian group and 
the characteristic submanifold of the underlying manifold \cite{walsh:bumping}, see also \cite[\textsection{2.3}]{lecuire:plissage}. 
In Section \ref{sec:prescription} we set up the third step by building finite index subgroups of fundamental groups of graph of groups 
with prescribed intersections with edge and vertex groups. The main tool used here is Wise's virtually special quotient theorem \cite[Theorem 15.6]{wise:qcvxh}, see also \cite[Theorem 12.1]{wise:qcvxh}. 
In the last section, we conclude the proof of our main theorems. First we show Theorem \ref{thm:quotientman} in Section \ref{sec:residually finite} 
using an induction argument and hierarchies in groups quasi-isometric to three-manifold groups. In Section \ref{sec:qi acylindrical}, 
we establish the second step of the proof of Theorem \ref{thm:main2} by proving the quasi-isometric rigidity of pared $I$-bundles and pared acylindrical Kleinian groups. Lastly, in Sections \ref{sec:qi kleinian} and \ref{sec:qi general}, we proceed with the proof of Theorems \ref{thm:main2} and \ref{thm:main1} as explained above.

 {\noindent\bf Acknowledgements.--} 
 We are grateful to Misha Kapovich for having brought this question to our attention. 
 We also feel particularly indebted to Dani Wise for his help regarding the construction of finite index subgroups
 and for his enthusiastic support on this project. It is our pleasure to thank Daniel Groves and Jason Manning for fruitful discussions. We also thank Jack Button for pointing out a mistake in the introduction of a previous version and Chris Hruska for noticing that arguments from the proof of Theorem \ref{thm:quotientman} can be used to remove the faithfulness assumption in the Bieberbach theorem (Theorem \ref{thm:bieberbach}). We are grateful to the referees for their thorough reading of this work and their many suggestions that have improved the paper.
This work was partially supported by the
 ANR project ``GDSous/GSG'' no. 12-BS01-0003-01.
 	
\section{Topology and geometry of three-manifolds}

\subsection{Quasi-isometric rigidity and geometric manifolds}   \label{sec:qi 3mfds}

To complete our introduction, we will now give more details on the quasi-isometric rigidity results that were mentioned earlier.

The work of Stallings on ends of groups \cite{stallings:yale} is a natural starting point. It leads to the quasi-isometric rigidity of virtually free groups, see  \cite[Theorem 20.45]{drutu:kapovich:ggt}. A group $G$ is said to {\it virtually} have a property if a finite index subgroup $H$ of $G$ has the said property.

\begin{theorem}		\label{thm:free}

The classes of virtually cyclic and virtually non-Abelian free groups are quasi-isometrically rigid.
\end{theorem}

Another class of groups for which the quasi-isometric rigidity has been established is the class of virtually nilpotent groups. The next result follows from Gromov's polynomial growth theorem \cite{gromov:expanding}:

\begin{theorem}[Groups of polynomial growth]\label{thm:nilpotent}
The class of virtually nilpotent groups is quasi-isometrically rigid.
\end{theorem}

Combining
Theorem \ref{thm:nilpotent} with Bass-Guivarc'h formula for the polynomial growth of nilpotent groups \cite{bass:degree, guivarch:croissance} 
and the work of Pansu on their asymptotic cones \cite{pansu:croissance}, one gets \cite[Theorem 16.26]{drutu:kapovich:ggt}:

\begin{theorem}		\label{thm:abelian}
The class of virtually Abelian groups is quasi-isometrically rigid, with one quasi-isometry class for each rank.
\end{theorem}

For a proof avoiding the classification of groups of polynomial growth, see \cite{cornulier:tessera:valette:gafa07}.

\subsubsection{Surface groups}

Next, we explain the quasi-isometric rigidity of surface groups.

\begin{theorem} \label{thm:qi surfaces ses}
Let $G$ be a group quasi-isometric to the fundamental group of a compact surface $S$ (with or without boundary).  Then there is a short exact sequence $$1\to F\to G\to Q \to 1$$ where  $F<G$ is a finite group and $Q$ has a finite index subgroup isomorphic to the fundamental group of a compact surface.
\end{theorem}

To make good use of the \v{S}varc--Milnor lemma, we want to put a convenient metric on a compact manifold. In dimension $2$, the Poincar\'e-Koebe uniformization theorem provides us with a spherical, Euclidean or hyperbolic metric for any compact surface. For better consistency with the $3$ dimensional case, let us put that statement in the perspective of Thurston's model geometry.

A {\it model geometry} $(G,X)$ is a manifold $X$ together with a Lie group $G$ of diffeomorphisms of $X$, such that:
\begin{enumerate}[(a)]
\item $X$ is connected and simply connected;
\item $G$ acts transitively on $X$, with compact point stabilizers;
\item $G$ is not contained in any larger group of diffeomorphisms of $X$ with compact stabilizers of points; and
\item there exists at least one closed manifold $M$ {\it modelled on $(G,X)$, meaning that there is a diffeomorphism 
from $M$ to $X/\Gamma$, 
where $\Gamma$ is 
a discrete subgroup of $G$ acting freely on $X$}.
\end{enumerate}

More generally, if $M$ is a compact manifold, we say that $M$ is {\it modelled on $(G,X)$} when there is a diffeomorphism from $int(M)$ to $X/\G$. Then $X/\G$ defines a {\it geometric structure} and we say that a manifold is  {\it geometric} if it has a geometric structure.

In dimension two there are three model geometries \cite[Theorem 3.8.2]{thurston:geometry}: $\SS^2$, $\E^2$ and $\HH^2$ together with their groups of isometries. Thus, the Poincar\'e-Koebe uniformization theorem says that all compact surfaces are geometric. Now the proof of Theorem \ref{thm:qi surfaces ses} boils down to the rigidity of discrete subgroups of isometries of $\SS^2$, $\E^2$ and $\HH^2$.

\demode{Theorem \ref{thm:qi surfaces ses}}
As we have already explained, $\pi_1(S)$ is isomorphic to a discrete subgroup $\Gamma$ of $Isom(X)$ with $X=\SS^2$, $\E^2$ or $\HH^2$.

If $X=\SS^2$, $G$ is finite and there is nothing to prove.

If $X=\E^2$, $\pi_1(S)$ is virtually Abelian of rank $1$ or $2$ depending on whether $S$ has a boundary or not and 
the conclusion follows from Theorem \ref{thm:abelian}; an alternate proof can be found in \cite{mess:cvaction}.

If $X=\HH^2$, $\Gamma$ can be chosen to be a lattice of $\PP SL_2(\R)$ ---a {\em Fuchsian group} of finite area. 
If $\Gamma$ is non uniform (equivalently if $\partial S\neq\emptyset$), it is a free group and the conclusion follows from Theorem \ref{thm:free}. Thus we are only left with the case where $G$ is quasi-isometric to a uniform lattice $K<\PP SL_2(\R)$. 
Then $K$ is word hyperbolic (see Definition \ref{dfn:hyperbolic group}) with boundary homeomorphic to $S^1$, so $G$  is as well. Therefore, $G$ admits a uniform action on $S^1$ and it follows from \cite{casson:jungreis} or \cite{gabai:S1} that this action is conjugate to that of a cocompact Fuchsian group.
Since $G$ is a convergence group, the kernel of the action is finite.
\endp

Combining Theorems \ref{thm:qi surfaces ses} and \ref{thm:quotientman}, we get:

\begin{theorem} \label{thm:qi surfaces}
Let $G$ be a group quasi-isometric to the fundamental group of a compact surface $S$ (with or without boundary).  Then $G$ has a finite index subgroup isomorphic to the fundamental group of a compact surface.
\end{theorem}

\subsubsection{Three-manifold groups}
For fundamental groups of three-manifolds, we summarized in the introduction the state of the art with the following statement:

\begin{theorem:zero}
Let $G$ be a group quasi-isometric to the fundamental group of a compact irreducible three-manifold $M$ with zero Euler characteristic. Then there is a short exact sequence $$1\to F\to G\to Q \to 1$$ where  $F<G$ is a finite group and $Q$ has a finite index subgroup isomorphic to the fundamental group of a compact three-manifold with zero Euler characteristic.
\end{theorem:zero}

A three-manifold is {\it irreducible} if every embedded sphere bounds a ball  (the rest of the terminology used below is given in Section \ref{3 mfds}).
As in the surface case, a first step in the proof consists in equipping three-manifolds with convenient metrics.

Thurston has shown \cite[Theorem 3.8.4]{thurston:geometry} that there are eight three-dimensional model geometries $(G,X)$ which are $\SS^3$, $\E^3$, $\HH^3$, $\SS^2\times\E^1$, $\HH^2\times\E^1$, $Nil$, $\widetilde{SL(2,\R)}$, $Sol$ together with their groups of isometries.

In general, three-manifolds are not geometric but the geometrization theorem (proved for Haken manifolds and stated as a conjecture in general by Thurston \cite{thurston:bams} and proved by Perel'man in general \cite{bessiere:et:all, kleiner:lott, morgan:tian}) asserts that they can be decomposed into geometric pieces.

\begin{theorem}[Geometrization]		\label{thm:geometrization}
Every oriented irreducible compact three-manifold can be cut along tori, so that the interior of each of the resulting manifolds has a geometric structure.
\end{theorem}

{\em A geometric decomposition} of a three-manifold $M$ is a collection of essential tori $T$ such that each component of $M\setminus T$ has a geometric structure. The geometrization theorem precisely states that such a geometric decomposition always exists. We say that a geometric decomposition $T$ is {\em minimal} if for any component $T_1$ of $T$, $T\setminus T_1$ is not a geometric decomposition.

Let us now explain the main steps of the proof of Theorem \ref{thm:zero exact seq}, starting with geometric three-manifolds.

\begin{theorem} \label{thm:list geom}
Let $G$ be a group quasi-isometric to the fundamental group of a geometric three-manifold $M$ with non-negative Euler characteristic. Then there is a short exact sequence $$1\to F\to G\to Q \to 1$$ where  $F<G$ is a finite group and $Q$ has a finite index subgroup isomorphic to the fundamental group of a geometric three-manifold with zero Euler characteristic.
\end{theorem}

\Proof
We have already seen that $M$ is modelled on one of the homogeneous spaces $X=\SS^3, \SS^2\times\E^1,\E^3, \HH^2\times\E^1,\HH^3, Nil, Sol, \widetilde{SL_2(\R)}$.

When $\pi_1(M)$ is finite, in particular when $X=\SS^3$ or when $M$ has positive Euler characteristic, we take the trivial group considered as $\pi_1(\SS^3)$ for $Q$. We now assume that $\pi_1(M)$ is infinite so that its Euler characteristic is zero.

When $X=\SS^2\times\E^1$, then $\pi_1(M)$ is virtually cyclic and the conclusion follows from Theorem \ref{thm:free}. 

When $X=\E^3$, it is a special case of Theorem \ref{thm:abelian}.

When $X=Nil$, it follows from Theorem \ref{thm:nilpotent} and results of Mal'cev, Guivarc'h \cite{guivarch:croissance} and Jenkins \cite{jenkins:growth}, cf. \cite[Theorem 1.7]{frigerio:qirighandbook}.

When $X=\HH^3$, since $M$ has zero Euler characteristic, its hyperbolic structure has finite volume (see \cite[2.3]{thurston:bams} for example). The conclusion
 is due to Sullivan \cite{DS3} and Cannon-Cooper \cite{cannon:cooper} when $\partial M=\emptyset$ and Schwartz \cite{schwartz:qirank1} when $\partial M\neq\emptyset$.

When $X=Sol$, it has been proved by Eskin, Fisher and Whyte \cite{eskin:fisher:whyte:qisol}. 

Finally, when $X=\HH^2\times\E^1$ and $X=\widetilde{SL_2(\R)}$, it is due to Rieffel \cite{rieffel:qiH2xR}, see also \cite[\textsection 5.2]{kapovich:leeb:qi3man}. Note that these two geometries are quasi-isometric. \endp

Notice that geometric three-manifolds with negative Euler characteristic are hyperbolic, i.e., modeled on $\HH^3$.

From a strong invariance of the geometric decomposition under quasi-isometry, Kapovich--Leeb \cite{kapovich:leeb:qi3man} deduce the quasi-isometric rigidity of fundamental groups of irreducible non-geometric three-manifolds with zero Euler characteristic. We will give more details in \textsection \ref{scn:torus decomposition} since we will use this invariance to prove Theorem \ref{thm:main1}.

\begin{theorem}[Kapovich--Leeb] \label{thm:non geom}
Let $G$ be a group quasi-isometric to the fundamental group of an irreducible non-geometric Haken compact three-manifold $M$ with zero Euler characteristic. Then there is a short exact sequence $$1\to F\to G\to Q \to 1$$ where  $F<G$ is a finite group and $Q$ has a finite index subgroup isomorphic to the fundamental group of an irreducible non-geometric compact three-manifold $M'$ with zero Euler characteristic.
\end{theorem}

By the geometrization theorem (Theorem \ref{thm:geometrization}), any non-geometric compact three-manifold is Haken (see the next section for a definition). Thus we get Theorem \ref{thm:zero exact seq} simply by combining Theorems \ref{thm:list geom} and \ref{thm:non geom}.

\subsection{Three-manifolds {and groups}}    \label{3 mfds}
In this section and the next one, we will review some classical definitions and results that are used in this paper. We start with three-manifold topology, basic references include \cite{jaco:three_manifold_topology,hempel:three_manifolds}.

\subsubsection{Three-manifold topology} \label{3 mfds 1}

Let $M$ be a compact orientable three-manifold with boundary or not. An embedded surface $(S,\partial S)\to (M,\partial M)$ is {\it incompressible} if the inclusion $i:S\to M$ gives rise to an injective morphism $i_*:\pi_1(S,x)\to\pi_1(M,x)$. The {\it double} of  a manifold $M$ with boundary is the union of $M$ and of a copy of itself glued along its boundary. An embedded surface $S$ in $M$ is {\it boundary incompressible} if its double is incompressible in the double of $M$. A surface $S$ in $M$ is {\it non-peripheral} if it is {\it properly embedded}, i.e., $S\cap\partial M=\partial S$, and if the inclusion $i:S\to M$
is not homotopic relative to $\partial S$ to a map $f:S\to M$ such that $f(S)\subset\partial M$. A surface $S$ is {\it essential}
if it is properly embedded, two-sided, incompressible, boundary incompressible, non-peripheral and does not bound a $3$-ball. An essential disc is also called a {\it compression disc}. 

The manifold $M$ is {\it irreducible} if it contains no essential sphere, equivalently every embedded sphere bounds a ball. We say that $M$ is {\it boundary irreducible} or {\it $\partial$-irreducible} 
if each component of $\partial M$ is incompressible. By results of Kneser and Stallings, $M$ is irreducible and boundary irreducible if and only if $\pi_1(M)$ is one-ended.

The irreducible manifold $M$ is {\it Haken} if it contains an essential surface. If $\partial M\neq\emptyset$ then $M$ is necessarily  Haken, see \cite[Chap. 13]{hempel:three_manifolds}.

An irreducible three-manifold $M$ is {\it atoroidal} if every subgroup of $\pi_1(M)$ isomorphic to $\Z\oplus\Z$ is conjugate to a subgroup of the fundamental group of a boundary component. As we will see in the next section this property characterizes hyperbolic manifolds. An {\it acylindrical compact manifold} is atoroidal,
has incompressible boundary and no essential annuli.

A {\it Seifert manifold} is a compact three-manifold which admits a foliation by circles. Even though this is not the classical definition, it is equivalent to it \cite{dbaepstein:seifert}. A {\it graph manifold} is a compact irreducible $\partial$-irreducible three-manifold with no atoroidal pieces in its torus decomposition, see \textsection \ref{scn:torus decomposition}.

A compact {\it pared manifold} $(M,P)$ is given by a compact irreducible three-manifold $M$ non-homeomorphic to a solid torus nor a thickened torus together with a paring $P\subset\partial M$ which is a finite
collection of pairwise disjoint incompressible annuli and tori  satisfying:
\begin{enumerate}[-]
\item  every Abelian, non cyclic subgroup of $\pi_1(M)$ is conjugate to a subgroup of the fundamental group of a component of $P$ and
\item   any incompressible cylinder $(C,\partial C)\subset (M,P)$ can be homotoped relatively to its boundary into $P$.
\end{enumerate}

Notice that by definition only atoroidal manifolds can have a paring. On the other hand, as shown by the following classical result, an atoroidal three-manifold has a natural paring.

\newtheorem*{lemman}{Lemma}

\begin{lemma}  \label{lma:natural pairing}
Let $M$ be an irreducible orientable compact atoroidal three-manifold and denote by $P$ the union of the incompressible tori in $\partial M$. Then $(M,P)$ is a pared manifold unless $M$ is 
a solid torus, an $I$-bundle over a torus or Klein bottle.
\end{lemma}

\Proof
Since $M$ is atoroidal, every Abelian, non cyclic subgroup of $\pi_1(M)$ is conjugate to a subgroup of the fundamental group of a component of $P$. It remains to check whether or not any incompressible cylinder $(C,\partial C)\subset (M,P)$ can be homotoped relatively to its boundary into $P$.
If $P=\emptyset$ there is nothing to prove. 
Let us now assume that every torus component of $\partial M$ is incompressible and that there is an incompressible cylinder $(C,\partial C)\subset (M,P)$ that can not be homotoped relative to its boundary into $P$.

We have two cases to consider. In the first case, $C$ joins two different components $T_1$ and $T_2$ of $P$, then a regular neighborhood ${\mathcal N}$ of $T_1\cup T_2\cup C$ is a product 
 $S^1\times S_{0,3}$ of a circle by a pair of pants. In the second case, the two boundary components of $C$ lie in the same component $T_1$ of $T$ and a regular neighborhoud ${\mathcal N}$ of $T_1\cup C$ is again a product $S^1\times S_{0,3}$. 

Let $T_1, T_2$ and $T_3$ be the components of $\partial{\mathcal N}$. By construction $T_1$ (and $T_2$ in the first case) is a component of $\partial M$.  
By the Loop theorem, if the inclusion $T_i\to M$ is not $\pi_1$-injective then $T_i$ bounds a solid torus for $i\in\{1,2,3\}$. Since $M$ is atoroidal, if $T_i$ is incompressible for $i=2$ or $3$, then  $T_i$ can be homotoped into a boundary component $T'_i$ of $M$. Then, by \cite[Theorem 10.2]{hempel:three_manifolds}, $T_i$ and $T'_i$ bound a thickened torus and $T_i$ is isotopic to $T_i'$. It follows that either $M$ is homeomorphic to ${\mathcal N}$ or $M$ is obtained by gluing a solid torus along $T_3$ and/or $T_2$. Since $C$ is incompressible, the product structure of ${\mathcal N}$ extends to foliations by circles. Thus we have established that $M$ is a Seifert fibered manifold, in particular $M$ has a finite cover $M'$ which is a circle bundle over a compact surface $F$ (see \cite[Theorem 12.2]{hempel:three_manifolds}) such that $M'$ and $F$ are orientable. Since $M$ is atoroidal, $M'$ is also atoroidal, hence $F$ has genus $0$ and at most two boundary components (see \cite[Example 1.20]{kapovich:book}). Now $M'$ is either a solid torus or a thickened torus. We have already ruled out the case of a solid torus so we are left with $M$ being an $I$-bundle over a torus or Klein bottle.
\endp

We say that $(M,P)$ is {\it boundary irreducible} if each component of $\partial M\setminus P$ is incompressible and {\it acylindrical} if there is no essential disc or cylinder in $M$ disjoint from $P$.

\subsubsection{Kleinian groups}	\label{sec:kleinian}

A {\it Kleinian group}  $K$ is a discrete subgroup of $\PP SL_2(\C)$ ---the group of orientation preserving isometries of the $3$-dimensional hyperbolic space $\HH^3$. 
An orientable compact three-manifold $M$ is {\it hyperbolizable} if its interior is homeomorphic to the quotient $\HH^3/K$ where $K$ is a torsion free Kleinian group. Such a manifold $M$ is irreducible and atoroidal and we say that $M$ is {\it uniformized} by $K$. Note that $K$ is isomorphic to the fundamental 
group of $M$, and that it is necessarily word hyperbolic if it contains no subgroup isomorphic to $\Z\oplus\Z$, see Section \ref{sec:wordhyperbolic}. On the other hand, the tameness theorem \cite{agol:tameness,calegary:gabai} asserts that, when $K$ is a finitely generated torsion free Kleinian group, $\HH^3/K$ is homeomorphic to the interior of a compact three-manifold $M_K$ (with fundamental group isomorphic to $K$) that we call 
the {\em Kleinian manifold of $K$}.

As Poincar\'e observed, we may identify the Riemann sphere with the boundary at infinity of  $\HH^3$ \cite{poincare:memoire:kleineen}. Then $K$ acts  on the Riemann
sphere $\cbar$ via M\"obius transformations. 
The latter action is usually not properly discontinuous: there is 
a canonical and invariant partition $$\cbar = \Omega_K \sqcup \La_K$$
where $\Omega_K$ denotes the  {\it domain of discontinuity}, also called {\it ordinary set}, which is the largest  open set of $\cbar$ on which $K$ acts properly discontinuously, 
and where $\La_K$ denotes the  {\it limit
set}, which is  the set of accumulation points in $\cbar$ of an (any) orbit $Kp\subset\HH^3$. When $\La_K$ has at least three points, implying that it is infinite (we say that $K$ is {\it non-elementary}), it is also the minimal non-empty $K$-invariant compact subset of $\cbar$.
The construction of the Kleinian manifold $M_K$ induces an embedding $(\HH^3\cup \Omega_K)/K\hookrightarrow M_K$ whose image is the complement of a subsurface of $\partial M_K$.

If $\La_K$ contains at least $2$ points, the group $K$ preserves the convex hull $\hbox{Hull}(\La_K)$ of its limit set in $\HH^3$. The group  $K$ is  {\it convex-cocompact} if its action is cocompact on $\hbox{Hull}(\La_K)$ and $K$ is {\it geometrically finite} if a regular neighborhood  of a fundamental domain for the action of $K$ on $\mathrm{Hull}(\La_K)$
has finite volume  or if the limit set consists of a singleton. 
Ahlfors \cite{ahlfors:fundamental} showed that if the limit set of a geometrically finite Kleinian group is not the whole Riemann sphere then it has measure 0 (this holds more generally for finitely generated Kleinian groups by \cite{canary:ends}, \cite{agol:tameness} and \cite{calegary:gabai}).   When $K$ is not Fuchsian, i.e., $\hbox{Hull}(\La_K)$ is not contained in a geodesic plane, 
and $K$ is torsion free, there is a homeomorphism between $(\HH^3\cup \Omega_K)/K$ and the convex core $\hbox{Hull}(\La_K)/K$ constructed using the closest point projection. As we have explained above $(\HH^3\cup \Omega_K)/K$ embeds in $M_K$, hence the {\em convex core} also embeds in the Kleinian manifold. When $K$ is geometrically finite, the image of this embedding is the complement of a paring $P$ of $M_K$ corresponding to the parabolics of $K$ (notice that this gives an alternate characterization of geometric finiteness). We say that a  torsion free geometrically finite Kleinian group $K$ {\em uniformizes} the pared manifold $(M,P)$ when $(\HH^3\cup \Omega_K)/K$ is homeomorphic to $M\setminus P$. 

Now it follows from the Loop theorem (see \cite[Theorem 4.2]{hempel:three_manifolds} for example) that $\partial M\setminus P$ is incompressible if and only if  each connected component of $\Omega_K$ is simply connected. This leads to:

\begin{fact}    \label{fact:lambda connected}
Let $K$ be a geometrically finite torsion free Kleinian group uniformizing the pared manifold $(M,P)$. Then $\partial M\setminus P$ is incompressible if and only if $\Lambda_K$ is connected.
\end{fact}

Thurston's uniformization theorem (extended to non Haken manifolds by Perel'man) gives a topological characterization of hyperbolic manifolds. We will use the following form, see \cite[Theorem B']{morgan:thurston}, and also \cite{otal:fibre, otal:haken, kapovich:book}:

\REFTHM{thm:thurston} Let $(M,P)$ be a Haken pared three-manifold. There is a geometrically finite, complete hyperbolic manifold $N$ whose convex core is homeomorphic to $M\setminus P$.
\ENDTHM

When $(M,P)$ is acylindrical, then a doubling argument shows that we can require the convex core to have totally geodesic boundary \cite[Thm 3]{thurston:bangor79}.

\REFTHM{thm:thurston acylindrical} Let $(M,P)$ be a Haken acylindrical pared three-manifold. There is a geometrically finite, complete hyperbolic manifold $N$ whose convex core is homeomorphic to $M\setminus P$ and has totally geodesic boundary.
\ENDTHM

We say that a Kleinian group is {\it minimally parabolic} if every parabolic subgroup is  
a rank $2$ Abelian subgroup. Combining Theorem \ref{thm:thurston} with Scott's core theorem we get that any finitely generated Kleinian group has a minimally parabolic version:

\begin{proposition}		\label{prop:minimal}
Any torsion-free finitely generated Kleinian group is isomorphic to a geometrically finite, minimally parabolic, Kleinian group.
\end{proposition}

\Proof
Let $K$ be a  torsion-free finitely generated Kleinian group. By Scott's core theorem \cite{scott:core}, 
the hyperbolic manifold $\HH^3/K$ contains a compact three-dimensional submanifold $C$ such that the inclusion is a homotopy equivalence. 
If $\partial C=\emptyset$, then $C=\HH^3/K$ and $\HH^3/K$ is compact. Thus, $K$ has no parabolic subgroup and there is nothing to prove. So let us assume that $\partial C\neq\emptyset$, in particular, $C$ is Haken.
A maximal Abelian non cyclic subgroup of $\pi_1(C)$ corresponds to a rank $2$ parabolic subgroup of $K$. By \cite[Theorem 2]{mccullough:submfds}  
we can change $C$ by a homotopy so that such a subgroup is the fundamental group of a component of $\partial C$. 
It follows that $C$ is atoroidal. Now if $P$ is the union of the tori in $\partial C$, 
then $(C,P)$ is a pared manifold and we conclude with Theorem \ref{thm:thurston}.
\endp

 \begin{rmk}    \label{torsion free etc}
For any finitely generated Kleinian group, Selberg's lemma ensures the existence  of a torsion-free subgroup of finite index. Combined with Proposition \ref{prop:minimal}, this ensures that a group $G$ that is quasi-isometric to a Kleinian group, is also quasi-isometric to a torsion free geometrically finite minimally parabolic Kleinian group.
\end{rmk}

\section{Hyperbolicity in the sense of Gromov}  \label{sec:hyperbolicity}

Hyperbolic spaces and groups were introduced by Gromov in \cite{gromov:hyperbolic} and have known many developments since. 
In this section, we will introduce hyperbolic spaces, hyperbolic and relatively hyperbolic groups and a few more objects associated to these groups. We will also establish various facts that will be used throughout the paper.
We start with hyperbolic spaces and word hyperbolic groups. Background on those includes \cite{gromov:hyperbolic,ghys:delaharpe:groupes,
kapovitch:benakli}.

\subsection{Hyperbolic spaces}
Let $X$ be a metric space. It is {\it geodesic} if any pair of points $\{x,y\}$ 
can be joined by a (geodesic) segment i.e, a map $\g:[0,d(x,y)]\to X$ such that $\g(0)=x$, $\g(d(x,y))=y$
and $d(\g(s),\g(t))=|t-s|$ for all $s,t\in [0,d(x,y)]$. The metric space $X$ is {\it proper} if closed balls
of finite radius are compact. 

A {\it triangle} $\De$ in a metric space $X$ is given by  three points $\{x,y,z\}$ and three (geodesic) segments (or sides)
joining them. 
Given a constant $\de \ge 0$, the triangle $\De$ is {$\de$-\it slim} if any side
of the triangle is contained in the $\de$-neighborhood of the two others. 

\begin{defn}[Gromov hyperbolic spaces] A geodesic metric space is (Gromov) hyperbolic if there exists $\de\ge 0$
such that every triangle is $\de$-{slim}.\end{defn}

Basic examples of hyperbolic spaces are the complete simply connected hyperbolic manifolds $\HH^n$,  $\R$-trees and their convex subsets. 
It follows from the shadowing lemma (see below) that, among geodesic metric spaces, hyperbolicity is invariant under quasi-isometry : 
if $X$, $Y$ are two quasi-isometric geodesic metric spaces, then $X$ is hyperbolic if and only if $Y$ is hyperbolic. A $(\lambda,c)$-{\it quasigeodesic} is the image of an interval {(that can be unbounded on both sides)} by a $(\lambda,c)$-quasi-isometric embedding.

\begin{lemma}[Shadowing lemma]		\label{lm:shadow}
Given $\de$, $\la$ and $c$, there is a constant $H=H(\de,\la,c)$ such that, for
any $(\la,c)$-quasigeodesic $q$ in a proper geodesic $\de$-hyperbolic metric space $X$, there is a geodesic $\g$
at Hausdorff distance at most $H$ from $q$.
\end{lemma}

\subsubsection{Boundaries of hyperbolic spaces}

A  proper (geodesic) hyperbolic space $X$ admits a canonical compactification $X\sqcup \partial X$ at 
infinity. This compactification can be defined by looking at the set of geodesic rays, i.e., isometric embeddings $r:\R_+\to X$,
up to bounded Hausdorff distance. The topology is induced by the uniform convergence on compact subsets {of $\R_+$: fix a base point $w\in X$; 
for each point $x\in X$, choose arbitrarily a geodesic segment $s_x:[0,d(w,x)]\to X$ between $w$ and $x$ so that $s_x(0)=w$ and extend it
to $r_x:\R_+\to X$ by setting $r_x(t)=x$ if $t\ge d(w,x)$. A sequence $(x_n)_n$ tends to a point at infinity if $d(w,x_n)\longrightarrow\infty$ and all the limits of {converging (uniformly on compact sets) subsequences of }$(r_{x_n})_n$ {have their images} at bounded 
Hausdorff distance to each other. This is independent from the choice of the base point}.
The boundary can be endowed with a family of {\it visual distances} $d_v$  compatible with its topology, i.e., {for any base point $w\in X$, there exists a constant $C>0$ }which satisfies
$$(1/C)  e^{-\ep d(w , (a,b))} \le d_v(a,b) \le C e^{-\ep d(w , (a,b))}$$ 
where  $\ep >0$ is a {\it visual parameter}
chosen small enough, and $(a,b)$ is any geodesic asymptotic to rays defining $a$ and $b$.

If $\Phi:X\to Y$ is a quasi-isometry between two proper hyperbolic spaces, then the shadowing lemma implies that 
$\Phi$ induces a homeomorphism
$\phi:\partial X\to\partial Y$.  Actually, the trace map at infinity of a quasi-isometry is 
{\it  quasi-M\"obius} \cite{vam},  i.e., there exists a \homeo $\te:\R_+\to\R_+$
such that, for any distinct points $x_1,x_2,x_3,x_4\in \partial X$, 
$$[\phi(x_1):\phi(x_2):\phi(x_3):\phi(x_4)]\le \te([x_1:x_2:x_3:x_4])$$
where{ $$[x_1:x_2:x_3:x_4]= \frac{d_v(x_1,x_2)\cdot d_v(x_3,x_4)}{d_v(x_1,x_3)\cdot d_v(x_2,x_4)}\,.$$}
Quasi-M\"obius maps are stable under composition.

Quasi-isometries provide natural examples of quasi-M\"obius maps, cf. \cite[Prop.\,4.6]{paulin:determined}:
\begin{theorem}[Efremovi{\v{c}}--Tikhomirova, \cite{efremovic:tihomirova}]\label{thm:qisomqm} A 
$(\la,c)$-quasi-isometry between proper hyperbolic spaces extends as a $\te$-quasi-M\"obius map between their boundaries, where
$\te$ only depends on $\la, c$, the hyperbolicity constants and the visual parameters. \end{theorem}

\subsubsection{Groups of isometries {and convergence actions}}    \label{secn:isometries}

{Let us consider a group $G$ acting  by homeomorphisms on a compact metrisable space $Z$ with at least $3$ points. We say that this action is a {\em convergence group action} if its diagonal action on the set of distinct triples is properly discontinuous \cite{gehring:martin:qcgroupsi}.}
The {\it limit set}  {$\La_G$} is by definition {the set of accumulation points of an (any) infinite orbit {$Gp\subset Z$}.
One may consult \cite{bowditch:convergence_groups}
for basic properties of convergence groups.

Let {us now consider} a group $G$ of isometries of a proper hyperbolic space $X$. It follows from 
Theorem \ref{thm:qisomqm} that the action of  $G$ extends to an action on $\partial X$ by
uniform quasi-M\"obius mappings.  

If the action of $G$ is furthermore properly discontinuous on $X$, then the action of $G$ on {$X\cup \partial X$ and on} $\partial X$ is {a convergence action} \cite{freden:hypcv1, tukia:convergence_groups}.

\subsection{Word hyperbolic groups}\label{sec:wordhyperbolic}
A properly discontinuous group action by isometries on a proper geodesic
metric space is {\it geometric} if it is {\it cocompact}, i.e., if the quotient is compact. 

\begin{defn}[Word hyperbolic groups]    \label{dfn:hyperbolic group} A group $G$ is word hyperbolic, or just hyperbolic for simplicity, 
if it admits a geometric action on a proper geodesic hyperbolic
metric space.
\end{defn}
We note that since the hyperbolicity of a proper geodesic hyperbolic
metric space is a quasi-isometric invariant property, the hyperbolicity of a group does not depend on the space it is acting upon.
In particular, by the \v{S}varc--Milnor lemma, $G$ is finitely generated and its hyperbolicity can also be read from any of its locally finite Cayley graphs.
Moreover, this implies the quasi-isometric rigidity of the class of word hyperbolic groups.

Fundamental groups of closed hyperbolic manifolds are obvious examples. Convex cocompact Kleinian groups are also hyperbolic by definition: 
a convex-compact Kleinian group $G$ has a properly discontinuous and cocompact action on the convex hull $\hbox{Hull}(\La_G)$ of its limit set 
which is a hyperbolic space in the sense of Gromov.

The definition and the \v{S}varc--Milnor lemma also imply that a word hyperbolic group $G$ admits a topological boundary $\partial G$ defined
by considering the boundary of any proper geodesic metric space on which $G$ acts geometrically.
In the case of a convex-cocompact Kleinian group $K$, a model for the boundary $\partial K$ is given by its limit set $\La_K$.

\medskip

The action of a hyperbolic group  on its boundary  is a {\em uniform convergence action}, i.e.,  
its diagonal action on the set of distinct triples is not only properly discontinuous but also cocompact, cf. \cite{bowditch:convergence_groups}.

 These properties characterize word hyperbolic groups and their boundaries:

\begin{theorem}[Bowditch \cite{bowditch:characterization}]\label{thm:def:hyp} 
Let $G$ be a convergence group acting on a perfect metrizable space $X$.
The action of $G$ is uniform on its limit set $\La_G$ if and only if $G$ is word hyperbolic and if,  furthermore, there exists an equivariant homeomorphism
between $\La_G$ and the boundary at infinity $\partial G$ of $G$. 
\end{theorem}

A general principle asserts that a word hyperbolic group is determined  by its boundary. More precisely, Paulin proved that the quasi-isometry class of a word hyperbolic group
is determined by its boundary equipped with
its quasiconformal structure \cite{paulin:determined}. 
\begin{theorem}[Paulin \cite{paulin:determined}]\label{thm:paulin}  Two non-elementary word hyperbolic groups are quasi-isometric if and only if there is a quasi-M\"obius homeomorphism between their boundaries. \end{theorem}

\subsection{Relative hyperbolicity} \label{sec:relhyp}
The idea behind relatively hyperbolic groups which will be defined next is that some metric spaces such as Cayley graphs are hyperbolic away from some ``codimension 1'' subspaces or subgroups.
Background on relatively hyperbolic groups includes \cite{gromov:hyperbolic,bowditch:rhg,rhuska:qcvxrhg} and the references therein. 
Let us first recall the definition of the {\it Busemann function} $\be_p(x,y)$ centered at a point $p\in \partial X$ at infinity 
 of a hyperbolic space $X$ measured at
two points $x,y\in X$:
$$\be_p(x,y)=\sup \left\{\lim_{t\to\infty} [ d(x,\g(t)) -t]:\g \hbox{ geodesic ray asymptotic to }p \hbox{ such that }\g(0)=y\right\}\,. $$

Given $p\in \partial X$, a base point $w\in X$ and $r\in\R$, the {\it horoball centered at $p$ of level $r$} is defined as
$$H_w(p,r)=\{x\in X: {\be_p(x,w)}\le r\}\,.$$

Let $G$ be a group and  let $\PP$ be a collection of subgroups.
The pair  $(G,\PP)$ is {\it relatively hyperbolic} if there exists a hyperbolic  proper geodesic metric space $X$
on which $G$ acts properly discontinuously by isometries and if there is a {$G$-invariant collection $\HHH$ of  pairwise disjoint horoballs}  
with the following properties:
\ben
\item any $P\in\PP$ is the stabilizer of the center (at infinity) of a horoball in $\HHH$;
\item the stabilizer of any center of a horoball  of $\HHH$ is conjugate to a subgroup from $\PP$;
\item the action of $G$ on $X\setminus Y$ is cocompact,
where we let $Y$ denote the union of the horoballs in $\HHH$.
 \een
 
Such an action of $G$ is called  {\it cusp uniform}. Note that we may assume that no two elements of $\PP$ are conjugate{, so that $\PP$ must be finite}. 
The subgroups in $\PP$ are called {\em peripheral subgroups}. {Note also that, in this case, the limit set of the action is the whole boundary $\partial X$.} 

\newtheorem*{rmkn}{Remark}

\begin{rmk} If $K$ is a geometrically finite Kleinian group, then $K$ is hyperbolic relative to its maximal parabolic subgroups and its action on $\mbox{\rm Hull}(\La_K)$ is cusp-uniform.\end{rmk}

\subsubsection{Horoballing, cusped spaces and Bowditch boundaries}
In this section, we introduce an explicit cusp uniform action for a relatively hyperbolic group. We first concentrate on the construction of horoballs resting on peripheral subgroups.
Let $(P,d)$ be a graph endowed with the path metric  that makes each edge isometric to $[0,1]$. Define its {\it horoballing space} $H_P$ following Groves and Manning \cite{groves:manning:dehn}
as the graph modelled on $P\times \N$ with additional edges
$$\begin{pmatrix} x\\ m\end{pmatrix} \sim \begin{pmatrix} y\\ m\end{pmatrix}\quad\hbox{if} \quad d(x,y)\le 2^m$$
and 
$$\begin{pmatrix} x\\ m\end{pmatrix}\sim \begin{pmatrix} x\\ m+1\end{pmatrix}$$
This new graph $H_P$ is also endowed with the path metric so that each edge is isometric to $[0,1]$. 
By \cite[Theorem 3.8]{groves:manning:dehn}, $H_P$ is hyperbolic.

Let $G$ be a finitely generated group, let $\PP=\{P_1,...,P_n\}$ be a (finite) family of finitely generated subgroups of $G$, and let $S$ be a  finite 
generating set for $G$ so that $P_i\cap S$ generates $P_i$ for each $i\in\{1,...,n\}$ , and denote by $\cay(G,S)$ the Cayley graph of $(G,S)$. 
For each $i\in\{1,...,n\}$, let $T_i$ be a left transversal for $P_i$, i.e., a collection of representatives for left cosets of $P_i$ {in} $G$ which contains exactly one element of each left coset. For each $i$, and each $t\in T_i$, 
let $\cay_{i,t}$ be the full subgraph of the Cayley graph $\cay(G,S)$ which contains $tP_i$. Each $\cay_{i,t}$ is isomorphic to 
the Cayley graph of $P_i$ with respect to the generators $P_i\cap S$. We define the {\it cusped space}
$$\cus(G,\PP,S)=\cay(G,S)\cup(\{\cup H_{\cay_{i,t}}|1\leq i\leq n, t\in T_i\}),$$ 
where the graphs $\cay_{i,t}\subset \cay(G,S)$ and $\cay_{i,t}=\cay_{i,t}\times\{0\}\subset H_{\cay_{i,t}}$ are identified in the obvious way. 
When the choice of the generating set does not matter, we simplify the notations to $\cay(G)$ and $\cus(G,\PP)$.

The cusped space provides a canonical way to construct a hyperbolic space on which a relatively hyperbolic group has a cusp uniform action. We follow Groves and Manning \cite{groves:manning:dehn} in using this space to give a different  characterization of relatively hyperbolic groups.

\begin{theorem}[Groves and Manning \cite{groves:manning:dehn}]
Let $G$ be a finitely generated group, let $\PP=\{P_1,...,P_n\}$ be a (finite) family of finitely generated subgroups of $G$. The pair $(G,\PP)$ is relatively hyperbolic if there is a generating set for $G$ so that $P_i\cap S$ generates $P_i$ for each $i\in\{1,...,n\}$ and that $\cus(G,\PP,S)$ is hyperbolic.
\end{theorem}
The horoballings {$H_{\cay_{i,t}}  \subset \cus(G,\PP,S)$, $t\in G$, $P_i\in P$} are horoballs and  the action of $G$ on $\cus(G,\PP,S)$ is cusp uniform.
The cusped space provides us a way to define a boundary for relatively hyperbolic groups. 

\begin{defn}[Bowditch boundary] Given a finitely generated hyperbolic group $G$ relative to a finite family of finitely generated
subgroups $\PP$, the Bowditch boundary $\bop G$ of $(G,\PP)$ is defined as the boundary $\partial \cus(G,\PP)$ of the cusped space.\end{defn}

Corollary \ref{cor:bowbdry} below shows that the boundary is well defined, as a metric space, up to a quasi-M\"obius change of metrics.
Note that the topology of the boundary is well-defined according to \cite[Theorem 9.4]{bowditch:rhg}. Reference to a generating set is thus unnecessary.

In the next sections we will discuss quasi-isometries between relatively hyperbolic groups. In particular we will show that the previous definition does not depend on the choice of the generating set $S$.

\subsubsection{Geometrically finite actions} \label{sec: geom fin actions}

Let $G$ be a convergence group acting on some metrizable compact space $Z$. An element $g\in G$ is {\it loxodromic} if it has infinite order and fixes exactly two points of $Z$. 
A subgroup $H<G$ 
is {\it parabolic} if it is infinite and has a unique fixed point. We refer to it as a {\it parabolic point}. {A parabolic subgroup}
contains no loxodromics \cite{bowditch:convergence_groups}.
The stabilizer  of a parabolic point is necessarily a parabolic group. There
is thus a natural bijective correspondence between parabolic points in
$Z$ and maximal parabolic subgroups of $G$. We say that a parabolic
group, $H$, with fixed point $p$, is {\it bounded} if the quotient $(\La_G \setminus \{p\})/H$
is compact. (It is necessarily Hausdorff.) We say that $p$ is a {\it bounded
parabolic point } if its stabilizer is bounded. A {\it conical limit point} is a point
$y\in Z$ such that there exists a sequence $(g_j)_{j\ge 0}$ in $G$, and distinct
points $a, b\in Z$, such that $g_j(y)$ tends to $a$ and $g_j(x)$ tends to $b$ for all $x \in Z\setminus \{y\}$.
We finally say that the action of $G$ on $Z$ is {\it geometrically finite} if every point of its limit set $\La_G$  is either conical or bounded parabolic
(they cannot be both simultaneously). When $G$ is a Kleinian group and $Z=S^2$, this definition is equivalent to the one given in \textsection \ref{sec:kleinian}, see for instance \cite{bowditch:gfps:jfa} and the references therein.

When a relatively hyperbolic group $(G,\PP)$ admits a cusp uniform action on a proper geodesic hyperbolic space $X$ then
its action on the boundary $\partial X$ is geometrically finite, cf. \cite{bowditch:rhg}.  The conjugates  of subgroups in $\PP$ are precisely the maximal parabolic subgroups. 
Conversely,  if a group $G$ admits a geometrically finite action on a metrizable space $Z$, then
the pair $(G,\PP)$ is  relatively hyperbolic \cite{yaman:crelle}, where  $\PP$ denotes a set of {conjugacy} representatives of maximal parabolic subgroups
of $G$. It turns out that the topology of $\La_G$ is independent of the space $Z$ as long as the action is geometrically finite with
the same parabolic subgroups  \cite[Theorem 9.4]{bowditch:rhg}.
 See \cite{rhuska:qcvxrhg} for more information.

A subgroup $H$ of a relatively hyperbolic group  $(G,\PP)$ is {\it elementary} if, whenever $(G,\PP)$ admits a cusp uniform action
on a proper geodesic hyperbolic space, the limit set of $H$ has at most two points.

\subsection{Horoballs}

Let us set some facts about the metrics on the combinatorial horoballs $H_P$ defined above.

\begin{fact}\label{fact:disthoro}
The distance on $H_P$ has the following behavior:
$$d\left[\begin{pmatrix} x\\ m\end{pmatrix},\begin{pmatrix} y\\ n\end{pmatrix}\right]  =  \max\{  |m-n|, 2\log_2 (1+d_P(x,y)) - (m+n)\}+O(1)$$
where $d_P$ denotes the metric on $P$.
\end{fact}

\Proof
By construction $H_P$ is a  geodesic space. We denote by $\ell(.)$ the length of a path,
i.e., the number of edges it contains. {A {\em horizontal path} will denote a path contained in some level $P\times\{ n\}$ for some $n\ge 0$ and a {\em vertical path} will be of the form $\{p\}\times [m,n]$ for some $p\in P$ and some integers $n\ge m\ge 0$.} Denote by $\left[\begin{pmatrix} x\\ n\end{pmatrix},\begin{pmatrix} y\\ n\end{pmatrix}\right]_h$ a horizontal path 
with minimal length and by $\left[\begin{pmatrix} x\\ m\end{pmatrix},\begin{pmatrix} x\\ n\end{pmatrix}\right]_v$ a vertical path.

It is easy to see that 
\begin{equation}    \label{eq:x=y}
d\left[\begin{pmatrix} x\\ n\end{pmatrix},\begin{pmatrix} x\\ m\end{pmatrix}\right]=|m-n|   
\end{equation}

When $x\neq y$, we remark that:
$$\ell\left(\left[\begin{pmatrix} x\\ n\end{pmatrix},\begin{pmatrix} x\\ n+1\end{pmatrix}\right]_v\cup
\left[\begin{pmatrix} x\\ n+1\end{pmatrix},\begin{pmatrix} y\\ n+1\end{pmatrix}\right]_h\right)\leq 
\ell\left(\left[\begin{pmatrix} x\\ n\end{pmatrix},\begin{pmatrix} y\\ n\end{pmatrix}\right]_h\cup
\left[\begin{pmatrix} y\\ n\end{pmatrix},\begin{pmatrix} y\\ n+1\end{pmatrix}\right]_v\right)$$
Therefore, if
$\begin{pmatrix} x\\ m\end{pmatrix}$ and $\begin{pmatrix} y\\ n\end{pmatrix}$ are two points in $H_P$, then they are joined by a geodesic
of the form $$\left[\begin{pmatrix} x\\ m\end{pmatrix},\begin{pmatrix} x\\ k\end{pmatrix}\right]_v\cup \left[\begin{pmatrix} x\\ k\end{pmatrix},\begin{pmatrix} y\\ k\end{pmatrix}\right]_h\cup\left[\begin{pmatrix} y\\ k\end{pmatrix},\begin{pmatrix} y\\ n\end{pmatrix}\right]_v$$ with $k\ge \max\{m,n\}$. Hence $$d\left[\begin{pmatrix} x\\m\end{pmatrix},\begin{pmatrix} y\\ n\end{pmatrix}\right]=\min_{k\geq\max\{m,n\}}\{2k-m-n+\lceil2^{-k}{d_P}(x,y)\rceil\}$$
with $\lceil a\rceil =p$ if $p-1< a\leq p$ for some $p\in\Z$. 
{Set $j=\lfloor\log_2({d_P}(x,y))\rfloor-k$ {then} 
$$2k-m-n+\lceil 2^{-k}{d_P}(x,y)\rceil=2\lfloor\log_2({d_P}(x,y))\rfloor-m-n-2j+\lceil 2^{j+\fraction}\rceil$$
where $\fraction=\log_2({d_P}(x,y))-\lfloor\log_2({d_P}(x,y))\rfloor$. {Set $f_\varepsilon(j)=-2j+\lceil 2^{j+\varepsilon}\rceil$. 
\begin{equation}    \label{value min}
\text{For } 0\leq\varepsilon<1,\text{let } g(\varepsilon)=\min\{f_\varepsilon(j),j\in\Z\}\in\{0,1,2\}
\end{equation}
This minimum is  always reached for $j=1$ (it might also be reached for $j=0$ or $j=2$ depending on the value of $\varepsilon$) and {$f_\ep$}
is decreasing for $j\leq 1$. Set $A=\lfloor\log_2({d_P}(x,y))\rfloor-\max\{m,n\}$. It follows that
\begin{equation}    \label{eq:k petit}
d\left[\begin{pmatrix} x\\m\end{pmatrix},\begin{pmatrix} y\\ n\end{pmatrix}\right]=|m-n|+ \lceil 2^{A+\fraction}\rceil \text{ when } 
{A\leq 1}
\end{equation}
and
\begin{equation}    \label{eq:k grand}
d\left[\begin{pmatrix} x\\ m\end{pmatrix},\begin{pmatrix} y\\ n\end{pmatrix}\right]=2\lfloor\log_2({d_P}(x,y))\rfloor-(m+n)+g(\fraction) \text{ when } 
A\geq 1. 
\end{equation}

We conclude the proof by combining equations (\ref{eq:x=y}), (\ref{eq:k petit}), (\ref{eq:k grand}) and (\ref{value min}) with the fact that $\log_2 a\leq \log_2 (1+a)\leq \max\{1,\log_2 a\} + 1$ for $a >0$.}
\endp

We recall the following fact that we leave as an exercise for the reader.

\begin{fact}\label{fact:3} Let $f:X\to Y$ and $g:Y\to X$, $\la\ge 1$ and $c,M\ge 0$ be such that
\ben
\item   for all $x,x'\in X$,  $ d_Y(f(x),f(x'))\le \la d_X(x,x')+ c$ and  for all $y,y'\in Y$,  $ d_X(g(y), g(y'))\le \la d_Y(y,y')+ c$;
\item for all $x\in X$, $d(g\circ f(x),x)\le M$ and for all $y\in Y$, $d(f\circ g(y),y)\le M$.
\een 
Then $f$ is a quasi-isometry.\end{fact}

{From this fact, we will deduce:}

\begin{fact}    \label{fact:close qi}
Consider a {$(\la,c)$}-quasi-isometry {$f:X\to X'$ between metric spaces} that maps a {path connected} subset $Y\subset X$ at Hausdorff distance {at most $D$} from {another path connected subset $Y'\subset X'$ for some $D>0$}. {Assume that there is a function $h:\R^+\to\R^+$ that is bounded on compact sets and such that given $x,y\in Y$, resp. $x',y'\in Y'$, then $d_Y(x,y)\leq h(d(x,y))$, resp. $d_{Y'}(x',y')\leq h(d(x',y'))$.} Then {$Y$ and $Y'$} {equipped with the length metrics $d_Y$ and $d_{Y'}$ induced by the inclusions} are quasi-isometric.  {The quasi-isometry constants only depend on $\la,c, D$ and $h$.}
\end{fact}

Fact \ref{fact:close qi} is also used to extend quasi-isometric embeddings between metric spaces to their horoballings, see also \cite[Lemma 6.2]{groff}.

{When needed, we will apply Fact \ref{fact:close qi} to subgroups as follows. Given a finitely generated subgroup $H$ of a finitely generated group $G$, consider a finite generating set $S$ for $G$ so that $S\cap H$ generates $H$. Let $X$ be the Cayley graph of $(G,S)$ and $Y$ the full subgraph with vertex set $H$.} {The subspace $Y\subset X$ is path connected and we may construct a function $h$ satisfying the assumption of Fact \ref{fact:close qi}.

\Proof
We assume that $f:X\to X'$ and its quasi-inverse $g:X'\to X$ are $(\la,c)$-quasi-isometries {and that $d(g\circ f(x),x)\le c$ for all $x\in X$}. Let us denote the Hausdorff distance by $d_H$ and
let us suppose that $d_H(f(Y),Y')\le D$. 

For any $x\in Y$, there exists $\vp(x)\in Y'$ such that $d(f(x),\vp(x))\le D$. This defines a map $\vp:Y\to Y'$. Given $x,y\in Y$ let $k\subset Y$ be an arc joining $x$ to $y$ such that $\ell(k)\leq \lfloor d_Y(x,y)\rfloor+1$. On $k$, we consider $p+1$ points $x=x_0, x_1,\ldots, x_p=y$ such that $p=\lfloor d_Y(x,y)\rfloor +1$ and $d_Y(x_i,x_{i+1})\leq 1$. By construction, we have $d(\vp(x_i),\vp(x_{i+1}))\leq 2D+\lambda+c$. It follows that $d_{Y'}(\vp(x_i),\vp(x_{i+1}))\leq\sup_{t\leq 2D+\lambda+c} h(t)=K$ and we get
{$$d_{Y'}(\vp(x),\vp(y))\leq p K\leq Kd_Y(x,y)+K.$$}

Moreover,
\begin{eqnarray*}
d_H(g(Y'),Y)  & \le & d_H(g(Y'), g(f(Y)) ) + d_H((g\circ f)(Y), Y) \\
& \le &  (\la d_H(Y',f(Y)) + c )+  c\\
&\le & \la D +2c
\end{eqnarray*}
Therefore, letting $\psi:Y'\to Y$ satisfy $d(\psi(x'),g(x'))\le \la D +2c$ for all $x'\in Y'$, we obtain as above {$K'$ such that:
$$ d_{Y'}(\psi(x'),\psi(y'))   \le  K' d(x',y') + K'$$}
{Additionally, for any $x\in Y$
\begin{eqnarray*}
d(\psi\vp(x),x)  & \le & d(\psi\vp(x), g\vp(x) ) + d(g\vp(x), gf(x)) + d(gf(x),x) \\
& \le &  (\la D +2c) + (\la  D + c ) +  c\\
&\le & 2\la D +{4}c
\end{eqnarray*}
We prove similarly that $d(\vp\psi,id)$ is bounded so that Fact \ref{fact:3} applies.}\endp

\begin{fact}\label{fact:6} Let $f:P\to Q$ be a quasi-isometric embedding. Then $F:H_P\to H_Q$ defined by
$$F\begin{pmatrix} x\\ m\end{pmatrix}= \begin{pmatrix} f(x)\\ m\end{pmatrix}$$ is a quasi-isometric embedding {where the constants
only depend on those of $f$}.\end{fact}

\Proof Let $\lambda,c$ be such that $\forall x,y\in P$, $\lambda^{-1} d(f(x),f(y))-c\leq d(x,y)\leq \lambda d(f(x),f(y))+c$. 

Then we have:
\begin{eqnarray*} 
d\left[F\begin{pmatrix} x\\ m\end{pmatrix},F\begin{pmatrix} y\\ n\end{pmatrix}\right] & \le & \max\{ |m-n| , 2  \log_2(1+ \la d(x,y)+c) - (m+n) \}+ O(1)\\
& \le & \max\left\{ |m-n|, 2\log_2 \left(\frac{1+c}{\la}+ d(x,y)\right)  - (m+n) \right\}  + O(1)\end{eqnarray*}
If $1+c\le\la$, we have $\log_2 \left(\dis\frac{1+c}{\la}+ d(x,y)\right)\leq \log_2 (1+ d(x,y))$.
Otherwise, 
\begin{eqnarray*}\log_2 \left(\frac{1+c}{\la}+ d(x,y)\right) & = &  \log_2\left(\frac{1+c}{\la}\right)+\log_2\left(1+\frac{\la}{1+c}d(x,y)\right)\\
& \leq & \log_2(1+d(x,y))+\log_2\left(\frac{1+c}{\la}\right)\,.\end{eqnarray*}
In both cases, we get 

\begin{eqnarray*} 
d\left[F\begin{pmatrix} x\\ m\end{pmatrix},F\begin{pmatrix} y\\ n\end{pmatrix}\right] & \le & d\left[\begin{pmatrix} x\\ m\end{pmatrix},\begin{pmatrix} y\\ n\end{pmatrix}\right] +O(1)\,.\end{eqnarray*}

On the other hand,  if $1+\la^{-1} d(x,y)-c>{1}$, we have
\begin{eqnarray*} 
d\left[F\begin{pmatrix} x\\ m\end{pmatrix},F\begin{pmatrix} y\\ n\end{pmatrix}\right] & \ge & \max\{ |m-n| , 2  \log_2(1+ \la^{-1} d(x,y)-c) - (m+n) \}+ O(1)\\
&\ge & d\left[\begin{pmatrix} x\\ m\end{pmatrix},\begin{pmatrix} y\\ n\end{pmatrix}\right] +O(1)\,,\end{eqnarray*}
{where we have used that $\log_2(1+ \la^{-1} d(x,y)-c) \ge \log_2 (1+ d(x,y)) + \log_2 1/(1+ c\la)$ since $d(x,y)\ge c\la$.   }
{If $1+\la^{-1} d(x,y)-c\le 1$, then both $d(x,y)$ and $d(f(x),f(y))$ are bounded, so the result follows. }\endp

The following fact is then easily deduced from Facts \ref{fact:3} and \ref{fact:6}:

\begin{fact}\label{fact 7} Let $f:P\to Q$ be a quasi-isometry. Then $F:H_P\to H_Q$ defined by
$$F\begin{pmatrix} x\\ m\end{pmatrix}= \begin{pmatrix} f(x)\\ m\end{pmatrix}$$ is a quasi-isometry, {where the constants
only depend on those of $f$}.\endp\end{fact}

We now prove an analogous statement for horoballs in real hyperbolic space.

\begin{fact}\label{fact:horohyp} Let $P$ be a graph endowed with the path metric so that each edge is isometric to $[0,1]$ and let $\vp:P\to  (\R^d, d_E)$ be a quasi-isometric embedding.
Consider a sequence $\{y_n\}\in\R^\N$ such that $1\leq y_n\leq L$ for any $n\in\N$ and some $L>0$ and define $\Phi:H_{P}\to \R^d\times\R_+^*\approx\HH^{d+1}$ by $\Phi(x,n)= (\vp(x),y_n 2^n)$. Then $\Phi$ is a quasi-isometric embedding. Furthermore, given $Y\supset \vp(P)$ and $c$ such that a $c$-neighborhood of $\vp(P)$ covers $Y$, then {$\Phi:H_{P}\to Y\times [1,\infty)$} is a quasi-isometry.
\end{fact}

\Proof 
Note that in the upper half-space model,
\begin{eqnarray*}
d_{\HH^{d+1}}((z,y),(z',y'))& =& \hbox{\rm arcosh} \left(1 + \frac{\|z-z'\|_E^2+(y-y')^2}{2yy'}\right)\\
&=&\ln\left(1 + \frac{\|z-z'\|_E^2+(y-y')^2}{2yy'}\right) + O(1)
\end{eqnarray*}

For any $a,b> 0$, $\max\{\ln a,\ln b\}\leq \ln(a+b)\leq \max\{\ln(2a),\ln(2b)\}$, hence $\ln(a+b)=\max\{\ln a,\ln b\}+O(1)$. Applied to the previous equality, this gives

\begin{eqnarray*}
d_{\HH^{d+1}}((z,y),(z',y'))& = & \max\left\{ \ln \left(\frac{\|z-z'\|_E^2}{2yy'}\right) , \ln\left(1 + \frac{(y-y')^2}{2yy'}\right) \right\}+ O(1)\\
& = &  \max\left\{ 2\ln (\|z-z'\|_E) -\ln {(yy')} , \ln\left(\frac{y}{2y'}+\frac{y'}{2y} \right) \right\}+ O(1)
\end{eqnarray*}

Let us assume that $y,y'\geq 1$, so that $\ln(yy')\geq 0$. When $\|z-z'\|_E\leq 1$, the term $2\ln (\|z-\nolinebreak z'\|_E) -\ln {(yy')}$ above is irrelevant. From the fact that, for any $a\geq 1$, $\ln a=\ln(1+a)+O(1)$, we deduce

\begin{eqnarray*}
d_{\HH^{d+1}}((z,y),(z',y'))& = &  \max\left\{ 2\ln (1 +\|z-z'\|_E) -\ln {yy'} , \left|\ln \left(\frac{y}{y'}\right)\right| \right\}+ O(1)
\,.
\end{eqnarray*}

Therefore, since $\{y_n\}$ is bounded,
\begin{eqnarray*}d_{\HH^{d+1}}((z,y_n 2^n),(z',y_m 2^m)) & = & 
\max\left\{ 2\ln (1 +\|z-z'\|_E) - (m+n)\ln 2) , |m-n| \ln 2 \right\}+O(1)\\
& = & \ln 2
\max\left\{ 2\log_2 (1 +\|z-z'\|_E) - (m+n)) , |m-n| \right\} + O(1)\,.\end{eqnarray*}
The rest of the proof follows as for Fact \ref{fact:6}.
\endp

\subsubsection{Pared groups}
\label{pared}

We adapt the notion of pared {three}-manifolds to groups. This point of view was first developed by Otal for free groups
in \cite{otal:equivlibre}.

Let $G$ be 
a finitely generated group. A {\it paring} will be given by a finite {\it almost malnormal} collection of subgroups $\PP_G=\{P_1,\ldots, P_k\}$; almost malnormal means that if $P_i\cap (gP_jg^{-1})$ is infinite for some $i,j\in\{1,\ldots, k\}$ and $g\in G$, 
then $i=j$ and $g\in P_i$. For example if $(G,\PP_G)$ is relatively hyperbolic, then it is a pared group. A pared compact {three}-manifold $(M,P_M)$ gives rise to a canonical pared group $(K,\PP_K)$ by letting $K=\pi_1(M)$ and $\PP_K$ denote a representative for each conjugacy class of subgroups corresponding to the fundamental groups of the annuli and tori composing $P_M$. Let us call this pared group the {\it pared fundamental group} of $(M,P_M)$.

To mimic more general $3$-manifolds with boundary we want to loosen the condition on $\PP_G$. An {\it adornment} is a finite collection of subgroups $\PP_G=\{P_1,\ldots, P_k\}$ that are non coarsely nested, i.e., given $P_i,P_j\in \PP_G$ and $g\in G$, if $g_iP_i$ is contained in some $D$-neighborhood of $P_j$ for some $D<\infty$ then $P_i=P_j$ and $g \in P_j$. It is not difficult to prove that a pared group is also adorned. A Seifert manifold $M$ with boundary admitting an $\HH^2\times \E^1$ structure gives rise to a canonical adorned group $(\pi_1(M),\PP_M)$ by letting $\PP_M$ denote the collection of the fundamental groups of the components of $\partial M$. The adorned fundamental group $(\pi_1(M),\PP_M)$ is not a pared group since $M$ is Seifert fibered and the fibers induce a normal cyclic subgroup that is a subgroup of all elements in $\PP_M$.

An isomorphism between {group pairs} $(G,\PP_G)$ and $(Q,\PP_Q)$ is an isomorphisms $G\to Q$ that maps the conjugacy classes in $\PP_G$ into the conjugacy classes in $\PP_Q$.

A quasi-isometry between two {adorned} groups $(G,\PP_G)$ and $(H,\PP_H)$ will be
given by a quasi-isometry $\Phi:G\to H$ {such for any $g\in G$ and $P\in\PP_G$ there is $h\in H$ and $Q\in \PP_H$ such that $\Phi(gP)$ is at bounded Hausdorff distance from $hQ$ and for any $h\in H$ and $Q\in \PP_H$ there is $g\in G$ and $P\in\PP_G$ such that $\Phi(gP)$ is at bounded Hausdorff distance from $hQ$. The  non coarse nestedness of the subgroups guarantees the uniqueness of $h$ and $Q$ in the first part and of $g$ and $P$ in the second.
Note that if $\PP_G$ is a paring for $G$, then the quasi-isometry implies that $\PP_H$ is also almost malnormal}. We will say that a
pared manifold $(M,P_M)$ is quasi-isometric to a pared group $(G,\PP_G)$ if there is a quasi-isometry between
its pared fundamental group $(K,\PP_K)$ and $(G,\PP_G)$.

{In the sequel, it will be convenient to write that two subsets are {\it at bounded distance} to mean they are at bounded Hausdorff distance.}

\begin{defn}[{induced adornment}]   \label{def:inducedgrouppair}
Let $(G,\PP)$ be {an adorned group} and assume that $H$ is a finite index subgroup of $G$, the {induced adornment} is
defined as follows. 

For every $P\in\PP$, let $T_P$ be a set 
of representatives of the double {cosets} $\{HgP, g\in G\}$ containing $1$ and set
$$\Q \fdefeq \{aPa^{-1}\cap H: P\in\PP, a\in T_P\}\,.$$
\end{defn}
Observe that since $H$ has finite index, each transversal $T_P$ is finite.

 Let us check that $(H,\Q)$ is an adorned group. Take $h\in H$, $a_i\in T_{P_i}$, $a_j\in T_{P_j}$ and $P_i,P_j\in\PP$ and assume that
$a_iP_ia_i^{-1}\cap H$ is contained in a finite neighborhood of $h(a_jP_ja_j^{-1}\cap H)$. Since $H$ has finite index, it follows that 
$a_iP_ia_i^{-1}$ is contained in a finite neighborhood of $h(a_jP_ja_j^{-1})$. Moreover, $a_iP_ia_i^{-1}$ is at bounded distance from $a_iP_i$, and similarly
for $j$ so that $a_iP_i$ is contained in a finite neighborhood of $ha_jP_j$. Therefore, $P_i=P_j$ and $ha_j \in a_i P_i$, so $a_j\in Ha_i P_i$ implying that
$a_i=a_j$ and $h\in a_iP_ia_i^{-1}\cap H$. This means that $(H,\Q)$ is an adorned group.

If  $HaP=HbP$ for some $P\in\PP$ and $a,b\in G$, then there exists $h\in H$ and $p\in P$ such that $b= hap$. Therefore
$$bPb^{-1}\cap H = hapPp^{-1}a^{-1}h^{-1}\cap H =  h(aPa^{-1}\cap H) h^{-1}$$
thus $bPb^{-1}\cap H$ is conjugate in $H$ to $aPa^{-1}\cap H$. This means that the choice made when picking $T_P$ only determines the representatives 
of the conjugacy classes of the elements in $\Q$. 

\begin{fact}\label{fact:inducedparing} 
Let $(G,\PP)$ be an adorned group and let  $H$ be a finite index subgroup of $G$. Let $(H,\Q)$ be the adorned group induced by $(G,\PP)$.  
If $G$ is finitely generated, then the canonical inclusion $H\hookrightarrow G$ extends to a quasi-isometry
between adorned groups $\vp : (H,\Q) \to (G,\PP)$.
Moreover, if $(G,\PP)$ is a pared group, then $(H,\Q)$ is also a pared group.
\end{fact}

\Proof 
Let us first prove that the canonical injection defines a quasi-isometry between {adorned} groups.
Let us fix a word metric $|\cdot |$ on $G$ provided by a finite generating set. We shall denote by $d(\cdot,\cdot)$ the Hausdorff distance between non-empty subsets of $G$.

Let us observe that as $\PP$ and $T_P$ are finite sets and each subgroup $aPa^{-1}\cap H$ has finite index in $aPa^{-1}$ for all $P\in\PP$ and $a\in T_P$, 
there is a constant $B$ such that $d(aPa^{-1}\cap H,aPa^{-1})\le B$ for all $P\in\PP$ and $a\in T_P$. 
We may as well assume that $|a|\le B$ holds for every $a\in T_P$. Thus, by the triangle inequality, we have
$$d(aPa^{-1}\cap H, aP) \le d(aPa^{-1}\cap H, aPa^{-1}) + d(aPa^{-1},aP)\le  B +|a| \le 2B\,.$$

Let us start by showing that each coset $h(aPa^{-1}\cap H)$ in $H$, where $P\in \PP$, $a\in T_P$ and $h\in H$, is at bounded distance from the coset $(ha)P$ in $G$. 
We have
$$d(h(aPa^{-1}\cap H), (ha)P)= d(aPa^{-1}\cap H, aP) \le  2B\,.$$

Conversely, if $gP$ is a coset, let us write $g=hap$ with $h\in H$, $a\in T_P$ and $p\in P$, and let us estimate $d(gP, h(aPa^{-1}\cap H))$:
$$d(gP, h(aPa^{-1}\cap H)) = d(aP,aPa^{-1}\cap H) \le 2B\,.$$
Therefore, we have a quasi-isometry between the {adorned groups} $(H,\Q)$ and $(G,\PP)$.

If $\PP$ is an almost malnormal collection in $G$, let us show that $\Q$ is also an almost malnormal collection in $H$. Consider $Q_1,Q_2\in\Q$, $h\in H$ and let us assume that $hQ_1h^{-1}\cap Q_2$ is infinite. We may find $P_1,P_2\in \PP$, $a_1\in T_{P_1}$ and
$a_2\in T_{P_2}$ such that  $Q_j= a_jP_ja_j^{-1}\cap H$ for $j=1,2$. Thus
$$hQ_1h^{-1}\cap Q_2 =  a_2[  (a_2^{-1}ha_1)P_1 (a_2^{-1}ha_1)^{-1}\cap P_2]a_2^{-1}\cap H$$
so the almost malnormality of $\PP$ implies that $P_2=P_1=P\in\PP$ and $a_2^{-1}ha_1\in P$. Hence, we may find $p\in P$
such that $h a_1= a_2p$. By definition of $T_P$, this implies that $a_1=a_2=a\in T_P$. Hence $Q_1=Q_2=Q$. Furthermore, we have
$a^{-1}ha\in P$, hence $h\in aPa^{-1}\cap H=Q$.
\endp
}

\subsection{Quasi-isometries between {adorned groups}}

Let us use the results of the previous section to study quasi-isometries between {adorned group}.

\begin{theorem}\label{thm:qisomgg} If two {adorned groups} $(G,\PP)$ and  $(G',\PP')$ are quasi-isometric then the cusped spaces $\cus(G,\PP)$ and $\cus(G',\PP')$ are also quasi-isometric.\end{theorem}

\Proof 
Let $\vp:(G,\PP)\to (G',\PP')$ be a quasi-isometry. Up to increasing the additive constant if necessary, we may assume that,  for any $g\in G$ and $P\in\PP$, there are $g'\in G'$ and
$P'\in\PP'$ such that $\vp: gP\to g'P'$ is a quasi-isometry, cf. Fact \ref{fact:close qi}. 

Set $\Phi:\cus(G,\PP)\to \cus(G',\PP')$ by letting $\Phi=\vp$ on $\cay(G)$ and extending $\vp$ by Fact \ref{fact 7} to $\Phi:H_{gP}\to H_{g'P'}$ for any $g\in G$ and $P\in\PP$.
We note that {there are constants $\la,c$ such that} each extension is a $(\la,c)$-quasi-isometry on the corresponding horoball, and that $\vp$ is also a $(\la,c)$-quasi-isometry
for the graph metric of $\cay(G)$. {Since the cosets are coarsely non-nested, it follows that $\Phi$ is quasi-surjective.} 

We conclude with a subdivision argument as in the proof of \cite[Theorem 6.3]{groff}.
Let $x,x'\in \cus(G,\PP)$ and let us consider a geodesic $[x,x']$. We may find a subdivision $(x_j), {0}\le j\leq p$, of the geodesic segment such that
each segment $[x_j,x_{j+1}]$ is {non degenerate (hence have length at least $1$) and is} either contained in $\cay(G)$
or in some $H_{g_jP_{k_j}}$. 
{Then we have $p\leq d(x,x')$ and}
\begin{equation*} 
d(\Phi(x),\Phi(x'))\le \sum d(\Phi(x_j),\Phi(x_{j+1})) \le \sum (\lambda d(x_j,x_{j+1})+c) \le  {(\lambda+c)d(x,x')}\,.
\end{equation*}

By symmetry, we obtain the same inequality for its quasi-inverse, hence Fact \ref{fact:3} concludes the proof (condition $(2)$ is satisfied by construction since $\Phi_{|G}:G\to G'$ and $\Phi_{|H_{gP}}:H_{gP}\to H_{g'P'}$ are quasi-isometries).
\endp

\begin{corollary} \label{cor:bowbdry}
The topology and quasi-M\"obius class of the boundary of the cusped space of a relatively hyperbolic group is independent from
the choice of the generating set $S$.\end{corollary}

\begin{corollary} \label{cor:qi pared hyperbolic}
If an {adorned group} is quasi-isometric to a relatively hyperbolic group, then it is hyperbolic relative to its {adornment}.
\end{corollary}

If $K$ is a geometrically finite Kleinian group, then $K$ is finitely generated, and relatively hyperbolic with respect to its maximal parabolic
subgroups \cite[Theorem 5.1]{farb:relhyp}. Let $\PP$ be a set of representatives of their conjugacy classes. We may then define $\cus(K){=\cus(K,\PP)}$ to be the cusped space
of $K$ by adding horoballs to the orbits of the elements of $\PP$.

\begin{proposition}\label{prop:qisomgfk} Let $K$ be a geometrically finite Kleinian group and let $\hbox{\rm Hull}(\La_K)$ be the convex hull of $\La_K$. Pick a point $o\in \hbox{\rm Hull}(\La_K)$. Then the map $k\mapsto k(o)$ extends to  a quasi-isometry
$\Phi: \cus(K)\to \hbox{\rm Hull}(\La_K)$ such that $\Phi(\partial \cus(K))=\La_K$.\end{proposition}

\Proof 
By assumption, we have $\Phi(k)=k(o)$ for any $k\in K$. We extend $\Phi$ to each horoball as in Fact \ref{fact:horohyp}. {Explicitly, we proceed as follows.}

Let $U$ be a {union} of $K$-invariant horoballs {in $\HH^3$} attached to the parabolic points of $K$. 
{Using Margulis' lemma} we may assume that they are pairwise disjoint {(see also \cite[Proposition VI.A.11]{maskit:kleinian_groups})} and that $o\not\in U$ by decreasing their size.
Fix $P\in \PP$ and denote by $z_P\in\cbar$ its parabolic fixed point. 
Represent $\HH^3$ by the upper half-space model $\R^2\times\R_+^*$ so that $o$ becomes $(0,0,1)$ and $z_P$ corresponds to the point at infinity. Thus, $\Phi(P)\subset \R^2\times\{y_0\}$, where $y_0=1$.  
Let $y_1$ be large enough so that $\R^2\times\{2y_1\}\subset U$ and set $y_n=y_1$ for $n\geq 1$. Let $p:\R^2\times\R_+^*\to \R^2$ be the  projection along the third coordinate and define $\Psi_P: H_P\to \HH^3$ by $\Psi_P(k,n)= (p\circ\Phi(k),y_n 2^n)$. Now define $\Phi$ on $H_P$ by $\Phi(k,n)=\Psi_P(k,n)$. We do the same construction for each parabolic subgroup $P\in \PP$ and we extend the result equivariantly to get a map $\Phi: \cus(K)\to Hull(\La_K)$. 
The next Claim will conclude the proof of Proposition \ref{prop:qisomgfk}. \endp

\noindent{\bf Claim. ---}  The map $\Phi: \cus(K)\to \hbox{\rm Hull}(\La_K)$ is a quasi-isometry.

\Proof
Up to taking smaller horoballs {that define} $U$, we may assume that for any parabolic subgroup {$P\in\PP$} 
and any $k\in P$, $\Phi(k,1)\subset\partial U$. 
Let $K_+\subset \cus(K)$ be the union of $K$ and the vertices at height $1$ in the horoballs and denote by $Y_K$ 
the maximal induced  subgraph of $\cus(K)$ with vertex set $K_+$. We denote by $Z_K$ the maximal subgraph of $\cus(K)$ 
whose vertices have height at least $1$ in some horoball. Thus we have $\cus(K)=Y_K\cup Z_K$ and $Y_K\cap Z_K$ 
is the maximal subgraph of $\cus(K)$ whose vertices have height exactly $1$ in some horoball.

Since {$K$ is geometrically finite, the actions of $K$ on $\overline{\hbox{\rm Hull}(\La_K)\setminus U}$ is cocompact  \cite[Theorem 4.73]{kapovich:book} and, by construction, its action on $Y_K$ is cocompact. It follows that} 
the restriction of 
$\Phi$ to $Y_K$ is a quasi-isometry between $Y_K$ and $\overline{\hbox{\rm Hull}(\La_K)\setminus U}$ endowed with the induced length metrics.

Let $U_i$ be a horoball {defining} $U$ and represent $\HH^3$ by the upper half-space model $\R^2\times\R_+^*$ so that $U_i=\{(x,y,z),z\geq 2y_1\}$. Then $\hbox{\rm Hull}(\La_K)\cap U_i$ has the form $\{\hbox{\rm Hull}(\La_K)\cap\partial U_i\}\times[2y_1,\infty)$. {It follows from the construction and the previous paragraph that $\Phi (Y_K\cap Z_K)$ lies at bounded Hausdorff distance from  $\hbox{\rm Hull}(\La_K)\cap\partial U$. By Fact \ref{fact:close qi}, $\Phi(Y_K\cap Z_K)$ and $\hbox{\rm Hull}(\La_K)\cap\partial U$, each equipped with the induced path metric, are quasi-isometric and by Fact \ref{fact:horohyp}, for any horoball $H_{kP}$, the restriction of $\Phi$ to $H_{kP}\cap Z_K$ is a quasi-isometry to the corresponding component of $\hbox{\rm Hull}(\La_K)\cap U$.}

In the previous paragraphs, we have seen that the restrictions $\Phi_{|Y_K}:Y_K\to \overline{\hbox{\rm Hull}(\La_K)\setminus U}$ and $\Phi_{|Z_K}:Z_K\to \hbox{\rm Hull}(\La_K)\cap U$ are quasi-isometries. It follows that there is $D$ such that the $D$-neighbourhood of $\Phi(\cus(K))$ covers $\hbox{\rm Hull}(\La_K)$.

Pick $k,k'\in \cus(K)$. With a subdivision as in the end of the proof of Theorem \ref{thm:qisomgg} we get:
$$ d(\Phi(k),\Phi(k'))\le 3\max\{\la,c\}(d(k,k')+1)\,.$$

On the other hand, we may decompose $[\Phi(k),\Phi(k')]$ into finitely many segments $[x_j,x_{j+1}]$, ${0}\leq j< p$ so that $x_0=\Phi(k)$, $x_p=\Phi(k')$,
$x_j\in\partial U$ for $0< j<p$ and {$(x_j,x_{j+1})\cap \partial U=\emptyset$} for $0\leq j<p$. For each index $0<j<p$, there is $k_j\in Y_K\cap Z_K$ 
so that $d(x_j,\Phi(k_j))\le c$ and we have $\frac{1}{\lambda} d(k_j,k_{j+1})\leq d(x_j,x_{j+1})+3c$ for any $j<p$, taking $k_0=k$ and $k_p=k'$. 
Let  {$T>0$} be the minimal distance between two horoballs in the family $U$. If {$(x_j,x_{j+1})\cap U=\emptyset$}, then 
$d(x_j,x_{j+1})\geq T$ and $d(x_j,x_{j+1})+6c\leq (\frac{6c}{T}+1)d(x_j,x_{j+1})$. 
Since there are at least $\lfloor\frac{p}{2}\rfloor$ such segments, we have {(shifting in the sum the constant $3c$ of segments off $U$ to 
a neighboring segment in $U$)}:

\begin{multline*}
d(k,k')\leq \sum d(k_j,k_{j+1})\leq \lambda\sum (d(x_j,x_{j+1})+3c)\\ \leq 3c\lambda+ \lambda\left(\frac{6c}{T}+1\right)\sum d(x_j,x_{j+1})=\lambda\left(\frac{6c}{T}+1\right) d(\Phi(k),\Phi(k'))+3c{\lambda}\,.
\end{multline*}

 \endp

\subsection{Quasiconvexity} 
Let $X$ be a proper, geodesic, hyperbolic metric space. A $K$-{\it quasiconvex} subset $Y\subset X$
has the property that any geodesic segment joining two points of $Y$ remains in the $K$-neighborhood of $Y$.
Note that quasiconvexity is a property invariant under quasi-isometries. 

Given a non trivial compact subset $\La\subset \partial X$, we define its {\it weak convex hull} (or {\it join}) $\CCC(\La)\subset X$ as the union of all geodesics joining pair of points of $\La$. Even though this set is usually not convex, it is shown in \cite{kapovich:short:greenberg} that it is uniformly quasiconvex. Note also that if $\La$ is a compact subset of $\cbar$,
then  the inclusion map   $\CCC(\La)\hookrightarrow  \hbox{\rm Hull}(\La)$
is a quasi-isometry in $\HH^3$. 

A subgroup $H$ of a hyperbolic group $G$ is {\it quasiconvex} if $H$ is quasiconvex in any locally finite Cayley
graph $\cay (G)$ of $G$. A subgroup $H$ of a relatively hyperbolic group $(G,\PP)$ is {\it relatively quasiconvex} if, for any cusped space $\cus(G,\PP)$
for $G$, there is a constant $L$ such {that} for any geodesic $\g\subset \cus(G,\PP)$ with endpoints in $H$, $\g\cap\cay (G)$ remains at distance
at most $L$ from $H$.

Moreover, according to \cite[Proposition 7.6]{rhuska:qcvxrhg}, if $H$ is a relatively quasiconvex subgroup of $G$ and $G$ has a cusp uniform action on $\cus(G,\PP)$, then either $H$ is finite or $H$ is parabolic or the action of $H$ {admits a cusp-uniform action on} some {hyperbolic proper geodesic metric space which is quasi-isometric to $\CCC(\La_H)$}. In the latter case, the maximal parabolic subgroups of $H$  define a finite number
of $H$-conjugacy classes and are  of the form $gPg^{-1}\cap H$
where $g\in G$, $P\in\PP$, and such that $gPg^{-1}\cap H$ is infinite. It follows that $H$ is hyperbolic relative to representatives of theses classes \cite[Theorem 9.1]{rhuska:qcvxrhg}. {Moreover, its limit set $\La_H\subset\partial\cus(G,\PP) $ is also the set of accumulation points of the image $H\hookrightarrow (\cus(G,\PP)\sqcup \partial\cus(G,\PP))$.}

\begin{proposition}\label{prop:qisomcvx}  Let $(G_1,\PP_1)$ and $(G_2,\PP_2)$ be two finitely generated relatively hyperbolic groups and let us consider
two infinite, non parabolic and relatively quasiconvex finitely generated subgroups $(H_1,\Q_1)$ and $(H_2,\Q_2)$ of $G_1$ and $G_2$ respectively.

Let $\vp:(G_1,\PP_1)\to (G_2,\PP_2)$ be a quasi-isometry between pared groups  together
with its extension  $\Phi: \cus(G_1,\PP_1)\to \cus(G_2,\PP_2)$ given by Theorem \ref{thm:qisomgg} and denote by $\partial\Phi:\partial \cus(G_1,\PP_1)\to\partial \cus(G_2,\PP_2)$ its boundary map.
If $\partial\Phi(\La_{H_1})=\La_{H_2}$ then $\vp(H_1)$ is at bounded distance from $H_2$ and   $(H_1,\Q_1)$ and
$(H_2,\Q_2)$ are quasi-isometric.
\end{proposition}

{We first establish some lemmas for its proof. Let us start with a geometric lemma.

\begin{lemma}\label{lma:truncatedlimitset} Let $X$ be a proper, geodesic, $\de$-hyperbolic metric space. For any $D>0$, there exists $D'=D'(\de,D)$ with the following property. 
Let $\La\subset \partial X$ be a compact subset with at least two points and $Y\subset X$ be the complement of a collection of pairwise disjoint  horoballs. Then
$$N_D(\CCC(\La))\cap N_D(Y) \subset N_{D'}(\CCC(\La)\cap Y)$$
where  $N_D(S)$ denotes the $D$-neighborhood of a subset $S$.
\end{lemma}

\Proof 
Let $z_0\in N_D(\CCC(\La))\cap N_D(Y)$. By definition, there exists $z_1\in \CCC(\La)$ at distance at most $D$ from $z_0$. If ever $z_1\in Y$, then we are done. Otherwise, $z_1$ belongs to a horoball {$H=H_w(p,r)$} centred at some point $p\in\partial X$. Let us consider a geodesic $(a,b)$ with endpoints in $\La$ that contains the point $z_1$. As $z_1\in H$, we may find $x \in [a,z_1)$ and $y\in [b,z_1)$ such that {$x,y\in \partial H\cup \{p\}$}. 
Since $H$ is a horoball, we may assume that $x\in X$. {Since $z_1\in H$, then $d(z_1,\partial H)= d(z_1,Y)\le 2D$. 
{Let $z_2\in\partial H$
realize this distance. It follows from \cite[Lemme 8.1]{ghys:delaharpe:groupes} that $\beta_p(.,w)$ is roughly Lipschitz, i.e., it holds $|\beta_p(x,w)-\beta_p(y,w)|\leq d(x,y)+O(\delta)$. In particular, since $x$ and $z_2$ belong to $\partial H$, we have $|\be_p(x,w)-r|\leq O(\delta)$ and $|\be_p(z_2,w)-r|\leq O(\delta)$. As $\be_p$ is almost a cocycle according to \cite[Proposition 8.2]{ghys:delaharpe:groupes}, we get $\be_p(x,z_2)=O(\de)$ and
$$\be_p(x,z_1) = \be_p(x,z_2) + {\be_p(z_2,z_1)} +O(\de) \le  d(z_1,z_2) +O(\de) \le  2D +O(\de)\,.$$}

If $y$ is at infinity, then we have $y=b=p$. In this
case $\be_p(x,z_1)= d(x,z_1)$, 
and, as $x\in \CCC(\La)\cap Y$, we obtain
$$d(z_0, \CCC(\La)\cap Y) \le d(z_0,x)\le d(z_0,z_1) + d(z_1,x)\le 3D + O(\de)$$
Let us now assume that $y\in X$ as well. 
Since triangles are slim, we may find two rays $[x,p)$ and $[y,p)$ so that the triangle $\{x,y,p\}$ is $\de$-slim. It follows that the point $z_1$ is $\de$-close to one of those
rays. Let us assume for the sake of simplicity that it is $[x,p)$. It follows that $2D\ge \be_p(x,z_1)  \ge d(x,z_1) + O(\de)\,.$ 
Therefore, as $x\in \CCC(\La)\cap Y$, we obtain $d(z_0, \CCC(\La)\cap Y)\le 3D + O(\de)$ and the lemma follows.\endp
}

We analyse {relatively quasiconvex} subgroups in the vicinity of their parabolic subgroups.

\begin{lemma}\label{lma:parabclsoe} {Let $G$ be a group acting properly discontinuously by isometries on a  proper hyperbolic geodesic metric space $X$}.
Let $H<G$ be a subgroup admitting a non-elementary  action on $\partial X$. 
{Fix  a parabolic point $a$ for $H$, and let  $Q=\stab_H a$ and $P=\stab_G a$. Assume that $a$ is a bounded parabolic point for $Q$.}
For any base point $o\in X$ and any $D>0$, there exists $D'$ such that $$N_D(Ho)\cap N_D(Po)\subset N_{D'}(Qo)\,.$$
\end{lemma}

\Proof
{Let us first notice that $a$ is also parabolic for $G$ as $H$ is a subgroup of $G$.} Consider a sequence $(x_n)_n\subset N_D(Ho)\cap N_D(Po)$ and $h_n\in H$ and $p_n\in P$ such that $d(h_n(o),x_n)\leq D$ and $d(p_n(o),x_n)\leq D$. We are going to show that $(h_n(o))_n$ lies in a bounded neighbourhood of $Qo$, the conclusion follows easily. 

Since $d(p_n^{-1}h_n(o),o) \le 2D$ and the action of $G$ is properly discontinuous, we may assume that there exists $g\in G$ such
that $h_n=p_ng$ for all $n$. If $(h_n)_n$ has finitely many elements, the conclusion is straightforward so we may assume the contrary.

We note that $(p_n)_n\subset P$ contains infinitely many distinct elements and the action of $P$ on $\partial X\setminus \{a\}$ is properly discontinuous {since $\La_P=\{a\}$ as $P$ is parabolic with fixed point $a$}{, cf. \cite[\S\S 3.1, 8.2]{gromov:hyperbolic} and \cite[Theorem 2L]{tukia:convergence_groups}}. Therefore,  for any  $z\in \partial X\setminus\{g^{-1}(a)\}$, 
the sequence  $(p_ng(z))_n=(h_n(z))_n$
tends to $a$.
Similarly, $(h_n^{-1}(z))_n$ tends to $g^{-1}(a)$ for all $z\in \partial X\setminus\{a\}$.

Since $Q$ is a bounded parabolic subgroup of $H$, there exists a compact subset $K\subset \La_H\setminus\{a\}$ that meets every $Q$-orbit in $\La_H\setminus\{a\}$. Let us pick $z\in\La_H\setminus \{g^{-1}(a)\}$. Since $h_n(z)\in \La_H\setminus\{a\}$, we may find a sequence $(q_n)$ in $Q$ such that $q_n(h_n(z))\in K$. If the sequence $(q_nh_n)_n$ is finite then $d(h_n(o),Qo)$ is bounded and we are done.

Let us assume that $(q_nh_n)_n$ is infinite. Up to passing to a subsequence,  we may assume that  $(q_nh_n)$ is a collapsing sequence, i.e., there are two points $\al,\be\in\La_H$ such that
$(q_nh_n)_n$ tends uniformly to $\be$ on the compact subsets of $\partial X\setminus\{\al\}$ {\cite[Prop.\,2.1]{bowditch:convergence_groups}}. Then  $(q_nh_n)^{-1}_n$ tends uniformly to $\al$ on the compact subsets of $\partial X\setminus\{\be\}$. Since $a$ is not conical, the limit of $((q_nh_n)^{-1}(a))_n$ has to be $\al$, hence $\al=g^{-1}(a)$ since $(q_nh_n)^{-1}(a)=g^{-1}(a)$ holds for all $n$.
It then follows that $\be \in K$ since $z\ne \al$. {Since $(q_np_n)_n$ is an infinite sequence in $P$, which is a parabolic subgroup with fixed point $a$, then $(q_nh_n(a))_n=(q_np_n(g(a)))_n$ tends to $a$. By construction $a\not\in K$, hence $(q_nh_n(a))_n$ can not tend to $\beta\in K$. It follows that $a=\alpha$.} 
Therefore, we conclude that $a=\alpha=g^{-1}(a)$ and $g\in P$. It follows that $h_n\in H\cap P=Q$ for all $n$. 
{Applying these arguments on subsequences, we get that for any subsequence such that $(q_nh_n)_n$ is infinite, there is a further subsequence such that $h_n\in H\cap P=Q$. Thus we have proved that either $(q_nh_n)_n$ is finite or there are only finitely {many} indices such that $h_n\not\in Q$ and the conclusion follows.}
\endp

\demode{Proposition \ref{prop:qisomcvx}} Let $j=1,2$ and let us denote by $\cay_j$ a locally finite Cayley graph  of $G_j$ used to build the cusped space $\cus_j=\cus(G_j,\PP_j)$. 
As mentioned above, the action of $H_j$ on $\CCC(\La_{H_j})\ (\subset \cus_j)$ is cusp uniform.  
Hence the action of $H_j$ on the truncated space $\CCC_T(H_j)= \CCC(\La_{H_j})\cap \cay_j$ is cocompact. It follows that the Hausdorff distance between $H_j\subset \cus_j$ and $\CCC_T(H_j)$ is bounded since they are both
$H_j$-invariant. 

By the shadowing lemma (Lemma \ref{lm:shadow}), $\Phi(\CCC(\La_{H_1}) )$ is quasiconvex, at bounded distance from $\CCC(\La_{H_2})$. By assumption, $\Phi(\cay_1)=\vp(\cay_1)$ is at bounded distance from $\cay_2$. It follows then from the previous paragraph that {there is some $D>0$ such that $\Phi(\CCC_T(H_1))= \vp(\CCC_T(H_1)) \subset N_D(\cay_2)\cap N_D(\CCC(\La_{H_2}))$. By Lemma \ref{lma:truncatedlimitset}, we obtain that it is contained in a neighborhood of $\CCC_T(H_2)$. Thus }
 $\Phi(\CCC_T(H_1))$ is at bounded distance from $\CCC_T(H_2)$ as well, implying that $\Phi(H_1)$ is at bounded
distance from $H_2$ in $\cus_2$. Let us note that the canonical injection $\cay_j\hookrightarrow \cus_j$, $j=1,2$, is {a coarse embedding} by Fact \ref{fact:disthoro}, so
we may conclude that $\vp(H_1)$ is at bounded distance from $H_2$ in $\cay_2$.

It remains to prove that the quasi-isometry preserves the parings.

Let $Q_1\in\Q_1$ be  a parabolic subgroup of $H_1$ with parabolic point $p$. There are $g_1\in G_1$  and $P_1\in\PP_1$  such that $Q_1= g_1P_1g_1^{-1}\cap H_1$. 
By assumption, we may find $P_2\in\PP_2$ and $g_2\in G_2$ such that $\vp(g_1P_1)$ is at bounded distance from $g_2P_2$. 
Let us consider $Q_2= g_2P_2g_2^{-1} \cap H_2$.  {The point $q=\partial\Phi(p)$}  is a parabolic for $G_2$ {with
stabilizer $g_2P_2g_2^{-1}$}. {But the point $q$} is also {a bounded} parabolic for $H_2$ since $H_2$ is geometrically finite,
so $Q_2$ is infinite.

Let us observe that, for $j=1,2$,  $g_jP_jg_j^{-1}$ is parabolic, so preserves horospheres centered at $p$ or {$q$}, and $g_jP_j$ is also a horosphere centered at the
same point,
so $g_jP_j$ and $g_jP_jg_j^{-1}$ lie at bounded distance.  On the one hand, $\vp(g_1P_1)$ is at bounded distance from $g_2P_2$, and $\vp(H_1)$ is at bounded distance from $H_2$. It follows that $\vp(Q_1)$ is {in a bounded neighborhood of} $H_2$, {of} $g_2P_2$ {and} $g_2P_2g_2^{-1}$.
But $H_2\cap g_2 P_2 g_2^{-1} = Q_2$, so  Lemma \ref{lma:parabclsoe} now implies that $\vp(Q_1)$ is {in a bounded neighbourhood of} $Q_2$. Therefore,
$\vp$ maps the peripheral structure for $H_1$ into a  bounded neighborhood  of the peripheral structure of $H_2$. 
By symmetry, we conclude  that $(H_1,\Q_1)$ and $(H_2,\Q_2)$ are quasi-isometric as pared groups.
\endp

Let $(G,\PP_G)$ be a relatively hyperbolic group and let $\PP\supset \PP_G$ be a paring for $G$ by infinite relatively quasiconvex subgroups
{(by definition $\PP$ is an almost malnormal finite collection of subgroups)}. 
Define on the Bowditch boundary $\partial_{\PP_G} G=\partial \cus(G,\PP_G)$  an equivalence relation $\sim_\PP$ as follows: let $x\sim_\PP y$
if, either $x=y$ or if there is a subgroup $P\in \PP$ and an element $g\in G$ such that
$\{x,y\}\subset  g(\La_{P})$. Set  $Q_\PP =\partial_{\PP_G} G/\sim_\PP$ to be the quotient of $\partial_{\PP_G} G$
by this relation $\sim_\PP$.

Being unable to find the following generalization of  \cite[Theorem 7.11]{bowditch:rhg} in the literature,
we sketch a proof of it.

\begin{proposition}    \label{prop:relstructure}
With the above notation, $(G,\PP)$ is relatively
hyperbolic and there is a $G$-equivariant homeomorphism between $Q_\PP$ and the Bowditch boundary of $(G,\PP)$.
\end{proposition}

We start with the following lemma.

\begin{lemma}\label{lma:relstructure}  Let $(G,\PP_G)$ be a relatively hyperbolic group and let {$\PP\supset \PP_G$ be a paring by
{infinite} relatively quasiconvex subgroups.} 
The following properties hold.
\ben
\item {For any $H\in \PP\setminus\PP_G$, t}he action of  $H$  is uniform  on $\La_H$ and cocompact on $\partial_{\PP_G}G\setminus \La_H$.
\item {For any $H\in \PP\setminus\PP_G$, w}e have $H=\stab_G\La_H$.
\item The collection of compact sets $\KKK= \{g\La_H {\colon} g\in G{,\ H\in\PP\setminus\PP_G}\}$ forms a null sequence of pairwise disjoint compact sets, i.e., 
for any distance on  $\partial_{\PP_G}G$, for any $\de >0$, there are only finitely many sets in $\KKK$ with diameter at least $\de$.
\een
\end{lemma}

\Proof Let $X=\cus(G,\PP_G)$, and let us consider the action of $G$ on $X\cup \partial X$. 
{Let us consider $H\in \PP\setminus \PP_G$.} We note that $H$ has no parabolic {subgroups}. If this was the case, then there would be some $g\in G$ and $P \in \PP_G$  such that $gPg^{-1}\cap H$ would
be infinite, which contradicts the almost malnormality assumption. Therefore, since $H$ is relatively quasiconvex, this implies that every point in $\La_H$
is conical.  We may then deduce that the action of $H$ is uniform on $\La_H$ by \cite[Thm.\,8.1]{bowditch:characterization} and cocompact on $\partial_{\PP_G}G\setminus \La_H$  by   \cite[Main Thm.\,(3)]{swenson:qcvxgps}. 
Furthermore, by the corollary of the main theorem in \cite{swenson:qcvxgps}, $H$ has finite index in $\stab_G \La_H$. If $g\in \stab_G\La_H$, then 
$gHg^{-1} \cap H$ is a finite index subgroup of $H$, hence infinite, so that $g\in H$ by almost malnormality.  We have proved (1) and (2).

By the corollary to \cite[Thm.\,13]{swenson:qcvxgps}, we also know that, {for any $K\in \PP\setminus \PP_G$}, we have {$g\La_H\cap \La_K= \La_{gHg^{-1}\cap K}$} so that either {$H=K$ and} $g\in H$  or {$g\La_H\cap \La_K=\emty$}.
This shows that the elements of $\KKK$ are pairwise disjoint.
 Let us fix $\de>0$. There exists $R>0$ such that   $d(e,gH)\le R$ whenever $\diam g\La_H\ge \de$  since $H$ is quasiconvex in $X$   \cite[Main Thm.\,(1)]{swenson:qcvxgps}. 
This implies that there are only finitely many such elements $g\in G/H$. This concludes the proof. \endp

\demode{Prop.\,\ref{prop:relstructure}} 
We proceed in two steps. We first establish that $\sim_\PP$ defines an upper semi-continuous decomposition of $\partial_{\PP_G} G$, i.e.,
the equivalence relation $\sim_\PP$ is closed. This  implies that $Q_\PP$ is Hausdorff and compact.
Then we prove that the action of $G$ on $Q_\PP$ is geometrically finite with the prescribed parabolic subgroups.

By Lemma \ref{lma:relstructure},  the limit
sets $\{g\La_P,\ P\in\PP, g\in G\}$ form a null sequence of pairwise disjoint sets, so  they define an upper semi-continuous decomposition of 
$\partial_{\PP_G}  G$. 
This shows that the quotient $Q_\PP$  is a Hausdorff compact space, and the group $G$ acts on $Q_\PP$. 

Lemma \ref{lma:relstructure} also implies (a) that   the action of any $P\in\PP\setminus\PP_G$  on $\partial_{\PP_G} G\setminus \La_P$ 
is cocompact, so they define bounded parabolic groups on $Q_\PP$ and  (b) that they  are  maximal parabolic subgroups.

Let us now check that  all the other points are conical. If we consider such a point $z$ with preimage $x\in\partial_{\PP_G} G$, then we may find distinct points 
$\al,\be\in\partial_{\PP_G} G$ and a sequence of elements $(g_n)$ such that $(g_n(x))_n$ tends to $\al$ while all the other sequences $(g_n(y))_n$
tend to $\be$. If $\al$ and $\be$ lie in different fibers of the projection map $\pi:\partial_{\PP_G} G\to Q_\PP$, then it follows that $z$ is also conical. 
If they belong to a common fiber,  they belong to some $g\La_P$, $P\in\PP\setminus\PP_G$, $g\in G$ and we may as well assume that $g$ is trivial. 
As $x\notin\La_P$ and the action of $P$ on $\partial_{\PP_G} G\setminus \La_P$ is cocompact, we may find $(h_n)_n$ in $P$ so that 
$(h_ng_n)(x)$ tends to a point $a\in \partial_{\PP_G} G\setminus \La_P$ {after extracting a subsequence if necessary}. By \cite[Prop.\,1.1]{bowditch:convergence_groups}, 
we may assume that $(h_n)_n$ has the convergence property, i.e., tends to a point $b\in\La_P$  uniformly on compact subsets of $\partial_{\PP_G}G\setminus\{\al\}$
(since $(g_n(x))_n$ tends to $\al$). 
Since $\beta\ne \al$, it follows and $(h_ng_n(y))$ tends
to $b$ for all $y\ne x$. Thus the limits of $(h_ng_n(x))$ and $(h_ng_n(y))_m$ will be in different fibers for all $y\ne x$.
Therefore, we may also conclude that $z$ is conical.

In conclusion, we have defined a geometrically finite action on $Q_\PP$ with the prescribed maximal parabolic subgroups.
{As the action of $G$ is minimal on $\partial\cus (G,\PP_G)$, it is also the case in its quotient $Q_{\PP}$}. Thus $(G,\PP)$ is relatively hyperbolic and by definition, $Q_\PP$ is equivariantly homeomorphic to the boundary of the cusped space $\partial\cus (G,\PP)$.
\endp
 !TEX root =intro.tex

\section{Canonical splittings}\label{sec:cansplit}

We describe well-known splittings of manifolds and their counterparts for finitely presented  groups and show the quasi-isometric invariance of those splittings.
This provides the first step in the proof of Theorems \ref{thm:main1} and \ref{thm:main2} as described in \textsection \ref{outline}.
We start with the definition of a graph of groups structure 
and see how they appear when splitting manifolds. 

 Let us recall that our main results are concerned with groups up to quasi-isometry, so up to finite index. 
In particular, since any manifold admits a finite
covering of degree at most two that is orientable, we may ---and will always--- assume that our manifolds are orientable.

\subsection{Graph of groups and manifold splittings}     \label{scn:graph of group}

Let $G$ be a group. 
We first recall the definition of a graph of groups
and set up some notations. We follow Serre, and define a graph as a pair of sets $(V, E)$, with a fixed-point free
involution $e\to\bar e$ on $E$ (exchanging an oriented edge with its reverse orientation),
and a  terminal map $t:E\to V$ mapping an edge to the vertex it is oriented to.

A {\it graph of groups} $\GGG=(\G,\{G_v\},\{G_e\}, G_e\hookrightarrow G_{t(e)})$ is
\begin{enumerate}[- ]
\item a graph $\Gamma=(V,E)$;
\item an assignment 
of a group $G_e$ or $G_v$ to each edge $e$ or vertex $v$ of $\G$,
satisfying $G_{\bar e} = G_e$; and
\item for each edge $e \in E$, a monomorphism $G_e\hookrightarrow G_{t(e)}$.
\end{enumerate}

A graph of groups as above defines a group $G$ up to
isomorphism called
the fundamental group of the graph of groups  \cite[\textsection 5]{serre:arbre}. It
is characterized by an action on a simplicial tree $T$,
called the {\it Bass-Serre tree} of the graph
 of groups, with the following properties. The action {has} no edge inversions and the orbit space
$T/G$ is isomorphic to $\G$. 
Moreover, for any vertex $v\in T$, $\stab (v)$ is isomorphic to $G_{p(v)}$, for any edge $e\in T$, $\stab (e)$ is isomorphic to $G_{p(e)}$ and the canonical injection $\stab (e)\subset \stab (t(e))$ projects to the injection $G_{p(e)} \hookrightarrow G_{t(p(e))}$, where $p:T\to \G=T/G$ denotes
the canonical projection. 

A {\it graph of groups structure} for a group $G$ is a graph of groups $\GGG=(\G,\{G_v\},\{G_e\}, G_e\hookrightarrow G_{t(e)})$
together with an isomorphism between $G$ and the fundamental group of $\GGG$ as defined above.

The structure is {\it finite} if $\Gamma$ is finite and is {\it trivial} 
if there is a vertex group equal  to $G$. 
Unless otherwise stated, graph of groups structures will be assumed to be finite and non trivial.
{We will also only consider graph of groups structures of a group $G$ for which the action of $G$ on its Bass-Serre tree $T$ is {\it minimal}, i.e., 
$T$ contains no proper $G$-invariant subtree.}
We say that $G$ {\it splits} over a subgroup $H$ (which can be trivial) if there is a graph of groups structure for $G$ in which $H$ is an edge group.

Let us first notice that a splitting of a compact {three}-manifold $M$ produces a graph of groups structure for its fundamental group as follows, see \cite{scott:wall} for details.

Let $M$ be a compact {three}-manifold and let $S$ be a finite collection of disjoint non isotopic essential surfaces. We define the tree $\TTT_S$ dual to $S$ as follows: vertices are lifts of the components of $M\setminus S$ to $\tilde M$ and there is an edge between two vertices if the closures of the corresponding components intersect. The action of $\pi_1(M)$ on $\tilde M$ by covering transformations induces an action of $\pi_1(M)$ on $\TTT_S$ by isometries. This action provides $Q=\pi_1(M)$ with a graph of groups structure $\GGG=(\Delta,\{Q_v\},\{Q_e\}, Q_e\hookrightarrow Q_{t(e)})$: $\Delta=\TTT_S/Q$; edge groups are fundamental groups of components of $S$ and vertex groups are fundamental groups of components of $M\setminus S$.

The components of $M\setminus S$ are not compact and their closures might give a different decomposition. To remedy this we set up the following definition:

\begin{defn}[Submanifold compactification]    \label{defn:compactification}
Let $N$ be a component of $M\setminus S$ and $\tilde N$ a lift of $N$ to $\tilde M$. We define the {\em compactification} $\bar N$ of $N$ as the quotient of the closure of $\tilde N$ under the action of its stabilizer in $\pi_1(M)$ (which is isomorphic to $\pi_1(N)$).
\end{defn}

\subsection{Splittings over finite groups}  \label{splittings finite}

Let $M$ be a compact {three}-manifold and let $S\subset M$ be an essential surface. As explained above, the action of $\pi_1(M)$ on the dual tree 
$\TTT_S$ gives rise to a graph of groups structure {on} $\pi_1(M)$. If $S$ is a union of pairwise  non isotopic spheres and discs, 
then the edge groups are trivial and if furthermore $S$ is maximal, the vertex groups are fundamental groups of  irreducible and 
boundary irreducible {three}-manifolds.

With the next proposition, we show how to reverse this construction starting from a graph of groups with vertex groups corresponding to irreducible {three}-manifolds.

\begin{proposition} \label{prop:gluing discs}
Let $G$ be a group with a graph of groups structure $$\GGG=(\G,\{G_v\},\{G_e\}, G_e\hookrightarrow G_{t(e)})$$ with finite graph $\G$ and with finite edge groups. Assume that each vertex group $G_v$ has a finite index normal subgroup $G'_v$ isomorphic to the fundamental group of an irreducible orientable compact {three}-manifold $M_v$ (with or without boundary), then $G$ is commensurable to the fundamental group of a compact {three}-manifold $M$.

Furthermore, if $\partial M_v$ is non-empty for every vertex $v$, then $M$ is irreducible. {If for every $v$ either $M_v$ is atoroidal and $\chi(M_v)<0$ or $M_v$ is a ball or a solid torus or a thickened torus, then we may choose $M$ irreducible and atoroidal.} 
\end{proposition}

\Proof
If $G_v$ is finite, we may as well assume that $G'_v$ is trivial. Otherwise, since {$M_v$} is orientable and irreducible, $G'_v$ is torsion free for any vertex $v$ {(see \cite[Corollary 3.4]{boileau:around:3mfd} or \cite[(C3), p.50]{aschenbrenner:friedl:wilton:book})}. 
For an edge $e=(v,w)$, since $G_e$ is finite, $G_v'\cap G_e$ and $G_w'\cap G_e$ are trivial. Consider the graph of groups 
$$\overline{\GGG}=(\G,\{\overline{G}_v,\{\overline{G}_e\}, \overline{G}_e\hookrightarrow \overline{G}_{t(e)}\})$$ where $\overline{G}_v= G_v/G_v'$ and $\overline{G}_e= \{1\}$.
Let $\overline{G}$ be the fundamental group of $\overline{\GGG}$, which is a finite graph of finite groups, hence is virtually free \cite[Thm. 7.3]{scott:wall} and
residually finite  \cite{stallings:graphs}.

The canonical projections $G_v\to\overline{G}_v$ define a projection $q:G\to\overline{G}$ such that for any vertex $v$, $\ker q\cap G_v=G_v'$.
Since $\overline{G}$ is residually finite, there is a finite index subgroup $\overline{K}$ which is disjoint from any non trivial element of any vertex group $\overline{G}_v$. 
 Let $\overline{Q}=\cap_{g\in \overline{G}} g\overline{K} g^{-1}$ be the normal core of $\overline{K}$  in $\overline{G}$ (which has finite index in $\overline{G}$) 
 and let $p:\overline{G}\to\overline{G}/\overline{Q}$ be the canonical projection. By construction, the kernel $Q$ of $p\circ q$ is a finite index normal subgroup of $G$ 
 such that $G_v\cap Q=G_v'$ for any vertex group $G_v$. It follows that $Q$ has a graph of groups structure 
 $(\G',\{G_v'\})$ with trivial edge groups.

For each edge $e=(v,w)$ of $\G'$ we proceed as follows: 
if $M_v$ and $M_w$ have non-empty boundaries we pick a disc on $\partial M_v$ and a disc on $\partial M_w$ and glue $M_v$ to $M_w$ along those discs 
(we may have $M_v=M_w$). If $M_v$ or $M_w$ has no boundary, we remove a ball from $M_v$ and a ball from $M_w$ and glue  
$M_v$ to $M_w$ along the resulting boundary spheres. We do this operation for every edge by choosing disjoint discs and balls. 
It is easy to deduce from van Kampen's theorem that the fundamental group of the resulting manifold $M$ has a graph of groups structure 
$(\G',\{G_v'\})$ with trivial edge groups \cite{scott:wall}. Hence $\pi_1(M)=Q$ and we are done.

Notice that if every manifold $M_v$ has non-empty boundary, then we have only glued along discs so $M$ is irreducible. {If $M_v$ has more than one boundary components then we pick all the discs corresponding to adjacent edges on the same component and if $\chi(\partial M)<0$, we choose a component with negative Euler characteristic. This will ensure the last part of the statement as follows.}

{Note that, by construction, every subgroup isomorphic to $\Z^2$ within a vertex group is conjugated to a subgroup of the fundamental group of a boundary torus. By Seifert--van Kampen Theorem, $\pi_1(M)$ is obtained by starting from $G'_v=\pi_1(M_v)$ for an arbitrary initial vertex $v$ and inductively doing an HNN extension relative to the trivial subgroup  or a free product with $G'_w$ for another vertex $w$. In a free product $G_1\ast G_2$ any subgroup isomorphic to $\Z^2$ lies in a conjugate of $G_1$ or $G_2$. Similarly, in an HNN extension relative to the trivial subgroup $G_1\ast\Z$ any subgroup isomorphic to $\Z^2$ lies in a conjugate of $G_1$. It follows that every subgroup of $\pi_1(M)$ isomorphic to $\Z^2$ lies in a conjugate of a vertex group $G'_v=\pi_1(M_v)$. Therefore, if each vertex manifolds $M_v$ is either atoroidal with $\chi(M_v)<0$ or homeomorphic to a ball, a solid torus or a thickened torus, we have constructed $M$ so that it is atoroidal.}
\endp

We already mentioned the fact that a compact {three}-manifold $M$ is irreducible and has incompressible boundary if and only if $\pi_1(M)$ is one-ended. Thus, when we decompose $M$ along a maximal union of essential spheres and discs, we get a graph of groups with trivial edge groups and one-ended or finite vertex groups.

A similar decomposition has been established by  Stallings  for finitely generated groups \cite{stallings:yale}.
If $G$ is  non-elementary and not one-ended, it  splits over a finite group, i.e., $G$ is the fundamental group of a graph of groups with finite edge groups. 
Such a graph of groups is called {\it terminal} if it is finite and the vertex groups are one-ended or finite. A group which has a terminal splitting is {\it accessible}. 
According to Dunwoody \cite{dunwoody:accessibility}, finitely presented groups are accessible.

Terminal graphs of groups are invariant under quasi-isometries \cite{papasoglu:whyte:splitf}, see below. 
Combining Proposition \ref{prop:gluing discs} and Theorem \ref{prop:qi-ds}, it suffices to prove Theorem \ref{thm:main1} and \ref{thm:main2} for one-ended groups. 

\begin{theorem}[Papasoglu \& Whyte \cite{papasoglu:whyte:splitf}]\label{prop:qi-ds} 
Let $G$ be an accessible group and let 
$$\GGG=(\G,\{G_v\},\{G_e\}, G_e\hookrightarrow G_{t(e)})$$ be a terminal graph
of groups decomposition of $G$. A group $G'$ is quasi-isometric to $G$ if and only if
it is also accessible and any terminal decomposition of $G'$ has the same
set of quasi-isometry types of one-ended factors and the same number of ends.\end{theorem}

\subsection{Torus decomposition}     \label{scn:torus decomposition}

The second splitting of $M$ that we will use is related to its characteristic torus decomposition. The torus decomposition together with the annulus decomposition (see \textsection \ref{scn:annulus decomposition}) compose the JSJ splitting of $M$. For more on the history of the next statement due to Johannson and Jaco-Shalen \cite{jaco:shalen:seifert,johannson:jsj}, see \cite[Theorem 3.4]{bonahon:3-var}. 

\begin{theorem}[Characteristic torus decomposition]	\label{thm:torus}
Let $M$ be an orientable compact irreducible {three}-manifold. Then, up to isotopy, there is a unique compact $2$-dimensional submanifold $T$ of $M$ 
with the following properties. 

\begin{enumerate}[(i)]
\item Every component of $T$ is an essential torus.
\item The compactification of every component of $M\setminus T$ either contains no essential torus or else admits a Seifert fibration.
\item Property (ii) fails when any component of $T$ is removed.
\end{enumerate} 
\end{theorem}

We call $T$ the {\it characteristic torus decomposition} of $M$. {It follows from the Torus theorem \cite[Theorem IV.4.1]{jaco:shalen:seifert} that a compact Haken {three}-manifold that contains no essential embedded torus is either atoroidal or Seifert fibered. Thus {when $M$ is Haken,} we can replace (ii) with "the compactification of every component of $M\setminus T$ either is atoroidal or else admits a Seifert fibration".} 
It follows from the works of Thurston and Perel'man that $T$ is also a geometric decomposition {and that non-Haken {three}-manifolds are geometric.}
Conversely geometric manifolds have empty characteristic torus decompositions 
except the quotients of $Sol$ \cite[Thm\,5.3]{scott:geo3}. Thus the uniqueness statement in Theorem \ref{thm:torus} implies the uniqueness 
(up to isotopy) of the geometric decomposition if we add a minimality assumption in the spirit of property {\em (iii).}

As explained at the beginning of \textsection \ref{sec:cansplit}, $T$ induces on $\pi_1(M)=Q$ a structure of graph of groups 
$\GGG=(\Delta,\{Q_v\},\{Q_e\}, Q_e\hookrightarrow Q_{t(e)})$ where edge groups are isomorphic to $\Z^2$ and vertex groups 
are fundamental groups of components of $M\setminus T$.

We would like to have a similar splitting for a group quasi-isometric to $\pi_1(M)$. Different versions of the
JSJ splitting have been given for finitely generated groups (see \cite{guirardel:levitt:jsjpanorama} for a survey) but they tend to be different from the one given 
by the torus decomposition above since $\pi_1(M)$ may also split over cyclic subgroups. Rather than using this general theory, we will follow \cite{kapovich:leeb:qi3man} and use the quasi-isometry to construct a splitting of $G$. For technical purposes, we define the {\em balanced torus decomposition $T_B$} by replacing each component of $T$ 
by two disjoint parallel copies of itself {and by adding a torus parallel to each torus of $\partial M$}, thus adding a {thickened} torus between 
adjacent components of $M\setminus T$ {and next to each torus of $\partial M$}.

\begin{theorem}[Quasi-isometric invariance of the torus decomposition]		
\label{prop:qi torus}
Let $G$ be a finitely generated group quasi-isometric to the fundamental group of a non-geometric irreducible $\partial$-irreducible orientable {three}-manifold $M$. Let {$T_B$ be the  balanced characteristic torus decomposition of $M$}. 

The group $G$ has a graph of groups structure $\GGG=(\G,\{G_v\},\{G_e\}, j_e: G_e\hookrightarrow G_{t(e)})$ and there is a map $W:\G^{(0)}\to\{\mbox{Components of } M\setminus {T_B} \}$ such that {for each vertex $v$, the group $G_v$ is quasi-isometric to $\pi_1(W(v))$ and}
\begin{enumerate}[-]
{
\item if $W(v)$ is a thickened torus, then $v$ has one or two adjacent edges and for each adjacent edge $e$, 
$j_e(G_e)$ has index $1$ and $2$ in $G_v$};
\item if $W(v)$ is not a thickened torus, then {the {adorned group} $(G_v,\{G_e,\ t(e)=v\})$} is quasi-isometric to $(\pi_1(W(v)),P(v))$ where $P(v)$ denotes the collection of the fundamental groups of the tori of $\partial\bar W(v)$;
\item
{if} two vertices $v,w$ are adjacent {then} $W(v)$ and $W(w)$ are adjacent.
\end{enumerate}
\end{theorem}

When $M$ has zero Euler characteristic this result follows from the work of Kapovich and Leeb \cite{kapovich:leeb:qi3man} and Theorem \ref{thm:quotientman}.
In the general case the conclusion can still be deduced from the arguments of \cite{kapovich:leeb:qi3man} and Theorem \ref{thm:quotientman} as we will now explain. {Let us finally recall that $(\pi_1(W(v)),P(v))$ is a (non-pared) adorned group when $W(v)$ is a Seifert piece and a pared group when it is hyperbolic.}

The starting point of Kapovich and Leeb's proof is that $M$ is nonpositively  curved {on the large scale}  \cite[Theorem 1.1]{kapovich:leeb:npc3man}. More precisely, if $M$ has zero Euler characteristic, there exist a nonpositively curved  compact {three}-manifold $N$ and a  bi-Lipschitz homeomorphism between
their  universal covers $\tilde M$ and $\tilde N$ that  preserves their torus decompositions.

For manifolds with negative Euler characteristic, we can start with a stronger statement, using the arguments of \cite[Theorems 3.2 and 3.3]{leeb:nonpositive}. 

\begin{proposition}\label{nonpositive}
Let $M$ be an irreducible {three}-manifold with non-empty boundary. Then $M$ admits a Riemannian metric with nonpositive curvature {such that every torus on the boundary is flat}.
{If $M$ is non-geometric, then we may assume that the characteristic tori are flat.}\end{proposition}

\Proof We will recall some of the arguments of the proof of \cite[Theorem 3.3]{leeb:nonpositive} to show that it easily extends to manifolds with negative Euler characteristic.

{Since $\partial M\neq\emptyset$ holds, $M$ is a Haken manifold \cite[Chap. 13]{hempel:three_manifolds}. Therefore, it follows from Theorem \ref{thm:torus} and its following paragraph that if} the characteristic torus decomposition is empty then $M$ is either hyperbolic or admits a Seifert fibration. 
In the latter case, $M$ has a finite cover $M'$ which is a circle bundle over a compact surface $F$ (see \cite[Theorem 12.2]{hempel:three_manifolds}) such that $M'$ and $F$ are orientable. Since {$M$ has non empty boundary, this is also the case for $F$, so that} $M'$ is a trivial bundle and it admits an  $\HH^2\times\E^1$  or $\E^3$ structure which projects to $M$. It follows that $M$ is modeled on $\HH^3$, $\E^3$ or $\HH^2\times \E^1$ and the conclusion is obvious.}

{If the characteristic decomposition $T$, and hence the balanced decomposition $T_B$, of $M$ is not empty, we denote}
 by $M_H$ the union of the atoroidal components of $M\setminus T_B$, {by $M_E$ the union of the Euclidean components and 
by 
$M_G$ the closure of $M\setminus \overline{M_H\cup M_E}$.} By construction, $M_G$ is a graph manifold 
with non-empty boundary and, by (the proof of) \cite[Theorem 3.2]{leeb:nonpositive}, 
$M_G$ admits a Riemannian metric with nonpositive curvature and with flat boundary. By the hyperbolization theorem, the interior of $M_H$ admits a complete hyperbolic metric. 
Now, notice that the proof of \cite[Proposition 2.3]{leeb:nonpositive} consists in changing the metric of hyperbolic {three}-manifolds only in the rank $2$ cusps. 
Following exactly the same proof,  we get that any flat metric on $T$ can be extended to a nonpositively curved metric on $M_H$. 
In particular, we can extend the nonpositively curved metric on $M_G$ to a metric on $M$  with nonpositive curvature.
\endp

{To simplify later arguments we will list some properties of the metric on $M$ that follow easily from the construction. The characteristic tori $T_B$ are flat and totally geodesic. The metric on the components of $M\setminus T_B$ of Euclidean type is actually Euclidean. The metric on the other components is locally isometric to $\HH^3$ or $\HH^2\times \E^1$ outside a small neighborhood of the characteristic tori. The boundary $\partial M$ can be decomposed into a union $\partial_{\chi=0} M$ of flat totally geodesic tori and a union $\partial_{\chi<0} M$ of convex surfaces with negative curvature. We denote by $\tilde T_B\subset {\tilde M}$ the preimage of $T_B$ under the covering projection and call {\it geometric component} the closure of a component of $\tilde M\setminus\tilde T_B$.} 

{We now explore the action of a self-quasi-isometry of $\tilde{M}$ following the work of Kapovich and Leeb.}

\begin{lemma}	\label{quasi-invariant}
{Any quasi-isometry $\vp: \tilde{M}\to\tilde{M}$ preserves the balanced decomposition up to finite Hausdorff distance, i.e., any geometric component of $\tilde M$ is
mapped within bounded Hausdorff distance to another geometric component of the same type -- Euclidean, Seifert, hyperbolic whose quotient has zero Euler characteristic or hyperbolic whose quotient has negative Euler characteristic -- and adjacent components are mapped to adjacent components (up to bounded Hausdorff distance).} 
\end{lemma}

\Proof
{
First, we define an equivalence relation on the flats of $\tilde M$ as follows: $F_1\sim F_2$ if and only $F_2$ lies in the $r$-neighborhood of $F_1$ for some $r>0$. By \cite[Theorem 4.6]{kapovich:leeb:qi3man}, a quasi-isometry $\varphi:\tilde M\to\tilde M$ induces a bijection $\Phi$ between equivalence classes. Following \cite{kapovich:leeb:qi3man}, a flat is {\em isolated} if it is a component of $\tilde T_B\cup\widetilde{\partial_{\chi=0} M}$, the preimage of $T_B\cup\partial_{\chi=0} M$ under the covering projection $\tilde M\to M$. As observed in \cite{kapovich:leeb:qi3man}, a flat is equivalent to an isolated flat if and only if no other flat intersects it transversally. 
By \cite[\textsection 3]{kapovich:leeb:qi3man}, $\varphi$ maps an equivalence class of isolated flats to a bounded neighborhood of an equivalence class of isolated flats. By construction, a connected component $A$ of $\tilde M\setminus \tilde T_B$ is Euclidean if and only if it is an equivalence class of isolated flats. It follows that the union of the Euclidean components of $\tilde M\setminus\tilde T_B$ is invariant under quasi-isometry.

On the other hand a family of flats ${\mathcal F}$ is equivalent to the boundary of the closure of a non-Euclidean component of $\tilde M\setminus\tilde T_B$ if and only if it has the following properties:
\begin{enumerate}[-]
\item {${\mathcal F}$ is an infinite family of pairwise non-equivalent flats, and}
\item Given $F_1,F_2\in{\mathcal F}$, any isolated flat separating $F_1$ and $F_2$ is equivalent to $F_1$ or $F_2$, and
\item Given $F_1\in{\mathcal F}$ and $F_2$ an isolated flat not equivalent to any element of ${\mathcal F}$, there is $F_3\in{\mathcal F}$ such that either $F_1$ separates $F_2$ and $F_3$ or $F_3$ separates $F_1$ and $F_2$.
\end{enumerate}
By \cite[Lemma 4.7]{kapovich:leeb:qi3man}, these two properties are invariant under quasi-isometries. Thus we get:

{\noindent\bf Claim.---} A family of flats bounding a non-Euclidean component of $\tilde M\setminus \tilde T_B$ are mapped by a quasi-isometry to a family of quasi-flats at bounded distance from a family of flats bounding a non-Euclidean component of $\tilde M\setminus \tilde T_B$. 

Since $M$ is compact, the closure {$\bar A$} of any component {$A$} of $\tilde M\setminus \tilde T_B$ lies in a bounded neighborhood of the union of its boundary flats. {It follows then from the above claim that any non-Euclidean geometric component is mapped by a quasi-isometry within bounded distance to another non-Euclidean geometric component.

Any non-Euclidean component {$A$} of $\tilde M\setminus \tilde T_B$ is biLipschitz equivalent to a manifold $\tilde V$ where $\tilde V$ is either a product $\tilde F\times \R$ of the universal cover $\tilde F$ of a compact hyperbolic surface with geodesic boundary and $\R$ equipped with the product metric or the universal cover of a compact hyperbolic $3$-manifold whose toric boundary components are flat. In particular the Hausdorff distance between two different components of $\partial \tilde V$ is infinite.
 
From the two previous paragraphs, we get that any non-Euclidean geometric component is mapped by a quasi-isometry within bounded distance to exactly one non-Euclidean geometric component.
 
From the definition of the balanced torus decomposition we get that two components {$A$} and {$A'$} of $\tilde M\setminus \tilde T_B$ are adjacent if and only if one of them, say {$A'$}, is Euclidean and contains a flat equivalent to a flat in {$\partial A$}. From the previous discussion, we deduce then that adjacency is being preserved by quasi-isometries.}

{By \cite[Lemma 3.2]{kapovich:leeb:qi3man}, the universal cover of a hyperbolic piece can not be quasi-isometric to the universal cover of a Seifert piece.}
We can infer that each geometric piece is mapped to another geometric piece of the same type up to bounded Hausdorff distance. In particular, hyperbolic pieces are mapped one another. {Given a geometrically finite representation $\rho:\pi_1(W)\to\PP SL_2(\C)$ uniformizing a compact orientable hyperbolic manifold $W$, $\chi(W)=0$ if and only if $\HH^3/\rho(\pi_1(W))$ has finite volume or equivalently if and only if its limit set $\La_\rho$ is the whole sphere at infinity. When $\rho(\pi_1(W))$ is minimally parabolic, by definition, $\La_\rho$ is the Bowditch boundary of $\pi_1(W)$ relative to its (conjugacy classes of) rank $2$ Abelian subgroups. It follows that vanishing Euler characteristic is a quasi-isometry invariant of compact hyperbolic $3$-manifolds and that hyperbolic pieces with negative Euler characteristic are preserved as well.}
\endp

To continue the proof of {Theorem} \ref{prop:qi torus}, we need to set up a few definitions.

Suppose that $G$ is a group and $\rho$ is a map from $G$ to the set of all $(K,\epsilon)$-quasi-isometries of a metric space $X$. We call $\rho$ a {\it quasi-action} of $G$ if for some constant $L$ and all $g_1,g_2\in G$ the quasi-isometries $\rho(g_1g_2)$ and $\rho(g_1)\circ\rho(g_2)$ are $L$-close. The quasi-action is called {\it quasi-transitive} if for some constant $M$ all orbits $\rho(G).x$ are $M$-close to $X$. The {\it kernel} of the action $\rho$ is the subgroup of $G$ which consists of elements whose action on $X$ is {at bounded distance from}  the identity.
A quasi-action is called {\it properly discontinuous} if for each bounded subset $C\subset X$ there are only finitely many elements $g_j\in G$ so that $\rho(g_j)(C)\cap C\neq\emptyset$.

We say that a collection ${\mathcal A}$ of subsets $A\subset X$ is {\em quasi-invariant} under the quasi-action $\rho$ if: 
\begin{enumerate}[-]
\item every bounded subset $B\subset X$ intersects only finitely many sets in ${\mathcal A}$;
\item any two distinct sets in ${\mathcal A}$ {are at} infinite Hausdorff distance;
\item there is a constant $H$ such that for all $g\in G$ and $A\in{\mathcal A}$ the set $\rho(g)(A)$ is $H$-Hausdorff close to another set in ${\mathcal A}$.
\end{enumerate}
We can define the {\it stabilizer} in $G$ of a set $A$ in a quasi-invariant collection ${\mathcal A}$. 
It consists of all elements $g\in G$ such that $\rho(g)(A)$ and $A$ {are at} finite Hausdorff distance 
{(this is relevant when the sets $A$ are unbounded as in our case of interest)}. 
Clearly the stabilizer is a subgroup of $G$ and it is easy to define its quasi-action on $A$. This quasi-action is properly discontinuous 
and quasi-transitive if the quasi-action of $G$ on $X$ has such properties, see \cite[Lemma 5.2]{kapovich:leeb:qi3man}.

\medskip

\demode{Theorem \ref{prop:qi torus}}
By Proposition \ref{nonpositive}, 
$M$ may be endowed with a metric of nonpositive curvature. The group $G$ is quasi-isometric to its universal cover $\tilde M$. Let $f_1:\tilde M\to G$ and $f_2:G\to\tilde M$ be {quasi-isometries} 
such that $d(f_2\circ f_1(x),x)\leq C$ and $d(f_1\circ f_2(y),y)\leq C$ for some $C\ge 0$ for any $x\in\tilde M$ and any $y\in G$. 
Given $g\in G$, it is easy to see that the map $\rho(g):\tilde M\to \tilde M$ defined by $\rho(g)(x)=f_2(g f_1(x))$ is a quasi-isometry and that $\rho$ is a quasi-transitive properly discontinuous quasi-action of $G$ on $\tilde M$. Let {$\TTT$} be the tree dual to {$\tilde T_B$}, as defined in \textsection \ref{scn:graph of group}. By Lemma \ref{quasi-invariant}, $\rho(g)$ induces an automorphism $\sigma(g)$ on $\TTT$ {with no edge inversion {(by the definition of the balanced torus decomposition, if two components {$A$} and {$A'$} of $M\setminus T_B$ are adjacent, then one is a thickened torus and the other one is either hyperbolic or Seifert)}.} 
Since $\rho$ is a quasi-action, {we obtain} $\sigma(g_1)\sigma(g_2)=\sigma(g_1g_2)$, and we get a simplicial action of $G$ on $\TTT$. The quotient $\G={\mathcal T}/G$ is finite and the action induces a graph of groups structure $\GGG=(\G,\{G_v\},\{G_e\}, G_e\hookrightarrow G_{t(e)})$ for $G$. 
{Let $v$ be a vertex of $\TTT/G$. It follows from the construction that the corresponding vertex group $G_v$ is the stabilizer of a component {$A(v)$} of $\tilde M\setminus\tilde T_B$ under the quasi-action.}

{
 Without loss of generality, we may assume that there is $x\in \tilde M$ such that  $f_1(x)$ is the neutral element in $G$. By definition, we have $\rho(g)(x)=f_2(gf_1(x))=f_2(g)$, in particular the image of $f_2$ is the orbit of $x$ under the quasi-action of $G$ on $\tilde M$. Since the quasi-action of $G_v$ on {$A(v)$} is quasi-transitive, \cite[Lemma 5.2]{kapovich:leeb:qi3man}, its image $f_2(G_v)$ is at bounded Hausdorff distance from {$A(v)$}. It follows from Fact \ref{fact:close qi} that $G_v$ is quasi-isometric to {$A(v)$}, and hence to the fundamental group {$\pi_1(W(v))$} of its projection {$W(v)\subset M$}. 
 {Similarly, an edge group $G_e$ is quasi-isometric to an isolated flat and if $e$ is adjacent to $v$ then this flat is equivalent to a
 boundary component of {$A(v)$. Thus when $W(v)$ is a thickened torus, then $v$ has one or two adjacent edges and for each adjacent edge $e$, 
 {$j_e(G_e)$ is a finite index subgroup of $G_v$ with index $1$ or $2$}. When $W(v)$ is not a thickened torus, 
 the stabilizers of adjacent edge groups  $\{G_e,\ t(e)=v\}$ are also the stabilizers of the flat components of $\partial A(v)$ and} we get, from the quasi-transitivity of the quasi-action, a quasi-isometry between $(G_v,\{G_e,\ t(e)=v\})$ and $(\pi_1(W(v)),P(v))$ where $P(v)$ denotes the collection of the fundamental groups of the tori of $\partial\bar W(v)$ 
$(W(v),\partial_{\chi =0}W(v))$. {As noticed in the proof of Lemma \ref{quasi-invariant} the Hausdorff distance between two flat components of $A(v)$ is infinite, hence $\{G_e,\ t(e)=v\}$ is a non coarsely nested collection.} By construction,}
 } 
 when two vertices $v,w$ are adjacent then {$W(v)$} and {$W(w)$} are adjacent. This ends the proof of Theorem \ref{prop:qi torus}. \endp

 We conclude this section by showing that the kernel of the quasi-action $\rho$ defined in the latter proof is finite.

\begin{fact}
 The kernel $\ker \rho$ of the quasi-action of $G$ on $\tilde{M}$  is finite.
\end{fact}

\Proof
 It follows from above that any element of the kernel fixes $\TTT$ pointwise, so 
 the kernel is in the intersection of all vertex stabilizers. 

If $v\in\TTT$ is a vertex of hyperbolic type, then the restriction of the quasi-action to  $G_v$ on   the corresponding component $X_{i(v)}$ of $\tilde{M}\setminus T_B$ extends to a genuine action on its Gromov boundary $\partial X_{i(v)}$ as a convergence group action \cite[Cor.\,3J, Theorem 3E]{tukia:convergence_groups}. Note that its Gromov boundary 
 $\partial X_{i(v)}$ is infinite, so, as the action of $\ker \rho$ is trivial on $\partial X_{i(v)}$, it must be finite since the quasi-action of $G$ is properly discontinuous.
 
 If there are no vertex of hyperbolic type, then, since $M$ is non-geometric, we may find three consecutive vertices $u,v,w\in\TTT$ 
 with $u,w$ corresponding to Seifert pieces and $v$ to a common thickened boundary flat. 
 It follows from  \cite[Prop.\,5.4]{kapovich:leeb:qi3man} and its subsequent paragraphs 
that the kernel stabilizes the fibers up to bounded Hausdorff distance of both $X_{i(u)}$ and $X_{i(w)}$. By definition of the torus decomposition, 
these fibers intersect transversely on $X_{i(v)}$. Therefore, the quasi-orbit of a point under  the kernel is also bounded on the latter, 
hence the kernel is finite.
\endp
}

\subsection{Annulus decomposition}  \label{scn:annulus decomposition}
Lastly, we will introduce the JSJ decomposition along annuli for atoroidal {three}-manifolds and a generalization to relatively hyperbolic groups. This will lead to splittings of groups quasi-isometric to Kleinian groups with the following properties.

\begin{theorem}[Geometric decomposition of quasi-Kleinian one-ended groups] \label{thm:dec JSJ}
Let $G$ be a finitely generated one-ended group quasi-isometric to a Kleinian group $K$. Then $G$ {is hyperbolic relative to some collection of subgroups $\PP$ made up of virtually rank $2$ Abelian subgroups and} {it admits} a graph of groups structure $\GGG=(\G,\{G_v\},\{G_e\}, G_e\hookrightarrow G_{t(e)})$ with the following properties:
\ben
\item edge groups are virtually cyclic;
\item edge groups incident to a {virtually} Abelian vertex group are all commensurable;
\item {virtually} Abelian vertex groups are virtually cyclic or virtually $\Z^2$;
\item {virtually} Abelian vertex groups are not adjacent to each other nor to themselves;
\item for every {non virtually} Abelian vertex group $G_v$, there is a pared compact hyperbolic {three}-manifold { ---which is either acylindrical or a pared $I$-bundle---} with pared fundamental group $(H_v,\PP_v)$ and a quasi-isometry between $(H_v,\PP_v)$ and $G_v$ equipped with the paring provided by adjacent edges and parabolic subgroups.
{\item every virtually rank $2$ Abelian subgroup of $G$ is contained in a conjugate of a vertex group;
\item if a maximal parabolic subgroup is not a vertex group, then it does not contain any edge group.}
\een
\end{theorem}

{By Remark \ref{torsion free etc}, we may assume that $K$ is torsion free, geometrically finite and minimally parabolic.}
To prove that proposition, we will build a splitting of $K$ in two different ways. 
On the one hand, we will introduce the second part of the JSJ decomposition of  the Kleinian manifold $M_K$: the characteristic annulus decomposition. This decomposition is obtained by doubling 
the manifold along its boundary and considering the restriction to the initial manifold of the torus decomposition of the double. Thus we get a family $A$ of annuli that cuts the manifold into pieces with some specific topological properties. As explained at the beginning of \textsection \ref{scn:graph of group} the action of $K=\pi_1(M_K)$ on the dual tree $\TTT_A$ to $A$ induces a graph of groups structure for $K$. {For practical purposes, we will refine $A$ into the finer {\it balanced annulus decomposition} $B$ 
with its dual tree $\TTT_B$ (see definition below).}

On the other hand, following \cite{bowditch:jsj} and \cite{papasoglu:swenson:continua}, we will build a tree $\TTT_{\La_K}$ 
from the topological features of the limit set $\La_K$. We will show that $\TTT_{\La_K}$ is isomorphic to {$\TTT_B$ (Cor. \ref{cor:commonjsj})} and that the action of $K$ on $\La_K$ induces the action of $K$ on {$\TTT_B=\TTT_{\La_K}$}. We will then see that the quasi-isometry between $G$ and $K$ extends to a homeomorphism between their Bowditch boundaries. This will yield an action of $G$ on $\TTT_{\La_K}$ and hence a splitting of $G$. Finally the quasi-isometric invariance of those splittings will provide us with the desired properties {\em (1)} to {(7)}.

\subsubsection{Annulus decomposition and JSJ tree}	\label{cut tree}

We start by introducing the characteristic annulus decomposition, which, together with the torus decomposition described 
in Section \ref{scn:torus decomposition}, form the JSJ decomposition defined by Johannson-Jaco-Shalen \cite{johannson:jsj,jaco:shalen:seifert}, 
see also \cite[Theorem 3.8]{bonahon:3-var}. As explained above it can be defined using the torus decomposition of the manifold obtained by doubling $M$ along its boundary, but it will be convenient to have a more straightforward definition.

\begin{theorem}[Characteristic annulus decomposition]	\label{thm:charac annulus}
Let $(M,P)$ be an orientable compact atoroidal boundary irreducible pared {three}-manifold. Then, up to isotopy, there is a unique compact $2$-dimensional submanifold $(A,\partial A)\subset (M,\partial M\setminus P)$ of $M$ such that:
\begin{enumerate}[(i)]
\item Every component of $A$ is an essential annulus.
\item The compactification of every component of $M\setminus A$, is either pared acylindrical or a pared $I$-bundle or a solid torus or a thickened torus.
\item Property (ii) fails when any component of $A$ is removed.
\end{enumerate}
\end{theorem}

We call $A$ the {\it characteristic annulus decomposition} of $M$. {Given a component $W$ of $M\setminus A$, with compactification $\overline W$, we set $P_W=(P\cap W)\cup (\partial\overline{W}\setminus (W\cap\partial M))$. When $W$ is not a solid or thickened torus, item (ii) of Theorem \ref{thm:charac annulus} says that {$(W,P_W)$} 
is a pared manifold (see definition in \textsection \ref{3 mfds 1}).}
A {\em pared $I$-bundle} is a pared manifold $(N,P)$ such that $N$ is homeomorphic to a product $F\times I$ over a compact surface $F$ by a homeomorphism that maps $P$ into $\partial F\times I$.

The action of $\pi_1(M)$ on the dual tree to $A$ induces a graph of groups structure for $\pi_1(M)$.
As previously mentioned, we will give an alternate construction of this action when $M$ is uniformized by a Kleinian group. To simplify the identification of the two constructions we add solid tori to $M\setminus A$ so that a component that is pared acylindrical or a pared $I$-bundle is only adjacent to solid tori and thickened tori. Concretely we replace each component $A_i$ of $A$ that does not lie in the closure of any solid torus or thickened torus component of $M\setminus A$ by two disjoint parallel copies of itself, thus adding a solid torus to $M\setminus A$. We call the resulting surface $B$ the {\it balanced annulus decomposition of $M$}, it has the following property:

\noindent {\it (ii') Every component of $B$ lies in the closure of a component of $M\setminus B$ whose compactification is a solid torus or a thickened torus}.

Notice that if we set 

\noindent {\it (iii') Property (ii) or (ii') fails when any component of $B$ is removed.} 

\noindent then the balanced annulus decomposition is uniquely defined (up to isotopy) by properties (i), (ii), (ii') and (iii').

To give an alternate definition of $\TTT_B$ using the limit set $\La_K$ of a Kleinian group $K$ uniformizing {$(M,P)$}, we will now, following \cite{papasoglu:swenson:continua}, define the subsets of $\La_K$ that will be used as vertices. Notice that since {$(M,P)$} is assumed to be boundary irreducible, $\La_K$ is connected {by Fact \ref{fact:lambda connected}}.

Given a continuum (i.e., a connected compact set) $X$, a point $x\in X$ is a {\it cut point} if $X\setminus \{x\}$ is not connected. We define an equivalence relation $\RRR$ on $\Lambda_K$. Each cut point is only equivalent to itself and if $a,b\in\Lambda_K$ are not cut points we say that $a\RRR b$ if they are not separated by any cut point, i.e. for any cut point $c\in \Lambda_K$, $a$ and $b$ lie in the same component of $\Lambda_K\setminus\{c\}$. By \cite[Lemma 32]{papasoglu:swenson:continua}, when $X$ is a Peano (i.e., locally connected) continuum the closure $Y$ of each non singleton equivalence class is a Peano continuum  without cut points.

Let $Y\subset X$ be the closure of a non singleton equivalence class for $\RRR$. A pair $\{a, b\} \subset Y$  is a {\it cut pair} if $Y \setminus \{a, b\}$ is not connected. A nonempty subset $A\subset Y$ {with at least two points} is called {\it inseparable} if no two points of $A$ lie in different components of the complement
of any cut pair. Every inseparable set is contained in a maximal inseparable set. 

A finite subset $S$ of a $Y$ is called a {\it cyclic subset} if either $S$ is a cut pair or there is an ordering $S  =\{s_j,\ j\in\Z/n\Z\}$, $n\ge 3$, and continua $M_j\subset Y$, $j\in\Z/n\Z$, such that
\begin{enumerate}[- ]
\item $M_i \cap M_{i+1} = \{s_i\}$, $i\in\Z/n\Z$,
\item $M_i \cap M_j = \emptyset$ whenever $|i - j| > 1$,
\item $\bigcup M_i = Y$.
\end{enumerate}

An infinite subset in which all finite subsets of cardinality at least $2$ are cyclic is also called
{\it cyclic}. A maximal cyclic subset with at least 3 elements is called a {\it necklace}.

In \cite{papasoglu:swenson:continua}, the authors use these subsets to construct an $\R$-tree associated to a Peano  continuum called the {\em combined tree}. We introduce a simplified version of this construction when the continuum is the limit set $\Lambda_K$ of a geometrically finite Kleinian group $K$.
Consider the set $V$ of cut points, {inseparable} cut pairs, necklaces and maximal inseparable sets of $\Lambda_K$ ({inseparable} cut pairs, necklaces and maximal inseparable sets are taken in the closures of non singleton equivalence classes of the relation $\RRR$ defined above). 
As we will see in Propositions \ref{prop:cutaccidental} and \ref{prop:annulus tree}, this set $V$ is countable. 
We define a graph $\TTT_{\Lambda_K}$  with vertex set $V$ by putting an edge between two vertices $v_1,v_2$ if $v_1\subset v_2$ as subsets of $\Lambda_K$. Rather than showing that this graph is a simplicial tree, we will directly show that {when $K$ is torsion-free,} it is isomorphic to the dual tree $\mathcal{T}_B$ to the balanced annulus decomposition of the Kleinian manifold $M_K$.

We conclude this subsection with a remark for the reader familiar with the work of \cite{papasoglu:swenson:continua}. We will see in Proposition \ref{prop:annulus tree} that the closure of a necklace intersects a maximal inseparable set only along an inseparable cut pair. Using \cite[Lemma 8 and Lemma 28] {papasoglu:swenson:continua} and the definition of $T$ in \cite[p. 1765]{papasoglu:swenson:continua}, one can see that $\TTT_{\La_K}$ is the {\em combined tree} associated to $\La_K$ defined on \cite[p. 1782]{papasoglu:swenson:continua}.

In the next sections we will describe the stabilizers of the sets constituting $V$ and their relation with the characteristic annulus decomposition.

\subsubsection{Cut points stabilizers} \label{scn:cutstab}
In this section we establish some properties of stabilizers of cut points. They will be used to prove that $\TTT_{\La_K}$ is isomorphic to $\TTT_B$ as well as properties {\em (1) and (2)} of Theorem \ref{thm:dec JSJ}. For the latter we will work in the general situation of relatively hyperbolic groups acting on their Bowditch boundary. {Let $(G,\PP)$ be a relatively hyperbolic group and assume that its Bowditch boundary $\partial_\PP G$ is connected. Then it
also is locally connected \cite[Thm 1.1]{dasgupta:hruska:jt2024}. Therefore cut points are bounded parabolic points  \cite[Thm 0.2]{bowditch:connectedness}, see also  \cite[Thm 1.1]{dasgupta:hruska:jt2024}. We study the dynamics of these parabolic subgroups intrinsically.}

\begin{proposition}\label{prop:cutpoint}
{Let $Y$ be a locally connected continuum, $p\in Y$ a cut point and assume that $H$ is a bounded parabolic convergence group fixing $p$.}
Let us denote by $\CCC$ the collection of connected components
of $X=Y\setminus\{p\}$. Then the action of $H$ on $\CCC$
has finitely many orbits and, for any $C\in \CCC$, the action of $\stab_HC$ is cocompact on $C$, hence {$\stab_HC$ is} infinite. Moreover, any compact subset of
$X$ meets at most finitely many components of $X$.
\end{proposition}

\Proof 
{Let $L$ be a compact subset whose $H$-translates cover $X$}.  Every point in $L$ admits a connected neighborhood in $X$, hence contained in a unique
element of $\CCC$. By compactness, this implies that only finitely many components of $\CCC$ intersect $L$. This shows that $\CCC$ is a finite union of
$H$-orbits. This also implies that   any compact subset $K$  of
$X$ meets at most finitely many components of $X$, since the fact that the action is properly discontinuous on $X$ implies that there are only finitely many
elements of $H$ that map points of $K$ in $L$. 

Let us denote by $C_1,\ldots, C_n\in\CCC$ the components that intersect $L$. If $C_i$ and $C_j$ are in the same orbit, we fix $h_{ij}\in H$ such that
$h_{ij}(C_j)=C_i$, with $h_{ii}=id$.  Given $C=C_i$, write $L_i=  \cup h_{ij}(C_j\cap L) $ where the union is taken over the components $C_j$ in the same $H$-orbit {as} $C_i$.
Note that $L_i$ is a compact subset of $C_i$, as a finite union of compact subsets of $C_i$. If $x\in C_i$, then we may find $g\in H$ such that $g(x)\in L$ by definition. If $g(x)\in C_j$, then
$h_{ij}(gx)\in L_i$. Since $h_{ij}g(C_i)=C_i$ by construction, we have proved that the action of $\stab_HC_i$ is cocompact on $C_i$.
Since $C_i$ is non-compact (the point $p$ is an accumulation point by definition), we may conclude that its stabilizer is infinite. 
\endp

A parabolic isometry $g$ in a {non-elementary} Kleinian group $K$ {with connected limit set} is an {\em accidental parabolic} if $g$ stabilizes a component $\OOO$ of the {ordinary set together with} a geodesic in $\OOO$ equipped with its hyperbolic metric {obtained from the {Riemann} uniformization theorem, since $\OOO$ is a proper open {simply connected} subset of the punctured Riemann sphere $\cbar$}.

\begin{proposition}\label{prop:cutaccidental} Let $K$ be a geometrically finite Kleinian group with connected limit set. 
A point $p\in\La_K$ is a cut point if and only if its stabilizer contains a primitive accidental parabolic.
\end{proposition}

\Proof {We first assume that the point $p$ is parabolic and that its stabilizer $H$ contains an accidental parabolic element $h\in K$. So let $\g\subset \Omega_K$ be a hyperbolic geodesic stabilized by the parabolic isometry $h$}.  Then $\g$ joins the fixed point $p$ of $h$ to itself and $\g\cup p$ is a Jordan curve separating $\La_K$. It follows that $p$ is a cut point.

If $h$ is not primitive, then there exists $g\in K$ and an iterate {$k> 1$} such that $g^k=h$.  
We first observe that $g\in H$ since $g^k$ fixes $p$. {Let $\OOO$ be the component of $\Omega_K$ containing $\g$. If $g$ does not stabilize $\OOO$, then $g(\gamma)\cap\gamma=\emptyset$ and since $g$ acts as a translation on the uniformization $\C$ of $\hat \C\setminus\{p\}$,  $g^i(\gamma)\cap\gamma=\emptyset$
{for all $i\ge 1$} contradicting the assumption $g^k=h$. Thus we have $g(\OOO)=\OOO$. Since $h$ stabilizes $\g$, it acts as a hyberbolic isometry (for the hyperbolic metric) on $\OOO$. Basic hyperbolic geometry tells us that $g$ must also act as a hyberbolic isometry on $\OOO$ with the same axis $\g$. In particular, $g$ is a primitive accidental parabolic.}

{We now establish the converse and suppose that $p$ is a cut point.} By Selberg's lemma, we may assume that $K$ is torsion free. Since $p$ is a cut point, one can find two disjoint closed sets $A$ and $B$ of $X= \La_K\setminus\{p\}$ which covers $X$. Note that,
on $\cbar$,  $\overline{A}\cap\overline{B}=\{p\}$. Therefore, by the separation theorem \cite[Thm VI.3.1]{whyburn:analytic_topology}, we may find
a Jordan curve $c$ that separates a point $a\in A$ from a point $b\in B$ such that $c\cap \La_K= \{p\}$.   

{As mentioned at the beginning of \S\,\ref{scn:cutstab}, it follows from  \cite{dasgupta:hruska:jt2024} and  \cite{bowditch:connectedness} that} the point  $p$ is parabolic; denote by $H$ its stabilizer (which is isomorphic to $\Z$ or $\Z^2$).

Let $\OOO$ be the component of $\Omega_K$ that contains the connected set $c\setminus\{p\}$ and $K_\OOO$ be its stabilizer.
Ahlfors finiteness theorem implies that
$\OOO/K_\OOO$ is a surface of finite area. Moreover, Thurston proved that $K_\OOO$ is geometrically finite as well, since it is finitely generated
\cite[Theorem 7.1]{morgan:thurston}. 
Since $p$ is not conical for $K$,
it cannot be conical for $K_\OOO$ either, so it is parabolic.

{Since $\La_K$ is connected, $\OOO$ is simply connected so
the Riemann mapping theorem provides us with a biholomorphic map $f:\DD^2\to\OOO$ from the unit disc to $\OOO$ and the action of $K_\OOO$ on $\OOO$ defines a representation $\sigma:K_{\OOO}\to \PP SL_2(\R)$.  This implies that $K_\OOO\cap H$ is a cyclic group.
By \cite[Cor. 4.2]{anderson:maskit}, $\partial \OOO$ is locally connected and by the Carath\'eodory-Torhorst theorem \cite[Theorem 5.5]{conway:functions:II}, $f$ extends continuously to $\partial\DD^2$ (see also \cite[Corollary p. 215]{floyd:completion}).
The curve $c$ disconnects $\partial\OOO$ since, otherwise, $c$ would bound a disc in $\OOO$, and, hence, would not separate $A$ and $B$.
It follows that $\partial\DD^2\setminus\{f^{-1}(p)\}$ is not connected which is only possible if $\sharp\{f^{-1}(p)\}>1$. 
Then we get from \cite[Theorem p.207]{floyd:completion} that $\sharp\{f^{-1}(p)\}= 2$.
Therefore the hyperbolic geodesic $\ell$ joining the two points of $f^{-1}(p)$ is also invariant and $f(\ell)\subset\OOO$ 
is invariant under the action of $K_\OOO\cap H$.}
This implies that $p$ is an accidental parabolic.
\endp

 \begin{corollary}\label{cor:cutaccidental}  Let $K$ be a geometrically finite Kleinian group with connected limit set. 
Assume that {the point $p=\infty \in \hat\C$}   is a cut point with stabilizer isomorphic to $\Z^2$. Then $X=\Lambda_K\cap \C$ has infinitely many components {and these components all have a common stabilizer which is cyclic}. Moreover, given a component $C$ of $\La_K\cap \C$ and $R>0$, there are only finitely many other 
components at Euclidean 
distance at most $R$ from $C$.  \end{corollary}

\Proof  Let $H$ be the stabilizer of $p$. As $p$ is at infinity, {a parabolic isometry fixing $p$ acts on $\C$ by translations. Since $H$ is discrete and isomorphic to $\Z^2$, its action on $\C$ is cocompact.} 

{It follows from Proposition \ref{prop:cutaccidental} that there is a component $\OOO$ {of $\Omega_K$} that is fixed by a cyclic group generated by a primitive element $h_1$.}
{As $h_1$ is primitive, Bezout's theorem implies the existence of a primitive element $h_2$ such  that $H$ is generated by $h_1$ and $h_2$.}

Moreover, $h_2$ acts freely on the $H$-orbit of $\OOO$, and so on the components
of $\La_K\cap \C$ since the $H$-orbit of $\OOO$ splits the plane in parallel (topological) strips  {in the following sense. By Proposition \ref{prop:cutaccidental}, There is an $h_1$-invariant curve in $\OOO$ that is translated by all the iterates of $h_2$. These curves cut the plane into strips.}

Let $C$ be a component of  $\La_K\cap \C$, and assume that $h\in H$ fixes $C$. Then we may find $k,\ell\in\Z$ such that $h=h_1^{k}h_2^{\ell}$.
It follows from above that $\ell=0$ and $h$ is an iterate of $h_1$. This shows that $\stab_HC\subset H_\OOO$. {We claim that} $h_1(C)=C$. {To see this, consider the canonical projection $\pi:\C\to \C/(h_1^k)$ to the cylinder. The set $C$ projects into a connected compact subset that separates both ends of the cylinder. The map $h_1$ acts as an automorphism of order $k$ and maps $\pi(C)$ to $\pi(h_1(C))$ that has the same property: if $h_1(C)\ne C$, then $\pi(h_1(C)))\cap \pi(C)=\emty$, but then $h_1$ would not be of finite order. It follows that $h_1(C)=C$ as claimed.}
Thus, every component $C$ has stabilizer the subgroup generated by $h_1$, i.e., $H_\OOO$. 

{\noindent\bf Claim.---} The sequence $(h_2^n)_n$ is uniformly convergent to $p$ on any component $C$ of $X$. 

If not, there would be a sequence of points $(x_k)$ in $C$ such that $h_2^{n_k}(x_k)$ remains in a compact subset of $X$ for some subsequence $(n_k)$.
Since the action of $\stab_HC$ is cocompact on $C$ {by Proposition \ref{prop:cutpoint}}, we may assume that $(x_k)$ remains in a compact subset of $C$ {precomposing $h_2^{n_k}$ by a suitable iterate of $h_1$}. Since the action is
properly discontinuous, this implies that  an iterate of $h_2$ fixes $C$, which is absurd.  
This proves the claim.

Fix $R>0$ and a component $C$ of $\La_K\cap\C$. Fix another component $C'$, and let assume that there are $m$ translates of $C'$ at
distance at most $R$ from $C$. We may find 
$(n_k)$, $1\le k\le m$, such that $\dist (C,h_2^{n_k}C')\le R$ for all $k$. Since the actions of the stabilizers  of components are cocompact on their components {(cf. Proposition \ref{prop:cutpoint})},
we may assume that these distances are realized in some compact subset of $X$. {Together with the fact that there are finitely many $H$-orbits of $X$ by Proposition \ref{prop:cutpoint}, t}he claim implies that $m$  has a uniform upper bound.  
So there are only finitely many components at distance at most $R$ from $C$.\endp

\begin{lemma}\label{lma:paraz2} Let $K$ be a Kleinian group with connected limit set $\La_K$ and let $G$ be a  group virtually isomorphic to $\Z^2$, 
acting by uniformly quasi-M\"obius maps on $\La_K$. We assume 
that the action of $G$ is parabolic with common fixed  point a cut point $p\in\La_K$, that the action is cocompact on $\La_K\setminus\{p\}$ and that the stabilizer of $p$ in $K$ {is a rank 2 Abelian group}. Then every stabilizer in $G$ of a connected component of $\La_K\setminus\{p\}$  is virtually cyclic, and they are
pairwise commensurable. 

Moreover, there exists a finite index free Abelian subgroup $G_A$  generated by $g_1,g_2\in G$ such that $g_1$ stabilizes each connected component and generates a subgroup acting cocompactly on each component and $g_2$ acts freely on the components of $\La_K\setminus\{p\}$.
\end{lemma}

\Proof
We may assume that the point $p$ is at infinity in $\cbar$ so that {$G$ acts by proper {uniformly} quasi-M\"obius maps on $\La_K\setminus\{p\}\subset \C$ endowed with the Euclidean metric. Notice that as  $x_4$ tends to infinity, the crossratio $[x_1:x_2:x_3:x_4]$ will degenerate into $d(x_1,x_2)/d(x_1,x_3)$, showing that} $G$ acts by uniform {\em quasisymmetric} maps on $X= \C\cap\La_K$.
This means that there exists an increasing homeomorphism $\eta:\R_+\to\R_+$ such that, for any $x,y,z\in X$, {any $t>0$} and any $g\in G$,
$$\hbox{if }  |x-y|\le t|x-z| \hbox{ then  } |g(x)-g(y)|\le \eta(t) |g(x)-g(z)|\,.$$

Since $p$ is a cut point, $X$ has at least $2$ connected components. Let $H\subset K$ denote the stabilizer of $p$. By assumption, it
is a rank $2$ Abelian group acting properly discontinuously by translations on $\C$. 
By Corollary \ref{cor:cutaccidental} there are constants $\de_\pm >0$ such that, for any component $C$ of $X$, $\dist(C,X\setminus C)\ge \de_-$ and for any $x\in C$,
$\dist(x,X\setminus C)\le \de_+$ hold. 

We first prove that we may deform the Euclidean metric on $X$ quasisymmetrically so that $G$ acts by isometries.
Let $x,y\in X$. We claim that $d_{X}(x,y)= \sup\{|g(x)-g(y)|,\ g\in G\}$ is finite. Let $C_x$ be the component of $X$ containing $x$.
Given $g\in G$, pick $z\in X\setminus C_x$ so that $|g(x)- g(z)|= \dist(g(x),X\setminus g(C_x))$,
then we have $$|g(x)-g(y)| \le \eta\left(\frac{|x-y|}{|x-z|}\right) |g(x)-g(z)|\le  \eta\left(\frac{|x-y|}{\de_-}\right) \de_+\,.$$

This implies that $(X,d_{X})$ is a metric space on which $G$ acts by isometries. Moreover, for any $x,y,z\in X$ and $\ep>0$, if $g\in G$
satisfies $d_{X}(x,y)\le (1+\ep)|g(x)-g(y)|$, then
$$\frac{d_{X}(x,y)}{d_{X}(x,z)}\le (1+\ep) \frac{|g(x)-g(y)|}{|g(x)-g(z)|} \le (1+\ep) \eta\left(\frac{|x-y|}{|x-z|}\right)\,.$$
Since $\ep$ is arbitrary, it follows that $id:(X,|\cdot|)\to (X,d_{X})$ is $\eta$-quasisymmetric. 

Since the action of $G$ is {cocompact on $X$ and} parabolic, Proposition \ref{prop:cutpoint} 
implies that there are only finitely many orbits of such connected components. Moreover, their stabilizer are virtually cyclic since $G$ has rank two
and there are infinitely many components, cf. Corollary \ref{cor:cutaccidental} {applied to the Kleinian group $K$}. 

If $C'$ is another component, then, for any $g\in\stab_G C$,
$d_X(C,g(C'))={d_X}(C,C')$. 
{By the definition of $d_X$, the Euclidean distance $d_e(C,g(C'))$ between $C$ and $g(C')$ is bounded
from above by $d_X(C,C')$ : $d_e(C,g(C'))\le d_X(C,g(C'))\le d_X(C,C')$.} {It follows then} 
from Corollary \ref{cor:cutaccidental} that the orbit of $C'$ under  $\stab_G C$ is finite. This implies that $\stab_G C\cap \stab_G C'$ has finite
index in $\stab_G C$. Therefore, stabilizers are commensurable.

Let $G'$ be a rank two free Abelian subgroup of finite index in $G$. Then stabilizers of components in the same $G'$-orbit are the same and cyclic {since they are infinite}. They
also act cocompactly on each of them. Therefore, there are only finitely many such stabilizers and their intersection is 
a cyclic group generated by some $g_1$ of finite index in any
stabilizer. Let us consider a primitive $g_2\in G'$ so that it generates with $g_1$ a subgroup $G_A$ of rank two. The element $g_2$ acts freely
on the components of $X$ by construction.
\endp

\subsubsection{Stabilizers and characteristic annulus decomposition} \label{sec:stab annulus}
Next we will study the relations between the stabilizers of the sets that make up $V$ and the balanced annulus decomposition. We start with the cut points.

\begin{lemma}   \label{lma:cut point}
Let $K$ be a {torsion free} 
geometrically finite Kleinian group uniformizing a pared manifold $(M,P)$ with associated representation $\rho:\pi_1(M)\to K$. {Assume that $\Lambda_K$ is connected and} let $B$ be the balanced annulus decomposition of $(M,P)$. 
A point $c\in \Lambda_K$ is a cut point if and only if there is a component $P_i$ of $P$ and a solid torus or a thickened torus $W$ in $M\setminus B$ containing $P_i$ such that $c$ is the fixed point of a conjugate of $\rho(\pi_1(W))=\rho(\pi_1(P_i))$.
\end{lemma}

{Notice that, by Fact \ref{fact:lambda connected}, $\partial M\setminus P$ is incompressible and we can use Theorem \ref{thm:charac annulus} to define the characteristic annulus decomposition and then the balanced annulus decomposition $B$.}

\Proof
Let $c\in\La_K$ be a cut point. By Proposition \ref{prop:cutaccidental}, $c$ is stabilized by a parabolic subgroup $H<K$. 

{Let $\tilde\gamma\subset\OOO_K$ be the union of all the geodesics in the ordinary set $\OOO_K$ that are stabilized by an element of $H$. It follows from \cite[Theorem p.207]{floyd:completion}, that each component of $\OOO_K$ can contain at most one of those geodesics. Then $\tilde\gamma$ projects to a multicurve $\gamma\subset\partial M$. Since $K$ uniformizes $(M,P)$ there is a component $P_i$ of $P$ such that $H$ is conjugate to $\rho(\pi_1(P_i))$ hence each component of $\gamma$ is homotopic to a closed curve on $P_i$.
By the Annulus theorem \cite[IV.3.1]{jaco:shalen:seifert} there is an essential annulus $E_{\gamma_k}$ joining each component $\gamma_k$ of $\gamma$ to $P_i$.  {By Proposition \ref{prop:cutaccidental} and Corollary \ref{cor:cutaccidental}, there is a single primitive element $h \in H$ such that every component of $\tilde{\g}$ has
$H$-stabilizer equal to $\langle h\rangle$, hence all the components of $\gamma$ are homotopic to the same simple closed curve on $P_i$.} A simple surgery argument allows us to construct a family of disjoint annuli $E_{\gamma_k}$ joining the components of $\gamma$ to $P_i$. Let us denote by $\Delta_i\subset M$ the $2$-complex made up of these annuli and $P_i$ and by $V_{P_i}\subset M$ a regular closed neighborhood of $\Delta_i$.

The frontier $Fr(V_{P_i})$ of $V_{P_i}$ is a union of essential annuli. It is a fundamental property of the JSJ decomposition that any essential annulus is isotopic to an annulus that is disjoint from $B$, see \cite[Theorem 3.9]{bonahon:3-var}. It follows then from the minimality of the characteristic annulus decomposition that $V_{P_i}$ is contained in a component $W_{P_i}$ of $M\setminus B$ up to isotopy. Let us assume that $V_{P_i}$ is not isotopic to its compactification $\overline{W}_{P_i}$. Then there is a component $F$ of $Fr(V_{P_i})$ that is not peripheral in $\overline{W}_{P_i}$, in particular $(\overline{W}_{P_i},(B\cap\overline{W}_{P_i})\cup P_i)$ is not pared acylindrical. If $W_{P_i}$ is an $I$-bundle $S\times I$ then $P_i$ and $F$ are union of fibers, up to isotopy, by \cite[Corollary 3.2 and Lemma 3.4]{waldhausen:irreducible}. Since $F$ is homotopic to a curve in $P_i$ this would contradict the fact that $F$ is 
{boundary incompressible and not peripheral}. Thus we have proved that $W_{P_i}$ is a solid or thickened torus, and by construction, $c$ is the fixed point of a conjugate of {$\rho(\pi_1(W_{P_i}))$}. We can further deduce from the maximality of $\tilde\gamma$ that $\overline{W}_{P_i}$ and $V_{P_i}$ are isotopic.}

Conversely, let $W$ be a solid torus or thickened torus in $M\setminus B$ containing a component $P_i$ of $P$. 
Then $\rho(\pi_1(W))$ is a parabolic subgroup stabilizing a point $c\in\La_K$. Let $F$ be a component of $W\cap\partial M$ 
which is disjoint from $P$ and let $\gamma$ be a simple closed curve on $F$. 
The curve $\gamma$ lifts in $\Omega_K$ to an arc  $\hat \gamma$ joining $c$ to itself. 
By construction,  $\hat\gamma\cup c$ separates $\La_K$. It follows that $c$ is a cut point.
\endp

Let $W_P\subset M\setminus B$ be the union of the connected components $W_i$ of $M\setminus B$ such that $\rho(\pi_1(W_i))$ is a parabolic subgroup {and $W_i$ is either a thickened torus or a solid torus (the latter case happens only if $K$ is not minimally parabolic). It follows from the minimality of the balanced annulus decomposition (Property (iii')) that no two components of $W_P$ are adjacent.} Recall that if $a,b\in\Lambda_K$ are not cut points  then $a\RRR b$ if they are not separated by any cut point.

\begin{lemma}   \label{Y}
{We keep the assumptions of Lemma \ref{lma:cut point}.} A subset $Y\subset \Lambda_K$ is the closure of a non-singleton equivalence class of $\RRR$ if and only if there is a connected component $W$ of $M\setminus W_P$ such that $Y$ is the limit set of a conjugate of $\rho(\pi_1(W))$.
\end{lemma}

\Proof
{
Given a component $W_{P_i}$ of $W_P$, in the proof of Lemma \ref{lma:cut point}, we constructed a $2$-complex $\Delta_i$ so that $W_{P_i}$ is isotopic to a neighborhood of $\Delta_i$. Using the same notations as in the proof of Lemma \ref{lma:cut point}, we will introduce another way to build $\Delta_i\setminus P_i$ using the Kleinian group $K$. 

Each component {$\tilde\gamma_k$} of $\tilde\gamma$ is stabilized by a subgroup $H_{\tilde\gamma}<H<K$ isomorphic to $\Z$ and the closure {$\tilde\gamma_k\cup c$} of {$\tilde\gamma_k$} bounds in $\HH^3\cup\hat\C$
a closed  disc {$\tilde D_{\tilde \gamma_k}$} whose interior lies in $\HH^3$ and that corresponds to the closure of the universal covering of  {$E_{\gamma_k}$} in $\HH^3\cup\Omega_K\approx \tilde{M\setminus P}$. As
these annuli are pairwise disjoint, these discs are also}
pairwise disjoint and form a $K$-invariant collection.
We denote by {$\tilde D_{P_i}=\bigcup \tilde D_{\tilde \gamma_k}$}
the union of all these discs. The projection $D_{P_i}$ of $\tilde D_{P_i}\cap (\HH^3\cup \Omega_K)$ to $(\HH^3\cup \Omega_K)/K\approx M\setminus P$ is a family of disjoint half infinite annuli joining the components of $\gamma$ to the cusp associated to $P_i$. {By construction $D_{P_i}$ is equal to $\Delta_i\setminus P_i$ as promised.}

Doing the same construction for each component of $P$ we get that every component $V$ of $( (\HH^3\cup \Omega_K)/K)\setminus\bigcup_i D_{P_i}$ is isotopic to a component of $M\setminus\bigcup_i \Delta_i$ and hence to a component $W$ of $M\setminus W_P$ according to the proof of Lemma \ref{lma:cut point}. It follows that if {$\hat V\subset\HH^3\cup\cbar$} is the closure of a lift $\tilde V$ of a component $V$ of $( (\HH^3\cup \Omega_K)/K)\setminus\bigcup_i D_{P_i}$, then $\La_K\cap\hat V$ is the limit set  $\Lambda_{K_W}$ of  a conjugate $K_W$ of $\rho(\pi_1(W))$.}

Since {$\partial W\setminus P$} is incompressible, $\La_{K_W}$ is connected. {Since $B\cap W$ is the balanced decomposition of $(W,P\cap W)$, $\La_{K_W}$ does not have any cut point  by Lemma \ref{lma:cut point}.
Hence $\La_{K_W}$ lies in the closure of a non-singleton equivalence class of $\RRR$. Given a point $x\in\La_K\setminus \La_{K_W}$, 
any point $y$ in $\La_{K_W}$ which is not a cut point of $\La_K$ is  separated from $x$ in $\HH^3\cup \cbar$ by a component of $\tilde D$. Hence $x$ and $y$ are separated in $\cbar$ by {the closure of a component} of $\tilde\g$ and in $\La_K$ by a cut point. It follows that $\La_{K_W}$ is the closure of a nonsingleton equivalence class of $\RRR$.}

By Lemma \ref{lma:cut point}, any point in $\La_K$ that is not a cut point lies in $\cbar\setminus \bigcup_{\tilde c\in\tilde\gamma} \tilde c$. Thus we have obtained above all the nonsingleton equivalence classes of $\RRR$.
\endp

Let $(W,P_W)$ be the compactification of a component of $M\setminus W_P$ equipped with its induced paring. It is easy to see that the balanced annulus decomposition of $W$ relative to $P_W$ is $(B\cap W)\setminus P_W$. 
Thus, after Lemmas \ref{lma:cut point} and \ref{Y} it only remains to study the relation between the sets in $V$ and $B$ for the subgroup $K_W$ of $K$ uniformizing such a submanifold $(W,P_W)$, 
i.e., with connected limit set without cut points.

\begin{proposition}     \label{prop:annulus tree}
Let $K$ be a {torsion free} geometrically finite Kleinian group uniformizing a pared manifold $(M,P)$ with associated representation $\rho:\pi_1(M)\to K$. Assume that $\La_K$ is connected without cut points and let $B$ be the balanced annulus decomposition of $(M,P)$. A subset $X$ of $\Lambda_K$ is:
\begin{enumerate}[- ]
    \item an inseparable cut pair if and only if there is a solid torus component $W$ of $M\setminus B$ such that $X$ is the limit set of a conjugate of $\rho(\pi_1(W))$,
    \item a necklace if and only if there is an $I$-bundle component $W$ of $M\setminus B$ such that $X$ is the limit set of a conjugate of $\rho(\pi_1(W))$,
    \item a maximal inseparable set {with at least $3$ points} if and only {if} there is a pared acylindrical component  $W$ of $M\setminus B$ such that $X$ is the limit set of a conjugate of $\rho(\pi_1(W))$.
\end{enumerate}
\end{proposition}

{Consider a Kleinian group $K$ uniformizing a pared manifold $(M,P)$ together with a homeomorphism identifying $M\setminus P$ with $(\HH^3\cup\Omega_K)/K$ and a representation $\rho:\pi_1(M)\to \PP SL_2(\C)$. Assume that $\La_K$ is connected and let $E$ be a connected component of $B$ such that $\rho(\pi_1(E))$ is generated by a parabolic isometry. Since $K$ uniformizes $(M,P)$, $E$ can be homotoped (without keeping its boundary on $\partial M$) in a component $P_i$ of $P$. Since $E$ is essential and {$(M,P)$ is a pared manifold}, at least one component of $\partial E$ is not homotopic to $P_i$ on $\partial M$. It follows that $\rho(\pi_1(E))$ is an accidental parabolic and, by Proposition \ref{prop:cutaccidental}, that $\La_K$ has a cut point.

Thus we have proved that when $\La_K$ is connected with no cut point, for any component $E$ of $B$, $\rho(\pi_1(E))$ is generated by a loxodromic isometry. In particular, since any rank $2$ Abelian subgroup of $K$ is parabolic, no component of $M\setminus B$ is a thickened torus.}

{Let $\PP$ be the paring in $K$ defined by the conjugacy classes of the parabolic subgroups. From Proposition \ref{prop:qisomgfk} we get a quasi-isometry $\Phi:\cus(K,\PP)\to \mathrm{Hull}(\La_K)$. If $E$ is a connected component of $B$, since $\rho(\pi_1(E))$ is generated by a loxodromic isometry, then $\tilde E\cap\mathrm{Hull}(\La_K)\subset\HH^3$ is quasiconvex for any connected component $\tilde E$ of the preimage of $E$. It follows that $\tilde W\cap\mathrm{Hull}(\La_K)$ is quasiconvex where $\tilde W\subset\HH^3\cup\Omega_K$ is a component of the preimage of a component $W$ of $M\setminus B$. Then $\Phi(\cus(K_{\tilde W},\PP_{\tilde W}))$ lies at bounded Hausdorff distance from $\tilde W\cap\mathrm{Hull}(\La_K)$ where $K_{\tilde W}<K$ is the stabilizer of $\tilde W$ and $\PP_{\tilde W}$ is the induced paring ({obtained from Definition \ref{def:inducedgrouppair} and Fact \ref{fact:inducedparing}}). From Fact \ref{fact:close qi}, we deduce that $\tilde W\cap\mathrm{Hull}(\La_K)$ is quasi-isometric to $\cus(K_{\tilde W},\PP_{\tilde W})$ and that the ideal boundary of $\tilde W\cap\mathrm{Hull}(\La_K)\subset (\HH^3\cup\Omega_K)$ is the limit set of {$K_{\tilde W}$} which is a conjugate of $\rho(\pi_1(W))$.
}

To prove {Proposition \ref{prop:annulus tree}} we use the work of Walsh \cite{walsh:bumping} where the author explains the relation between the bumping sets of the connected components of $\Omega_{K}$ and the characteristic submanifold of $\HH^3/K$. 

First we use the separation theorem to establish the relation between the cut pairs in $\La_{K}$ and the bumping set of $\Omega_{K}$.

\begin{lemma}   \label{cut pair}
Let $X\subset\cbar$ be a {locally connected} continuum without any cut point. A pair {$\{x,y\}\subset X$}  is a cut pair if and only if there are at least two components $\OOO$ and $\OOO'$ of $\Omega=\hat \C\setminus X$ such that ${\{x,y\}}\subset \overline{\OOO}\cap\overline{\OOO'}$.
\end{lemma}

\Proof
Since {$\{x,y\}$} is a cut pair, one can find two closed {subsets} $A$ and $B$ of {$X$ such that $X=A\cup B$ and $A\cap B=\{x,y\}$}. 
By the separation theorem \cite[Theorem VI.3.1]{whyburn:analytic_topology}, 
there is a Jordan curve $c$ that separates a point $a\in A$ 
from a point $b\in B$ such that $c\cap X={\{x,y\}}$. The set $c\setminus{\{x,y\}}$ has two connected components $k$ and $k'$. {Let us show that $k$ and $k'$ do not lie in the same connected component of $\Omega$. Otherwise, we would find an arc $k''$ in $\Omega$ joining $k$ and $k'$ so that $c\cup k''$ is a $\te$-curve. Among the boundaries of the three components of its complement, the two that contain $k''$ only intersect $X$ at a single point. But, as $A\cup B$ intersects at least two of the three components of the complement of $c\cup k''$, this would imply that $x$ or $y$ is a cutpoint --- a contradiction.} The two components $\OOO$ and $\OOO'$ of $\Omega$ containing 
$k$ and $k'$ respectively satisfy the conclusion.

{
Conversely, if there are two components $\OOO$ and $\OOO'$ of $\Omega=\hat \C\setminus X$ such that $\{x,y\}\subset \overline{\OOO}\cap\overline{\OOO'}$, consider $c=k\cup k'\cup\{x,y\}$ where $k\subset\OOO$ and $k'\subset\OOO'$ are the geodesics joining $x$ to $y$ ($\OOO$ and $\OOO'$ are Jordan domains 
equipped with their Poincaré metrics, i.e., the hyperbolic metrics obtained from the {Riemann} uniformization theorem). Since the two components $\CCC_1$ and $\CCC_2$ of $\cbar\setminus c$ intersect both $\OOO$ and $\OOO'$, both $\CCC_1$ and $\CCC_2$ intersect $X$. In particular $c$ separates $X$ and $\{x,y\}$ is a cut pair.}
\endp

Using this lemma and results of \cite{walsh:bumping}, see also \cite[\textsection 2.3]{lecuire:plissage}, we can study the stabilizer of a necklace.

\begin{lemma}   \label{necklace}
Under the hypothesis of Proposition \ref{prop:annulus tree}, a subset $X$ of $\Lambda_{K}$ 
is a necklace if and only if there is an $I$-bundle component $W$ of $M\setminus B$ such that $X$ is the limit set of a conjugate of $\rho(\pi_1(W))$.
\end{lemma}

\Proof
{Let $X$ be a necklace and let us consider a finite cyclic subset $S  =\{s_j,\ j\in\Z/n\Z\}$, $n\ge 3$, together with defining continua $M_j\subset \La_K$, $j\in\Z/n\Z$. }
By definition, any pair of points in a necklace is a cut pair. {We can argue in a similar manner to the proof of Lemma \ref{cut pair} to show that}
there are open Jordan domains $\{D_j,\ j\in\Z/n\Z\}$ such 
that $M_j\subset D_j\cup\{s_{j-1},s_j\}$, $\partial\overline{D_j}\setminus \{s_{j-1},s_j\}$ consists of two arcs in different components of $\Omega_K$, $\overline{D_j}\cap \overline{D_k}\ne \emty$ only if $|j-k|\le 1$ and $\overline{D_j}\cap \overline{D_{j+1}}= \{s_j\}$.
By construction, $\cbar \setminus \cup_{j\in\Z/n\Z} \overline{D_j}$ has exactly two components, {and since $\La_K\subset\bigcup_j M_j$, they}  are contained in $\Omega_K$ . As $D_1\setminus\La_K$ lies in two different components of $\Omega_K$, we may conclude that there are two components of $\Omega_K$ that contain $S$ on their boundaries. All the other components of $\Omega_K$ are contained in one of the domains $D_j$, so cannot contain $S$ on their boundaries.
It follows then that $X$ lies in the frontiers of two components of $\Omega_K$. By \cite[Theorem 3.1]{walsh:bumping}, $X$ is the limit set of {a conjugate of} $\rho(\pi_1(W))$ for an $I$-bundle component $W$ of $M\setminus B$. 
Conversely, for such a component $W$ of $M\setminus B$ which is a pared $I$-bundle, every pair of points in the limit set of $\rho(\pi_1(W))$ is a cut pair. It follows then from \cite[Corollary 18]{papasoglu:swenson:continua} that this limit set is a necklace.
\endp

\begin{lemma}   \label{ins cut pair}
Under the hypothesis of Proposition \ref{prop:annulus tree}, a subset $\{a,b\}$ of {$\La_K$} is an inseparable cut pair if and only if there is a solid torus component $W$ of $M\setminus B$ such that $\{a,b\}$ is the limit set of a conjugate of $\rho(\pi_1(W))$.
\end{lemma}

{Before explaining the proof, let us notice that since no component of $M\setminus B$ is a thickened torus, any component of $B$ lies in the boundary of the closure of a component of $M\setminus B$ that is a solid torus.}

\Proof
{By Lemma \ref{cut pair}, $\{a,b\}$ is a cut pair if and only if there are two components $\OOO$ and $\OOO'$ of $\Omega_K$ such that $\{a,b\}\subset Y= \overline{\OOO}\cap\overline{\OOO'}$. By \cite[Theorem 3.1]{walsh:bumping} (see also \cite[\textsection 2.3]{lecuire:plissage}), the set $Y$ is the limit set $\Lambda_{K_W}$ of a subgroup $K_W<K$ which uniformizes the compactification $W$ of a component of $M\setminus B$. By Lemma \ref{cut pair}, every pair of points in $Y$ is a cut pair and by \cite[Corollary 18]{papasoglu:swenson:continua} either $Y$ is an inseparable cut pair or a necklace. 

If $Y$ is an inseparable cut pair, then $\pi_1(W)$ is conjugated to the stabilizer of $Y=\{a,b\}$. This can only happen if $W$ is a solid torus.

If $Y$ is a necklace, $W$ is an $I$-bundle over a compact surface $F$ by Lemma \ref{necklace}. {Set $P^0_W=P\cap W$ and denote by $P_F$ the projection of $P_W^0$ on $F$ along the fibers. Let $H_F$ be a Fuchsian group uniformizing $(F,P_F)$, then $H_F<\PP SL_2(\R)<\PP SL_2(\C)$ uniformizes {$(W,P_W^0)$} as a Kleinian group. On the other hand, as previously mentioned, the subgroup $K_W$  of the
Kleinian group $K$ also uniformizes {$(W,P_W^0)$}. Let $(G_W,\PP_W)$ be the pared fundamental group of $(W,P^0_W)$. Proposition \ref{prop:qisomgfk} provides us with well-defined homeomorphisms between $\La_{H_F}$ and $\partial\cus(G_W,\PP_W)$ on the one hand and between $\La_{K_W}$ and $\partial\cus(G_W,\PP_W)$ 
on the other hand from which we get a homeomorphism between $\La_{H_F}$ and $\La_{K_W}$. The action of $H_F$ on $\HH^2$ defines a hyperbolic metric on {$F\setminus P_F$}. {Since $a\neq b$, they are joined by a geodesic in $\HH^2$. Let $l$ be the projection on $F\setminus P_F$ of this geodesic.} 
Notice that any geodesic in $F$ crossing $l$ defines a cut pair in $\La_{K_W}$ separating $a$ and $b$ in $\La_K$. On the other hand if $\{x,y\}$ is a cut pair separating $a$ and $b$, by \cite[Lemma 17]{papasoglu:swenson:continua}, $Y\cup\{x,y\}$ is cyclic and by maximality $\{x,y\}\subset Y$. It follows that $l$ is homotopic to a component of $\partial F\setminus P_F$ if and only if $\{a,b\}$ is inseparable. {Since $W$ is the compactification of a component of $M\setminus B$, the closures of the components of $\partial W\setminus\partial M$ are components of $B$. In particular if $l$ is homotopic to a component of $\partial F$ then $\{a,b\}$ is the limit set of a conjugate of $\rho(\pi_1(E))$ for a component $E$ of $B$.}
As mentioned before the proof, any component of $B$ is adjacent to a solid torus and the conclusion follows.} 

{
Conversely, if $W$ is a solid torus and $\{a,b\}$ is the limit set of $\rho(\pi_1(W))$, consider two simple closed curves $\gamma_1,\gamma_2$ on different components of {$W\cap\partial M$ }{
so that $\gamma_1\cup\gamma_2$ is isotopic on $\partial M$ to the boundary of a component of $B$. Then $\gamma_1\cup\gamma_2$ bounds an essential annulus $E$ (isotopic to the aforementioned component of $B$). If $\gamma_1\cup \gamma_2$ also bounded an annulus $E'\subset \partial M$, then, since $M$ is atoroidal, either $E\cup E'$ would bound a solid torus or $E\cup E'$ would be homotopic to a surface immersed in $\partial M$. In both cases this would contradict the fact that $E$ is essential (in particular boundary incompressible). Therefore $\g_1\cup \g_2$ does not bound an annulus in $\partial M$.}

If $\tilde\gamma_1$ and $\tilde\gamma_2$ are lifts of $\gamma_1, \gamma_2$ with endpoints $a$ and $b$, then the closure of $\tilde\gamma_1\cup\tilde\gamma_2$ separates $\cbar$. {If they were in the same component of $\Omega_K$, the curves $\gamma_1$ and $\gamma_2$ would be homotopic within $\partial M$, {contradicting the fact that they do not bound an annulus in $\partial M$,} so $\tilde\gamma_1$ and $\tilde\gamma_2$ lie in different components of $\Omega_K$.}
Hence $\{a,b\}$ is a cut pair and it follows from the paragraph {before last} that it is inseparable.}
\endp

The next lemma will allow us to show that limit sets of pared acylindrical submanifolds are inseparable.

\begin{lemma}   \label{lma:inseparable}
{Under the hypothesis of Proposition \ref{prop:annulus tree}, a} pair $\{a,b\}\subset\La_K$ is separable if and only if there is an essential annulus $E\subset M$ such that $E$ intersects the projection in $M$ of any {arc} in $\HH^3$ joining $a$ to $b$.
\end{lemma}

{When we say that two essential annuli are isotopic we require the isotopy to only go through essential annuli.}

\Proof
If there is such an essential annulus $E$ then it has a lift $\tilde E$ whose endpoints separate $a$ and $b$. This is a direct consequence of the following claim which will also be used to prove the converse.

\noindent{\bf Claim.---}
Let $E$ be an essential annulus, then $E$ intersects the projection in $M$ of any {arc} in $\HH^3$ joining $a$ to $b$ if and only if there is a component $\tilde E\subset\HH^3\cup\Omega_K\approx\widetilde{M\setminus P}$ of the preimage of $E$ such that {the closure of $\partial\tilde{E} =\tilde{E}\cap\Omega_K$ in $\cbar$}
separates $a$ and $b$.

\Proof
{
Let $E\subset M\setminus P$ be an essential annulus. If one boundary curve of $\partial E$ was homotopic into $P$, then $\rho(\pi_1(E))$ would define an accidental parabolic, hence a cut point of $\La_K$ according to Proposition \ref{prop:cutaccidental}, a contradiction. Therefore its core defines a non-trivial conjugacy 
class of a loxodromic element $g\in K$. 
Let $\tilde E$ be the lift in $\HH^3\cup\Omega_K$ of $E$ that is invariant under $g$.} 
{The closure of $\tilde{E}\cap \Omega_K$ 
is a Jordan curve in $\cbar$ that bounds the disc $\tilde E\cap\HH^3$. Let $\tilde\EEE$ denote the set of the closures of the lifts of $E$ and 
$\CCC$ the collection of their Jordan boundary curves on $\cbar$. Note that the elements of $\tilde \EEE$ are pairwise
disjoint. With the claim in mind, it is also convenient to
 consider a small toric neighborhood $T$ of $E$ whose universal cover $\tilde T$ may be identified with $\R\times [-1,1]\times [-1,1]$ and let $\tilde \TTT$
be the collection of the closures of its lifts. Its elements are pairwise disjoint as well and their traces on $\cbar$ correspond to $\CCC$.

We claim that $\tilde \TTT$ is a null sequence when we consider the hyperbolic space as the unit Euclidean ball, i.e., for any $\de>0$,  there are only finitely many components of diameter at least $\de$. To see this, let us assume there exists a sequence of distinct translates $(g_k(\tilde{T}))_k$ that 
has diameter bounded from below. Note that as the stabilizer of $\tilde T$ is cyclic (hence
quasiconvex) and malnormal, the set of its fixed points $c\cap \La_K$, $c\in\CCC$, is a null-sequence according to Lemma \ref{lma:relstructure} (3). 
Therefore, we may assume that the set of fixed points $g_k(\hbox{Fix}(g))$ is convergent towards a point $x\in \La_K$. 
Let $p\in\cbar$ be disjoint from $x$ and all elements of $\tilde\TTT$ and let us consider the hyperbolic space in the upper half-space model $\C\times \R_+$
by sending the point $p$ at infinity. We may assume that the fixed points of $g$ are $\{(-1,0),(1,0)\}$. Let $\g\subset\HH^3$ be the hyperbolic geodesic 
joining both points. For each $k$, consider the hyperbolic isometry $h_k$ that fixes the point $p$ and that maps back $g_k(\g)$ to $\g$. As $h_k g_k$ maps $\gamma$ to itself and the action of 
$g$ on $\g$ is cocompact, we may precompose each $h_k g_k$ by a suitable power of $g$ so that $\{h_kg_k\}$ is relatively compact in $\mathrm{PSL}(2,\C)$.
Therefore, the Euclidean
diameters of $h_kg_k(\tilde T)$ are uniformly bounded in $\C\times\R_+$. As each $h_k$ is affine, 
it follows that the Euclidean diameter of $g_k(\tilde T)$ is of
order the Euclidean diameter of $g_k(\hbox{Fix}(g))$, hence tends to $0$ in $\C\times\R_+$. This shows that $\tilde\TTT$ is a null sequence.

Let us now establish the claim. If some $c\in\CCC$ separates $\{a,b\}$ in $\cbar$, then any curve in $\HH^3$ joining them will intersect $\tilde{E}$
by the Jordan theorem.  Conversely, assume that $\{a,b\}$ is not separated by any  element of $\CCC$.}
{We will construct inductively an arc joining $a$ and $b$ that is disjoint from the preimage of $E$. Let us first enumerate $\tilde \TTT=\{\tilde T_j\}$ so that their diameters are decreasing with respect to the Euclidean metric when we identify
$\HH^3$ with the unit ball. Let $\kappa_1:[0,1]\to \BB^3$ be an arc joining $a$ and $b$ that is properly embedded in the unit ball and assume that $\kappa_j$ has been constructed with the same properties. 
If $\kappa_j\cap \tilde T_j=\emty$, then set $\kappa_{j+1}=\kappa_j$. Otherwise, let $s_j,t_j$ be the first and last point of $\kappa_j\cap\tilde T_j$. 
Since $\{a,b\}$ is not separated by any  element of $\CCC$, $s_j$ and $t_j$ lie on the same component of $\partial \tilde T_j\cap\HH^3$. 
The next arc $\kappa_{j+1}$ is obtained by replacing in $\kappa_j$ the arc joining $s_j$ and $t_j$ with an arc in {$\partial\tilde T_j\cap\HH^3$} with the same endpoints. {Since $\tilde \TTT$ is a null sequence of pairwise disjoint sets each of which contains a point in the unit sphere, for any interval $[x,y]\subset (0,1)$, $\kappa_1([x,y])$ intersects finitely many sets of $\tilde\TTT$. Thus, for any interval $[x,y]\subset (0,1)$, there is $J$ such that {$\kappa_1([x,y])\cap T_j = \emptyset$} 
for $j\geq J$, and so} $\kappa_{j|[x,y]}=\kappa_{J|[x,y]}$ for any $j\geq J$. 
It follows that the sequence $(\kappa_j)_j$ converges to a properly embedded arc $\kappa:[0,1]\to\BB^3$ that joins $a$ and $b$. By construction, $\kappa$ is disjoint from the preimage of $E$.}
\endp

Let $\{x,y\}$ be a cut pair separating $a$ and $b$. By the separation theorem, there is a Jordan curve $c\subset\cbar$ separating $a$ and $b$ such that $c\cap\La_K=\{x,y\}$. Let $k$ and $k'$ be the components of $c\setminus\{x,y\}$ and denote by $\OOO$ and $\OOO'$ the component of $\Omega_K$ containing $k$ and $k'$ respectively. Assume that $k$ and $k'$ are geodesic with respect to the Poincaré metrics.

If $\{x,y\}$ is an inseparable cut pair, then, by Lemma \ref{ins cut pair}, $\{x,y\}$ is the limit set of a conjugate of $\rho(\pi_1(W))$ for a solid torus component $W$ of $M\setminus B$. It follows that $k$ and $k'$ are invariant under the action of this conjugate and so they project to simple closed curves {that} 
bound an essential annulus. The conclusion follows from the claim above.

Otherwise, by \cite[Theorem 3.1]{walsh:bumping} and Lemma \ref{necklace}, $X=\overline\OOO\cap\overline{\OOO'}$ is the limit set of a conjugate of $\rho(\pi_1(W))$ where $W\subset M$ is an essential $I$-bundle. 

If $\{a,b\}\subset X$, let $(a,b)\subset\HH^3$ be the geodesic joining $a$ and $b$ and let $l$ be its projection to the base surface $F$ of $W$ {as defined in the proof of Lemma \ref{ins cut pair}. Since $\{a,b\}$ is not an inseparable cut pair, it follows from Lemma \ref{ins cut pair} that $l$ is not peripheral. In particular, there is a simple closed curve $\gamma\subset F$ that crosses $l$. The $I$-bundle over $\gamma$ provides us with the desired annulus.}

If $\{a,b\}\not\subset X$, say $a\not\in X$, let $\HHH$ and $\HHH'$ be the convex hulls of $X$ in $\OOO$ and $\OOO'$ respectively, equipped with their Poincaré metrics
{and observe that $c\subset (\HHH\cup\HHH'\cup X)$}. 
{The point $a$ lies in a component $\AAA$ of $\cbar\setminus(\HHH\cup\HHH'\cup X)$ and the frontier of $\AAA$ is a Jordan curve $c'$ that separates $a$ and $b$ since
either $c'$ coincides with $c$ or it separates $a$ from $c$ as $c\subset (\HHH\cup\HHH'\cup X)$, hence separates $a$ from $b$ as well.}
By \cite[Theorem 3.1]{walsh:bumping}, the projection of $c'\setminus X$ bounds an annulus and the conclusion follows from the claim above.
\endp

\begin{lemma}   \label{max ins}
Under the hypothesis of Proposition \ref{prop:annulus tree}, a subset $X$ of {$\La_K$ with at least $3$ points}
is a maximal inseparable set {if and only if} there is a pared acylindrical component $W$ of $M\setminus B$ such that $X$ is the limit set of a conjugate of $\rho(\pi_1(W))$.
\end{lemma}

\Proof
{Let $W$ be a component of $M\setminus B$, let $\tilde W\subset\tilde M$ be a connected component of its preimage and denote by $\La_{\tilde W}\subset\cbar$ the ideal boundary of $\tilde W\cap\mathrm{Hull}(\La_K)$. As explained earlier $\tilde W\cap \mathrm{Hull}(\La_K)$ is quasiconvex, hence $\La_{\tilde W}$ is the limit set of a conjugate of $\rho(\pi_1(W))$.

Given two points {$a,b\in \La_{\tilde W}$}, we can use the same construction as in Lemma \ref{lma:inseparable} to build an arc $l\subset \tilde W$ joining $a$ and $b$. As mentioned in the proof of Lemma \ref{lma:cut point}, any essential annulus can be made disjoint from $B$ by an isotopy (see \cite[Theorem 3.9]{bonahon:3-var}). When $W$ is acylindrical, it follows then from Lemma \ref{lma:inseparable} that $\La_{\tilde W}$ is inseparable.

On the other hand, let $X\subset\La_K$ be an inseparable set and let $\tilde B\subset\HH^3\cup\Omega_K$ be the preimage of $B$. {Replacing $E$ with $B$ in the claim inside the proof of Lemma \ref{lma:inseparable}, we may conclude that} 
any two points in $X$ are joined by a line that is disjoint from $\tilde B$.  Thus there is a component $\tilde W$ of {$(\HH^3\cup\Omega_K)\setminus\tilde B$}  
such that 
{$\tilde W$} contains a line joining $(a,b)$ for any pair $(a,b)\in X^2$. It follows that $X$ is a subset of the ideal boundary $\La_{\tilde W}$ of {$\tilde W$}. 
We have established at the beginning of the proof that $\La_{\tilde W}$ is the limit set of a conjugate of $\rho(\pi_1(W))$. If $X$ has at least $3$ points, $W$ can not be a solid or thickened torus and by Lemma \ref{necklace} it can not be an essential $I$-bundle when $X$ is also inseparable.

We have shown that a subset $X$ of $\La_K$ with at least $3$ points
is inseparable if and only if there is a pared acylindrical component $W$ of $M\setminus B$ such that $X$ lies in the limit set of a conjugate of $\rho(\pi_1(W))$. It follows that the maximal inseparable sets with at least $3$ points are exactly the limit sets of conjugates of fundamental groups of pared acylindrical components of $M\setminus B$.}
\endp

This concludes the proof of Proposition \ref{prop:annulus tree}.

Let us now recall the definition of $\TTT_{\La_K}$. Let $K$ be a {torsion free} geometrically finite Kleinian group {uniformizing a pared manifold $(M,P)$} and consider the set $V$ of cut points, {inseparable} cut pairs, necklaces and maximal inseparable sets {with at least 3 points} of $\Lambda_K$. By \cite[Thm 0.2]{bowditch:connectedness} and Proposition \ref{prop:annulus tree}, $V$ is countable. Define a graph $\TTT_{\Lambda_K}$  with vertex set $V$ by putting an edge between two vertices $v_1,v_2$ if $v_1\subset v_2$ as subsets of $\Lambda_K$. By Lemma \ref{lma:cut point}{, Lemma \ref{Y}} and Proposition \ref{prop:annulus tree}, $\TTT_{\La_K}$ is the tree dual to the balanced annulus decomposition of the pared manifold $(M,P)$ uniformized by $K$. Thus we have:

\begin{corollary}\label{cor:commonjsj}
{There is a $K$-equivariant isomorphism between the trees $\TTT_{\La_K}$ and $\TTT_{B}$.}
\end{corollary}

From now on, we will call $\TTT_{\La_K}$ the {\em JSJ tree} of $K$.
Let us notice that any homeomorphism of $\Lambda_K$ induces an isometry on $\TTT_{\Lambda_K}$.

\subsubsection{Quasi-isometries and JSJ splittings}

Let $G$ be a finitely generated group quasi-isometric to a one-ended Kleinian group $K$. {By Remark \ref{torsion free etc}}, we may assume that $K$ is torsion free, geometrically finite and minimally parabolic. To obtain the splitting from Theorem \ref{thm:dec JSJ}, we will prove that the Bowditch boundary of $G$ is homeomorphic to $\La_K$. Thus we get a simplicial action of $G$ on the JSJ tree $\TTT_{\La_K}$ and a splitting of $G$. We will then deduce properties {\em (1)} to {{\em (7)}} from the results of \textsection \ref{scn:cutstab}, \textsection \ref{sec:stab annulus} and  the quasi-isometric invariance of the JSJ splitting.

 Since $K$ is a geometrically finite Kleinian group, it is hyperbolic relative to its parabolic subgroups. 
 Since it is minimally parabolic, they are isomorphic to $\Z^2$. By \cite[Theorem 1.6]{drutu:sapir}, $G$ is hyperbolic relative to a collection {$\PP$} 
 of subgroups each of which is  quasi-isometric to $\Z^2$. It follows then from  Theorem \ref{thm:abelian} that the peripheral subgroups of $G$ 
 are virtually $\Z^2$. Furthermore, by results of Osin and Drutu-Sapir the peripheral subgroups are quasi-isometrically embedded 
 \cite[Corollary 8.3]{rhuska:qcvxrhg} and, by \cite[Theorem 1.7]{drutu:sapir}, the quasi-isometry maps each coset of a peripheral subgroup 
 in a {bounded} neighborhood of a coset of a parabolic subgroup and conversely. Thus we get:

\begin{lemma}\label{lma:parabolics}
Let $G$ be a one-ended finitely generated group quasi-isometric to a  geometrically finite minimally parabolic Kleinian group $K$, then $G$ is  hyperbolic relative to virtually Abelian rank $2$ groups. Furthermore the quasi-isometry preserves  the peripheral structures.
\end{lemma}

By Theorem \ref{thm:qisomgg}, the  quasi-isometry in Lemma \ref{lma:parabolics} can be extended to a quasi-isometry between cusped spaces $X_G$ and $X_K$ which, by Proposition \ref{prop:qisomgfk}, induces a quasi-isometry between $X_G$ and the convex hull of $\Lambda_K$. Such a quasi-isometry extends to a homeomorphism from the Bowditch boundary $\bop G$ to the limit set $\La_K$ of $K$. Thus we have:

\begin{lemma}	\label{lma:homeo boundary}
Let $G$ be a one-ended finitely generated group quasi-isometric to a  geometrically finite minimally parabolic Kleinian group $K$, then $G$ is relatively hyperbolic and its Bowditch boundary is homeomorphic to $\La_K$.
\end{lemma}

The action of $G$ on its Bowditch boundary induces a simplicial action on the JSJ tree $\TTT_{\La_K}$. Notice that by definition a vertex group of the JSJ splitting of $G$ is the stabilizer of a cut point, an inseparable cut pair, a necklace or a maximal inseparable set in $\bop G$. {Recall that two vertices $v_1,v_2$ are adjacent if and only if $v_1\subset v_2$ as subsets of $\Lambda_K$ (or the opposite). This can only happen if $v_1$ is a cut point or cut pair and $v_2$ is a necklace or maximal inseparable set {with at least 3 points}. It follows that the action of $G$ on $\TTT_{\La_K}$ has no edge inversion and that an edge group stabilizes two subsets $X,Y$ of $\Lambda_K$ where $Y$ is a cut point or cut pair and $X$ is a necklace or maximal inseparable set {with at least 3 points} containing $Y$.}
{The action of $G$ on $\TTT_{\La_K}$ induces a splitting, which we call the {\em JSJ splitting of $G$}. This in turn fits in the more general setting of JSJ decompositions of finitely generated groups established in \cite{guirardel:levitt:jsjpanorama}. The last ingredient we need to prove Theorem \ref{thm:dec JSJ} is the quasi-isometric invariance of the JSJ splitting.}

\begin{proposition}[Quasi-isometric invariance of the JSJ splitting]	\label{prop:qi JSJ}
Let $G$ be a one-ended group quasi-isometric to a minimally parabolic {geometrically finite} Kleinian group $K$ and let $X$ be a necklace or a maximal inseparable subset {with at least 3 points} of $\bop G\approx \La_K$. Let $(G_X,\PP_X)$, resp. $(K_X,\Q_X)$, be the stabilizer of $X$ in $G$, resp. in $K$, equipped with the paring induced 
 by the stablizers of cut points, {inseparable} cut pairs and parabolic points in $X$. 
Then $(G_X,\PP_X)$ is quasi-isometric to $(K_X,\Q_X)$.
\end{proposition}

Using Proposition \ref{prop:qisomcvx}, it is easy to show that $G_X$ is quasi-isometric to $K_X$, once we have established that $X$ is their limit set.

\begin{fact}    \label{fact:limit set}
Let $X\subset\bop G$ be a necklace or a maximal inseparable set {with at least 3 points} and let $G_X<G$ be its stabilizer. Then $X$ is the limit set of $G_X$.
\end{fact}

\Proof 
{If the JSJ decomposition is trivial, then there is nothing to prove. So let us assume that $X$ is a proper subset of $\bop G$.} 
By Proposition \ref{prop:annulus tree}, we know that   $X$ corresponds to the limit set  of $K_X$, so that it contains a dense collection
of cut points and/or inseparable cut pairs. But, for  each cut point or   inseparable cut pair $Y\subset X$,
{it follows from the definition of $\TTT_{\La_K}$ that there is an edge stabilizer in $G$ which stabilizes both $X$ and $Y$.}
Since $G$ is one-ended, edge groups are infinite, cf. \cite[4.A.6.6]{stallings:yale}. It follows that $Y\subset\La_{G_X}$, and by density of such sets in $X$,
that $X\subset \La_{G_X}$. Since $X$ is stabilized by $G_X$,  we finally get 
 $X = \La_{G_X}$.
\endp

\demode{Proposition \ref{prop:qi JSJ}}
By Lemmas \ref{lma:parabolics} and \ref{lma:homeo boundary}, there is a quasi-isometry $\vp:G\to K$ that preserves the parabolic structure
and that defines a homeomorphism $\partial\vp:\bop G\to \La_K$.
By {\cite[Proposition 3.4]{guirardel:levitt:splittings} (see also \cite{bigdely:wise:quasiconvexity} and
\cite[Theorem 7.1]{groff}), $G_X$ and $K_X$ are relatively quasiconvex. It follows from Fact \ref{fact:limit set} and Proposition \ref{prop:qisomcvx} that $\vp(G_X)$ is at bounded distance from $K_X$ and Fact \ref{fact:close qi} ensures that $G_X$ and $K_X$ are quasi-isometric. It remains to prove that the quasi-isometry preserves the parings.

Proposition \ref{prop:qisomcvx} already ensures that the parabolic structure is preserved, so we only  need to consider the stabilizers of {inseparable} cut pairs.
Let $P\in \PP_X$ be the stabilizer in $G_X$ of an inseparable cut pair $C\subset X$.  We know that $P$ is virtually cyclic and $P\subset G$ is a quasigeodesic with end points $C$. 
Then $\vp(P)$ is a quasigeodesic with endpoints $\partial\vp(C)$. Since $\partial\vp(C)$ is an inseparable cut pair, there are $k\in K$ and $Q\in \Q_X$ such that $kQk^{-1}$ stabilizes $\partial\vp(C)$ and $\partial\vp(X)$. Since $Q$ is cyclic, $kQ$ is a quasigeodesic with endpoints $\partial\vp(C)$. Now $\vp(P)$ and $kQ$ are quasigeodesics with the same endpoints, hence their distance is bounded.

This proves that {$(G_X,\PP_X)$ and $(K_X, \Q_X)$} are quasi-isometric as pared groups.
\endp

We are now ready to prove Theorem \ref{thm:dec JSJ}.

\demode{Theorem \ref{thm:dec JSJ}}
Let $G$ be  a finitely  generated  one-ended  group  quasi-isometric  to  a {torsion free} geometrically finite minimally parabolic Kleinian  group $K$, {cf. Remark \ref{torsion free etc}}. {By Lemma \ref{lma:parabolics}, $G$ is hyperbolic relative to some collection of subgroups $\PP$ made up of virtually Abelian rank $2$ subgroups.}
Lemma \ref{lma:homeo boundary} provides an action of $G$ on the JSJ tree $\TTT_{\La_K}$. We will show that the induced JSJ splitting has Properties {\em (1)} to {\em (7)}}, starting with the easiest properties.

Since $G$ is hyperbolic relative to subgroups virtually isomorphic to $\Z^2$ (Lemma \ref{lma:parabolics}), {virtual} Abelian subgroups are finite, virtually cyclic or virtually $\Z^2$. Since $G$ is one-ended, edge and vertex groups are infinite, cf. \cite{stallings:yale}. Property {\em (3)} follows.

{Let us establish Property {\em (1)}. As previously mentioned an edge group stabilizes in $\bop G$, two subsets $X,Y$ where $Y$ is a cut point or cut pair and $X$ is a necklace or maximal inseparable set {with at least 3 points} containing $Y$. We treat two cases depending on whether $Y$ is an inseparable cut pair or a cut point.
If an infinite subgroup of $G$ stabilizes a pair of points $\{a,b\}$, then it is virtually cyclic since it is elementary and 2-ended. Thus, this edge group is virtually cyclic, i.e. Property {\em (1)} holds in that case. 
If an infinite subgroup of $G$ stabilizes a single point, it is a subgroup of a peripheral subgroup of $G$ which is virtually $\Z^2$ by Lemma \ref{lma:parabolics}. Hence, if an edge group stabilizes a cut point $Y=\{a\}$ and  a necklace or a maximal inseparable set $X$ then it also stabilizes the connected component of $\bop G\setminus\{a\}$ containing $X$. Property {\em (1)} follows in that case from Lemma \ref{lma:paraz2}.

Let us prove Property {\em (2)}. Let $v$ be a vertex with (infinite) virtually Abelian group $G_v$. By Property {\em (3)} established above, it contains a finite index subgroup isomorphic to $\Z$ or $\Z^2$. Any edge group incident to a virtually cyclic vertex group  is commensurable  
to that vertex group since both are 2-ended; Property {\em (2)} follows in that case. Lemma \ref{lma:paraz2} {combined with the argument in the previous paragraph also establishes Property {\em (2)} when $G_v$ is commensurable to $\Z^2$.}

A group of isometries of a Gromov hyperbolic space is virtually Abelian only if it is elementary, i.e., its limit set is made up of at most two points. A consequence of  Fact \ref{fact:limit set} is that {virtual} Abelian vertex groups are precisely the stabilizers of cut points and inseparable cut pairs. Property {\em (4)} follows then from the definition of the JSJ tree.

Property {\em (5)} follows from Proposition \ref{prop:qi JSJ} and the fact that the JSJ tree is dual to the balanced annulus decomposition {by Corollary \ref{cor:commonjsj}.}

{Finally by Lemma \ref{lma:parabolics}, the maximal parabolic subgroups of $K$ are mapped by the quasi-isometry at bounded Hausdorff 
distance from those of $G$, so both
$G$ and $K$ have  the same parabolic points. Since {\em (6)} and {\em (7)} hold for $K$, by Corollary \ref{cor:commonjsj}, they also hold for $G$.}
\endp

\section{Prescription of subgroups up to {finite index}}  \label{sec:prescription}

We show how to construct {finite index} subgroups of JSJ decompositions with prescribed
vertex and edge groups. 

\subsection{Finite index subgroups with prescribed peripheral subgroups}
Let $G$ be a group. A subgroup $H<G$ is {\it separable} if, for any $g\in G\setminus H$, there exists a finite index subgroup $G'<G$ which
contains $H$ but not $g$. The group $G$ is {\it residually finite} if $\{1\}$ is separable; in other words, for any $g\ne 1$,  there exists a finite index subgroup 
$G'<G$ disjoint from $g$. Equivalently, for any $g\ne 1$,  there exists a normal finite index subgroup $G'<G$ disjoint from $g$. A group is {\em LERF} (Locally Extended Residually Finite) if any finitely generated subgroup is separable.
Two groups $G$ and $Q$ are {\it commensurable} when $G$ has a finite index subgroup isomorphic to a finite index subgroup of $Q$. Notice that if $Q$ is residually finite or virtually torsion free, then $G$ is as well.

\begin{defn}[Deep residually finite Dehn fillings]  \label{dfn:dehn filling}
A  pared group $(G,\PP)$, $\PP=\{P_1,\ldots,P_n\}$,  {\em has deep residually finite Dehn fillings} if it satisfies the 
following property $(\ast)$:
\begin{enumerate}[$(\ast)$]
\item
for each $P_j\in\PP$, there exists a  finite index subgroup $P_j^\circ <P_j$ such that, whenever $P^c_j < P^\circ_j$
is a  normal finite index subgroup of $P_j$ for each $j$, the quotient $\overline{G}=G/\ll P^c_j,\ P_j\in\PP\gg$ is residually finite and $\overline{P}_j=P_j/P^c_j$ embeds
in $\overline{G}$, where $\ll P^c_j,\ P_j\in\PP\gg$ denotes the smallest normal subgroup that contains $\{P^c_j,\ P_j\in\PP\}$.
\end{enumerate}
\end{defn}

Let us explain the analogy with {three}-manifolds. A Dehn filling of a compact {three}-manifold $M$ consists in gluing a solid torus along a toroidal boundary component of $M$. On the level of fundamental groups, if $P<\pi_1(M)$ is the fundamental group of the toroidal boundary component, a Dehn filling yields a subgroup $P^c <P$ in each $P$ such that the fundamental group of the filled manifold is isomorphic to $\pi_1(M)/ \ll P^c,\ P\in\PP\gg$. 
{Note that in this topological setting, the subgroup $P^c$ is cyclic with infinite index in $P$.}

We will use this property to find finite index pared subgroups with prescribed parings through the following claim:

\begin{cl}\label{cl:dehn} Let $(G, \{P_1,\ldots ,P_n\})$ be a pared group with deep residually finite Dehn fillings.

Then, whenever $P_j^c < P_j^{\circ}$
is a normal finite index subgroup of $P_j$ for each $j$,  there exists a normal finite index
pared subgroup $(H,\Q)$ of $(G,\PP)$ such that, for any $g \in G$ and any $j\in\{1,\ldots, n\}$,
$H \cap g P_j g^{-1} = g P_j^c g^{-1}$.
\end{cl}

\Proof  
By construction $\overline{P_j}=P_j/P_j^c$ is finite. By property $(\ast)$, $\overline{G}$ is residually finite and there is a morphism $f:\overline{G}\to Q$ to a finite group $Q$ such that
$1\notin f(\overline{P_j}\setminus\{1\})$ for all $j$. Let $H$ be the kernel of $G\to \overline{G}\to Q$. Then $H$ is   a normal finite index
subgroup of $G$ such that, for any $g \in G$ and any  $ j\in\{1,\ldots, n\}$,
{$H$}$\cap g P_{j} g^{-1} = g P_{j}^c g^{-1}$. We obtain a paring from Fact \ref{fact:inducedparing}. \endp

Our chief concern for the deep residually finite Dehn filling property ($\ast$) is to Kleinian groups, but this notion could be of interest for other classes. It is inspired by the work of Wise \cite{wise:qcvxh}.

\begin{theorem}[Residually finite Dehn filling] \label{thm:rfdf} {Let $(G,\PP)$ be a pared group and $(K,\PP_K)$ a finite index subgroup with its induced paring coming from Fact \ref{fact:inducedparing}. If $(K,\PP_K)$ is a geometrically finite Kleinian group with its usual paring, then $(G,\PP)$ has deep residually finite Dehn fillings.}
\end{theorem}

\Proof 
{We start by observing that,  by \cite[Theorem 17.14]{wise:qcvxh}, $K$, and hence $G$, is virtually compact special. Then using a generalization of the malnormal special quotient theorem \cite[Lemma 15.6]{wise:qcvxh}, see also \cite[Theorem 2]{einstein}, we establish the existence of $P_j^\circ <P_j$ such that, whenever $P^c_j < P^\circ_j$
is a  normal finite index subgroup of $P_j$ for each $j$, the quotient $\overline{G}=G/\ll P^c_j,\ P_j\in\PP\gg$ is also virtually compact special. The conclusion that $\overline{G}=G/\ll P^c_j,\ P_j\in\PP\gg$ is residually finite follows from the residual finiteness of compact special groups, \cite{haglund:wise:special} and \cite{malcev} and the fact that $P_j^\circ$ can be chosen so that $\overline{P}_j=P_j/P^c_j$ embeds
in $\overline{G}$ follows from \cite[Theorem 1.1]{osin:dehn}.}
\endp

 \subsection{Graph of groups structure with prescribed vertex groups}
  
We now use the previous results to  ``choose'' the vertex and edge groups in a graph of groups. 

\begin{proposition} \label{separability}
 Let $G$ be a finitely generated group.
We assume that
$G$ is  the fundamental group of a finite graph of groups
$\GGG=(\G,\{G_v\},\{G_e\}, j_e: G_e\hookrightarrow G_{t(e)})$ with the following properties:
\begin{enumerate}[(i)]
\item The set of vertices admits a partition $V(\G)= A\sqcup B$ into two types.
\begin{enumerate}
\item To a vertex of type $A$ corresponds  a group $G_v$ which admits a paring  $(G_v,\PP_v)$  containing the adjacent edge groups and which has the 
deep residually finite Dehn fillings property ($\ast$).
\item  To a  vertex of type $B$ corresponds a group $G_v$ such that any finite index subgroup of an adjacent edge group is separable in $G_v$.
\end{enumerate}
\item No edge has both extremities in $B$.
\end{enumerate}

Given a finite index  subgroup $H_v<G_v$ (resp. $H_e<G_e$) for each vertex group  $G_v$ (resp. edge group $G_e$), $G$ contains a normal finite index subgroup $G'$ 
which {inherits from $\GGG$ a finite graph of groups structure} $\GGG'=(\G',\{G_w'\},\{G_f'\}, G_f'\hookrightarrow G_{t(f)}')$ 
with the following properties.
\ben
\item Vertex groups are conjugate (within $G$) to {finite index}  subgroups of $H_v$.
\item Edge groups are conjugate (within $G$) to {finite index} subgroups of $H_e$.
\item {For a type B vertex $w$ of ${\mathcal G}'$, if two adjacent edge groups are commensurable, then they are equal.}
\een

{If furthermore we assume :\\
\indent\indent\indent (i) (b') For each group $H_v$ corresponding to a vertex of type $B$ if an adjacent edge group $j_e(H_e)\cap H_v$ is virtually cyclic then it is contained in a unique maximal virtually cyclic subgroup of $H_v$,\\ \\
then we may choose $G'$ so that

\indent\indent\indent (4) {for any vertex $w$ of type $B$,} if an incident edge group $G'_f$ is cyclic, then it is {\it primitive} in $G'_w$, i.e., if $g^n\in j_f(G'_f)$ with $g\in G'_w$ and $n\geq 1$ then $g\in j_f(G'_f)$.}
\end{proposition}

\Proof 
For each vertex subgroup $G_v$, we will find a finite index subgroup $G_v'$ of $H_v$, normal in $G_v$,  so that if $e\in E$ is an edge, then
$j_e^{-1}(G'_{t(e)})=   j_{\bar e}^{-1}(G'_{t(\bar e)})$.
 This will allow us to construct a quotient of $G$ in order to form a finite index subgroup of $G$  that  combines these subgroups together. 

Considering instead their normal core, we may assume that the subgroups $H_v$ and $H_e$ are normal (and of finite index) in $G_v$ and $G_e$ respectively, {and 
we may also suppose that $H_e=H_{\bar e}$}.

For $v\in A$, property ($\ast$) provides us for each $P\in\PP_v$ a finite index subgroup $P^\circ <P$ which fixes us the necessary deepness to perform controlled 
Dehn fillings. {We may assume that $P^\circ$ is also normal in $P$.}  { For $P= j_e (G_e)$, $v=t(e)$, we will also write $P_e= P$ and $P_{e}^\circ = P^\circ$}. {It is important here to remember that the edges are oriented and that the treatments 
for $e$ and $\bar e$ are independent since they depend on their terminal points.} {In particular $P_{\bar e}^\circ$ might not be defined, let alone equal to $P_e^\circ$.}

 Let $e\in E$; by condition (ii) of the statement, either a single extremity is in $A$,  and then, we may assume that it is $t(e)$, or both extremities are in $A$.
 {In the latter case, set $K_e= K_{\bar e}= H_e\cap  j_e^{-1}(H_{t(e)}\cap P_{e}^\circ) \cap j_{\bar e}^{-1}(H_{t(\bar e)}\cap P_{\bar e}^\circ)$.  In the former case, set $K_e= H_e\cap j_e^{-1}(H_{t(e)}\cap P_{e}^\circ)$. Note that, in both cases,  $K_e$ is a finite  index subgroup of $j_e^{-1}(P_e^\circ)$, contained in $H_e$ and normal in $G_e$ {(as intersection of normal subgroups)}.

{Fix $v\in B$ and let $e$ be an incident edge so that $t(e)=v$. Let  $I$ be the (finite) collection of edges $f\in t^{-1}(v)$ incident to $v$ for which $j_f(G_f)$ is 
commensurable to $j_e(G_e)$. Let $L_e =L_I = H_v \cap ( \cap_{f\in I} j_f(K_f))$; this is a finite index subgroup of $j_e(G_e)$. Therefore, since $v$ is in $B$, there is 
 a normal finite index subgroup $G_{v,e}=G_{v,I} <H_v$ in $G_v$ such that  $G_{v,e}\cap j_e(K_e)<L_e$ is a normal finite index subgroup of $j_e(G_e)$. Note that, {under the condition (i) (b'),}
 if $G_e$ is virtually cyclic, {since, by assumption, $G_{v,e}\cap j_e(K_e)$ is separable in {$H_v$}, up to taking a finite index subgroup,} we may assume that $G_{v,e}\cap j_e(K_e)$ is generated by a primitive element within $G_{v,e}$. }
 
Now, we define $G_v'= \cap_{t(e)=v} G_{v,e}$.  This is a normal finite index subgroup of $G_v$ such that, for any edge $e$ with $t(e)=v$,
$G_e' = j_e^{-1}(G'_v)$ is a  normal finite index subgroup of $G_e$ contained in $ j_{\bar e}^{-1} (P_{t(\bar e)}^\circ)$. We let
$G_{\bar e}'=G_e'$. {It follows from the previous paragraph} that if two edges $e_1$ and $e_2$ have commensurable  groups $G_{e_1}$ and $G_{e_2}$  in $G_v$, then $G_{e_1}'$ and $G_{e_2}'$ now
coincide in $G_v'$ and if {$j_e(G_e')$} is cyclic, then {it} is generated by a primitive element.}

If $e$ is an edge for which both extremities are in $A$, we set $G_e'=K_e$. 

Let us now consider $v\in A$. Let $e$ be an edge with $t(e)=v$ and let $P_e^c= j_e(G_e')$. For the peripheral subgroups $P$ in $\PP_v$ which are not 
edges, we consider any normal finite index subgroup $P^c< P^\circ$ of $P$. The deep residually finite Dehn filling property provides us through Claim \ref{cl:dehn} with 
a normal finite index subgroup $K_v{\subset G_v}$
such that $K_v\cap P= P^c$ for all $P\in\PP_v$. Set $G_v'=K_v\cap H_v$ {that is
by construction a finite index normal subgroup of $G_v$.} We note that since $j_e(K_e)\subset H_{t(e)}$ for every edge $e\in E$,
we obtain for each $v\in A$ and $e\in t^{-1}(v)$,  $G_v'\cap j_e(G_e)=  P_e^c= j_e(G_e')$. 

By construction, $j_e^{-1}(G'_{t(e)})=  G_e'= j_{\bar e}^{-1}(G'_{t(\bar e)})$ holds for each edge $e$. We may consider the graph of groups $$\overline{\GGG}=(\G,\{\overline{G}_v\},\{\overline{G}_e\}, \overline{G}_e\hookrightarrow \overline{G}_{t(e)})$$
where $\overline{G}_v= G_v/G_v'$ and $\overline{G}_e= G_e/ G_e'$.  Let $\overline{G}$ be the fundamental group of $\overline{\GGG}$,
which is a finite graph of finite groups, hence is virtually free and residually finite. 

The canonical projections $G_v\to\overline{G}_v$ define a projection $p:G\to\overline{G}$. Since $\overline{G}$ is residually finite, there is a morphism 
$\varphi:\overline{G}\to F$ to a finite group $F$ that maps a set of transversals for each $G'_v<G_v$ to {distinct} elements. 
Let $G' =\ker (\varphi\circ p)$ be the kernel of $(G\to\overline{G}\to F)$. Then $G'\cap G_v=G'_v$ for any vertex $v$. The action of $G'$ on the Bass-Serre tree
of $G$ defines a finite graph of groups structure for which vertex groups are conjugate to some subgroup $G'_v$ and edge groups are conjugate to  some subgroup $G_e'$.
\endp

We draw the following application to {three}-manifolds:

\begin{proposition}\label{prop:fcover} 
Let $M$ be a compact irreducible {non-geometric three}-manifold. There is a finite cover $N\to M$
such that  each Seifert piece in the characteristic decomposition of $N$ is the product of $S^1$ with an orientable  compact surface.
\end{proposition}

\Proof Taking a degree two cover, we may assume that $M$ is orientable. {Let us introduce
the {\em GH-decomposition} $T_{GH}\subset T_B$, which is the minimal union of components of $T_B$ such that the compactification of a component of $M\setminus T_{GH}$ is either a graph manifold or atoroidal (and hence hyperbolic). This decomposition}  provides us with a graph of groups
structure (which may be trivial)  such that each {vertex group is the fundamental group of a graph manifold or a hyperbolic manifold}. {By \cite[Lemma 3.7]{scott:geo3} and \cite[Proposition 4.4]{luecke:wu} (see also \cite[Lemma 2.1]{kapovich:leeb:npc3man}), if $N_i$ is graph manifold appearing in this decomposition then $N_i$ has a finite cover $M_i$ such that all the pieces in the characteristic torus decomposition of $M_i$ are products.} 
By \cite{hamilton:assh} and Theorem \ref{thm:rfdf}, the hypotheses of Proposition \ref{separability} are fulfilled, and  it   produces a normal finite index subgroup $G'$ so that each {Seifert} piece corresponds to a product $\Z\times \pi_1(S)$
where $S$ is an orientable surface. The corresponding cover satisfies the conclusions of the proposition.\endp

\section{Quasi-isometric rigidity of {three}-manifolds groups} \label{sec:results}

In this section, we  gather our previous results to establish Theorems \ref{thm:main1}, \ref{thm:main2} and \ref{thm:quotientman}, starting with the latter since it is used in the two other proofs.

\subsection{Quotients by finite subgroups vs finite index subgroups}	
\label{sec:residually finite}

Our main use for separability is the following proposition, see  \cite[Prop.\,7.2]{ph:unifplanar} for a proof.

\REFPROP{prop:indexsbgroup} Let $A'<A<G$ be groups with $[A:A']<\infty$ and $A'$ separable in $G$.
{Then there exists a normal finite index subgroup $H<G$ such that $H\cap A\subset A'$.}
\ENDPROP

Our goal here is to present some results to conclude that groups quasi-isometric to {three}-manifold groups are residually finite, cf. Theorem \ref{thm:quotientman}.

\begin{lemma}	\label{lem:separated}
Let $G$ be a finitely generated group and $p:G\to Q$ a surjective morphism with finite kernel. Assume that $Q$ has a graph of groups structure $(\G,\{Q_v\},\{Q_e\}, i_e:Q_e\hookrightarrow Q_{t(e)})$ with the following properties.
\ben
\item every finite index subgroup of an edge group $Q_e$ is separable in $Q$;
\item for every vertex $v$, $p^{-1}(Q_v)$ is residually finite. 
\een Then $G$ and $Q$ are commensurable.
\end{lemma}

\Proof
We will produce a morphism $q:G\to\overline{G}$ onto a virtually free group which is injective on $\ker p$, and use
the separability  of free groups to conclude.  The construction of the projection $q$ will result from taking  compatible quotients of the vertex
groups of $G$ onto finite subgroups.

The action of $Q$ on its Bass-Serre tree yields through $p$ an action of $G$, inducing 
a graph of groups structure $$(\Gamma,\{G_v=p^{-1}(Q_v)\},\{G_e=p^{-1}(Q_e)\}, j_e:p^{-1}(Q_e)\hookrightarrow p^{-1}(Q_{t(e)}))$$ for $G$ such that $i_e\circ p= p\circ j_e$ for every edge (this can be deduced for example from \cite[Theorems 12 and 13]{serre:arbre}).

By assumption each vertex group $G_v=p^{-1}(Q_v)$ is residually finite, hence there is a normal finite index subgroup $G'_v$ of $G_v$ such that $G_v'\cap  \ker p = \{1\}$.

For any edge $e$, set $G_e'=  j_e^{-1}(G_{t(e)}')\cap j_{\bar e}^{-1}(G_{t(\bar e)}')$. By construction,
the restriction $p:G_e'\to Q_e$ is injective and
the group $Q_e'= p(G_e')$ is a normal finite index subgroup of $p(G_e)=Q_e$. 
By assumption, $Q_e'$ is separable and  by Proposition \ref{prop:indexsbgroup}, there is a normal finite index subgroup $N_e$ in $Q$ such that  $N_e\cap Q_e\subset Q'_e$. 

Let $Q'= \cap_e N_e$, which has finite index in $Q$ since $N_e$ has finite index in $Q$ for any edge $e$ and $\G$ is finite. For each vertex $v$, set  $H_v= p^{-1}(Q')\cap G_v'$. This subgroup $H_v$ is normal and has finite index in $G_v$ since the same holds for  $G_v'$ in $G_v$  and for $Q'$ in $Q$; 
furthermore, we have $H_v\cap  \ker p\subset G_v'\cap  \ker p = \{1\}$.

Let us fix  an edge $e$. Since $p$ is injective on 
$j_e^{-1}(G'_{t(e)})\cup j_{\bar e}^{-1}(G_{t(\bar e)}')$ and $Q'\cap p(G_e)\subset Q_e'$,
we have  
$$j_e^{-1}(H_{t(e)}) = (p\circ j_e)^{-1}(Q')\cap G_e'= 
(p\circ j_{\bar e})^{-1}(Q')\cap G_e' =j_{\bar e}^{-1}(H_{t(\bar e)})$$ and this subgroup that we name $H_e$ is normal in $G_e$.

This implies that we may define a graph of groups $$\overline{\GGG}=(\G,\{\overline{G}_v\},\{\overline{G}_e\}, \overline{G}_e\hookrightarrow \overline{G}_{t(e)})$$
where $\overline{G}_v= G_v/H_v$ and $\overline{G}_e= G_e/H_e$.  
Let $\overline{G}$ be the fundamental group of $\overline{\GGG}$,
which is a finite graph of finite groups, hence  $\overline{G}$ is virtually free \cite[Thm. 7.3]{scott:wall} and residually finite. 

The canonical projections $G_v\to\overline{G}_v$ define a projection $q:G\to\overline{G}$ and $q$ is injective on $\ker p$.
Since $\overline{G}$ is residually finite, there is a morphism $r:\overline{G}\to F$ to a finite group which  is  injective on $q(\ker p)$.
The kernel $H= \ker (r\circ q)$ is a finite index subgroup of $G$ which avoids $\ker p\setminus\{1\}$, so  
$H$ embeds in $Q$ as a finite index subgroup.
\endp

Let $G$ be a group. {\it A hierarchy of length $0$} for $G$ is the graph of groups  made of the single vertex group $G$.
Let $n>0$, a {\it hierarchy of length at most $n$} for $G$ is a non-trivial finite graph of groups structure for $G$ together with a hierarchy of length at most $n-1$ for each vertex group. 
The hierarchy has {\it length $n$} if it has length at most $n$ and at least one vertex group has a hierarchy of length $n-1$.

Let $G_v$ be a vertex group of the graph of groups structure of $G$. We say that $G_v$ is {\it at level $1$} of the hierarchy. The groups at level $1$ for the hierarchy of $G_v$ are {\it at level $2$} of the hierarchy for $G$ and so on. The vertex groups with a hierarchy of length $0$ are called the {\it terminal groups} of the hierarchy.

Using Lemma \ref{lem:separated} inductively on a complete hierarchy, we get the following:

\begin{lemma}	\label{lem:hierarchy residually finite}
Let $G$ be a group and $p:G\to Q$ a surjective morphism with finite kernel to a group $Q$ which has a hierarchy of length $n$. A group $H$ at level at most $n-1$ in the hierarchy of $Q$ comes with a graph of groups structure $(\Gamma,\{H_v\},\{H_e\}, H_e\hookrightarrow H_{t(e)})$. Assume that for any such group $H$ (including $Q$ itself) every finite index subgroup of an edge group $H_e$ is separable in $H$. Assume furthermore that for any terminal group $T$ in the hierarchy for $Q$, $p^{-1}(T)$ is residually finite. If $Q$ is residually finite, then $G$ and $Q$ are commensurable.
\end{lemma}

\Proof
From the hierarchy for $Q$ the morphism $p$ induces a hierarchy of length $n$ for $G$ with the property that if $G_v$ is a group at level $k$ of the hierarchy for $G$ then $G_v=p^{-1}(Q_v)$ for some group $Q_v$ at level $k$ of the hierarchy for $Q$.

We are going to prove Lemma \ref{lem:hierarchy residually finite} with a finite recurrence on the level in reverse order. Our induction hypothesis is:
\begin{enumerate}[($P(k)$)]
\item Any group $G_v$ at level $k$ of the hierarchy for $G$ is residually finite and commensurable to $Q_v=p(G_v)$.
\end{enumerate}

$P(n)$ is satisfied by assumption.

Assume $P(k)$ holds for a given $k$ ($1\leq k\leq n$) and let us show that $P(k-1)$ holds. Let $G_v$ be a group at level $k-1$ of the hierarchy for $G$, by construction, $Q_v=p(G_v)$ is a group at level $k-1$ of the hierarchy for $Q$. Since $P(k)$ holds, by the assumption made on the hierarchy for $Q$, $G_v$ and $Q_v$ satisfy the assumption of Lemma \ref{lem:separated}. It follows that $G_v$ is commensurable to $Q_v=p(G_v)$. Since $Q$ is residually finite, the same is true for $Q_v$, hence for $G_v$.
\endp

We may now prove Theorem \ref{thm:quotientman}.

\begin{theorem:quotientman} Let $G$ be a finitely generated group and $p:G\to Q$ a surjective morphism with finite kernel. 
If $Q$ has a finite index subgroup isomorphic to the fundamental group of a compact $2$- or {three}-manifold $M$ then $G$ is commensurable to $Q$.
\end{theorem:quotientman}

\Proof 
Up to taking a finite index subgroup of $G$, we may assume that $Q$ is  isomorphic to the fundamental group of $M$ and that $M$ is orientable. If $G$ is finite, $G$ and $Q$ are commensurable to the trivial group, so we may assume that $G$ and $Q$ are infinite. 

If $M$ is a surface, it is easy to find a hierarchy for $Q$ with trivial terminal groups {and finitely generated edge groups} (see \cite{scott:wall} for example). According to \cite[Theorems 3.3]{scott:almostgeo,scott:almostgeo2}, $Q$ is LERF and the conclusion follows from Lemma \ref{lem:hierarchy residually finite}.

If $M$ is a compact irreducible {three}-manifold, we know from Agol \cite{agol:virtualHaken} that $M$ is virtually Haken. So, up to taking a finite index subgroup of $G$, we may assume that $M$ is Haken.

Let $T$ be the characteristic torus decomposition of $M$. If $T$ is empty,  then $M$ is Seifert or atoroidal (hence hyperbolic by Theorem \ref{thm:thurston}). 
Fundamental groups of Seifert and hyperbolic manifolds are LERF by \cite[Theorems 4.1]{scott:almostgeo,scott:almostgeo2} for Seifert manifolds and \cite[Cor.\,9.4]{agol:virtualHaken} for hyperbolic manifold (see also \cite[Corollary 17.4]{wise:qcvxh}). 
The Haken hierarchy \cite{haken:hierarchy, haken:surfaces} in $M$ provides us with a hierarchy for $Q$ with trivial terminal groups {and finitely generated edge groups}. 
Again we conclude with Lemma \ref{lem:hierarchy residually finite}.

If $T\neq\emptyset$, which also includes the case of {\em Sol } manifolds, then it induces a graph of groups structure $(\Delta,\{Q_v\},\{Q_e\}, Q_e\hookrightarrow Q_{t(e)})$ for $Q$ with edge groups $Q_e$ isomorphic to $\Z^2$ and vertex groups $G_v$ isomorphic to fundamental groups of Seifert or hyperbolic manifold $M_v$. By \cite[Theorem 1]{hamilton:assh}, any subgroup of an edge group is separable in $Q$.  By the hyperbolization theorem, every vertex manifold $M_v$ is geometric, in particular $Q_v$ is linear and hence residually finite by Malcev's theorem \cite{malcev}.
We have already proved Theorem \ref{thm:quotientman} for Seifert and hyperbolic manifolds so we know that $p^{-1}(Q_v)$ is comensurable to $Q_v$.  Thus $p^{-1}(Q_v)$ is residually finite and the conclusion follows from Lemma \ref{lem:separated}.

If $M$ contains an essential sphere, $Q$ has a graph of groups structure with trivial edge groups and vertex groups isomorphic 
to fundamental groups of compact irreducible {three}-manifolds. It follows from the previous paragraphs that Theorem \ref{thm:quotientman} holds 
for compact irreducible {three}-manifolds. Since fundamental groups of compact {three}-manifolds are residually finite \cite[Corollary 1.2]{hempel:residual}, 
the conclusion follows again from Lemma \ref{lem:separated}.  
\endp

{Lemma \ref{lem:hierarchy residually finite} can also be used to prove a strong version of the Bieberbach theorem where the classical faithfulness assumption is removed. {Notice that this generalization is already well known} and was used in the proof of the quasi-isometric rigidity of virtually Abelian groups by Cornulier-Tessera-Valette \cite{cornulier:tessera:valette:gafa07}. }

{
\begin{theorem}	\label{thm:bieberbach}
If a group $G$ acts properly discontinuously and cocompactly  by isometries on the Euclidean space $\E^n$, then $G$ is virtually $\Z^n$. 
\end{theorem}

\Proof
Let $K\triangleleft G$ be the kernel of this action, which is finite since the action is proper. Then $G/K$ acts properly, faithfully and cocompactly on $\E^n$, hence by the Bieberbach theorem, it is virtually $\Z^n$. Thus we have a morphism $p:G'\to \Z^n$ with finite kernel where $G'$ is a finite index subgroup of $G$. For any $k\geq 1$, $\Z^k$ is an HNN extension of $\Z^{k-1}$, i.e., has a graph of groups structure for which the graph has a single vertex and a single edge, 
and where both vertex and edge groups are  $\Z^{k-1}$. This leads us 
to a hierarchy for $\Z^n$ that satisfies the assumptions of Lemma \ref{lem:hierarchy residually finite}. It follows then that $G'$ and $\Z^n$ are commensurable and that $G$ is virtually $\Z^n$.
\endp
}

\subsection{Acylindrical hyperbolic manifolds and $I$-bundles}\label{sec:qi acylindrical}

The goal of this section is to establish the quasi-isometric rigidity of non-Abelian vertex groups of a JSJ splitting in the sense of \textsection \ref{scn:annulus decomposition}. More concretely we are going to show the quasi-isometric rigidity of pared fundamental groups of pared $I$-bundles and of acylindrical hyperbolic pared manifolds.

\begin{lemma}   \label{lma:ibundle}
Let $(G,\PP)$ be a pared group quasi-isometric to the pared fundamental group of a pared $I$-bundle. Then $(G,\PP)$ has a finite index pared subgroup which is the pared fundamental group of a pared $I$-bundle.
\end{lemma}

\Proof
Let $(W,\PP_W)=(F\times I,\partial F\times I)$ be a pared $I$-bundle whose pared fundamental group $(K,\PP_K)$ is quasi-isometric to $(G,\PP_G)$. If $\partial F=\emptyset$, the conclusion follows from Theorem \ref{thm:qi surfaces}. Otherwise $K$ is a free group and by Theorem \ref{thm:free}, $G$ has a free finite index subgroup $H$. The induced paring  $(H,\PP_H)$ given by Fact \ref{fact:inducedparing} is thus quasi-isometric to $(W,\PP_W)$. Considering a finite volume hyperbolic structure on the interior of $F$ and applying Proposition \ref{prop:qisomgfk}, we see that $\partial\cus (K,\PP_K)$ is homeomorphic to $S^1$. By Theorem \ref{thm:qisomgg} $\partial\cus (H,\PP_H)$ is also homeomorphic to a circle and by Proposition \ref{prop:relstructure} and \cite[Theorem 2]{otal:equivlibre}, $(H,\PP_H)$ is the pared fundamental group of a pared $I$-bundle.
\endp

To prove the quasi-isometric rigidity of pared fundamental groups of acylindrical hyperbolic manifold{s}, we first establish the rigidity of their action on their limit set.

\begin{theorem}\label{cor:bkm} 
If a  finitely generated group $G$ acts minimally as a geometrically finite convergence group 
by quasi-M\"obius mappings  on the limit set of a geometrically finite Kleinian group $K$ with infinite covolume
whose convex core has totally geodesic boundary, then {$G$ has a finite index subgroup which is isomorphic to and has the same action on $\Lambda_K$ as a finite index subgroup of $K$.} 
\end{theorem}

We start with a preliminary rigidity result concerning quasi-M\"obius mappings.
A {\it Schottky set} is the complement of a family of at least three pairwise disjoint open round discs of $\cbar$.

\begin{theorem}[Bonk, Kleiner and Merenkov  \cite{bonk:kleiner:merenkov}]\label{thm:bkm} Any quasi-M\"obius selfmap of a Schottky set
of measure zero is the restriction of a M\"obius transformation. \end{theorem}

The proof of Theorem \ref{cor:bkm} is essentially the same as Corollary 3.9 in  \cite{ph:abc-qrig}.

\demode{Theorem \ref{cor:bkm}}
Let $\La_K$ be the limit set of $K$, since $K$ is geometrically finite and its convex core has totally geodesic boundary, $\La_K$ is a Schottky set and has measure $0$ by \cite{ahlfors:fundamental} {and \cite[Proposition 4.2]{marden:kg}}.
Let $F$ be the kernel of the action of $G$ on $\La_K$ and $G'=G/F$. Since $G$ acts properly discontinuously on distinct triples in $\La_K$, $F$ is finite.
Let $G_M$ denote the set of quasi-M\"obius selfhomeomorphisms of $\La_K$.
Note that $G_M$ contains $K\cup G'$. 

According to Theorem \ref{thm:bkm}, the action of $G_M$  extends to an action of M\"obius transformations on $\cbar$. 
It is clearly discrete since any sequence which tends uniformly to the identity will have to eventually stabilize at least  three circles, 
implying that such a sequence is eventually the identity.
Since the limit sets  of $K$ and $G_M$ coincide and since $K$ is geometrically finite, $K$ has finite
index in $G_M$, cf. \cite[Theorem 1]{susskind:swarup}. Since the action of $G'$ on $\La_K$ is minimal, $\La_K$ is the limit set of $G'$. By assumption this action is geometrically finite, hence again by \cite[Theorem 1]{susskind:swarup}, $G'$ has finite
index in $G_M$. {Thus $G'\cap K$ has finite index both in $G'$ and $K$. By Theorem \ref{thm:quotientman}, $G'\cap K$ has a finite index subgroup that embeds in $G$ and the conclusion follows.} 
\endp

We generalize Theorem \ref{cor:bkm} as follows:

\begin{theorem}[pared quasi-isometric rigidity] \label{thm:paredmain} A pared group $(G,\PP_G)$  is quasi-isometric
to {the fundamental group of} an acylindrical hyperbolic pared {three}-manifold $(M,P_M)$ if and only if there is a compact hyperbolic pared {three}-manifold $(N,P_N)$ whose pared 
fundamental group is isomorphic to a finite index pared subgroup of $(G,\PP_G)$.
Moreover, $G$ is commensurable to $\pi_1(M)$. \end{theorem}

Notice that the {closed} 
and finite volume cases follow from \cite{cannon:cooper} and \cite{schwartz:qirank1} respectively and that the case of  an acylindrical pared free group has been established by Otal \cite{otal:equivlibre}.

\demode{Theorem \ref{thm:paredmain}}
Let $(G,\PP_G)$ be a pared group quasi-isometric to the pared fundamental group of an acylindrical pared hyperbolic manifold $(M,P)$. 
By Theorem \ref{thm:thurston acylindrical}, there is a geometrically finite Kleinian group $K$ whose convex core is homeomorphic to $M\setminus P$ and has totally geodesic boundary. 
As previously mentioned, when $\HH^3/K$ has finite volume, the conclusion follows from \cite{cannon:cooper} and \cite{schwartz:qirank1}. So let us assume that the volume of $\HH^3/K$ is infinite. 
Let $\PP_K$ be the paring of $K$ given by its parabolic subgroups. {With respect to the quasi-isometry between $(G,\PP_G)$ and $(K,\PP_K$),} 
the image of any coset of an element of $\PP_G$ is at bounded Hausdorff distance from a coset of an element of $\PP_K$ and conversely. 
By Theorem \ref{thm:qisomgg}, the  quasi-isometry 
between $(G,\PP_G)$ and $(K,\PP_K)$ extends to a quasi-isometry between cusped spaces $\cus(G,\PP_G)$ and $\cus(K,\PP_K)$. 
It follows that $\cus(G,\PP_G)$ is hyperbolic and that $G$ is hyperbolic relative to $\PP_G$. By Proposition \ref{prop:qisomgfk} 
the quasi-isometry between cusped spaces $\cus(G,\PP_G)$ and $\cus(K,\PP_K)$ induces a quasi-isometry between $\cus(G,\PP_G)$ 
and the convex hull of $\Lambda_K$ which extends to a homeomorphism $\Phi:\partial_{\PP_G} G  =  \partial\cus(G,\PP_G)\to\La_K$.  Since $G$ is relatively hyperbolic, its action on $\partial_{\PP_G} G$ (and hence on $\La_K$) is minimal and a geometrically finite convergence action. 
By Theorem \ref{cor:bkm}, there is a finite index subgroup $Q$ of $G$ isomorphic to a finite index subgroup $H$ of $K$ {and the isomorphism $\varphi:Q\to H$ conjugates the actions on $\Lambda_K$. Since the parings correspond to maximal parabolic subgroups, $\varphi$ preserves the parings}.
\endp

\subsection{Quasi-isometric rigidity of Kleinian groups} \label{sec:qi kleinian}
This section is devoted to the proof of Theorem \ref{thm:main2}.
We apply the previous sections to build a Kleinian group from a quasi-isometry $G\to K$. Theorems \ref{thm:dec JSJ} and \ref{thm:paredmain} 
and Lemma \ref{lma:ibundle} produce a graph of groups structure for $G$ with virtually Kleinian vertex groups. 
Then we use Proposition \ref{separability} to get a finite index subgroup which can be obtained by gluing fundamental groups of compact 
{three}-manifolds along their boundaries. This first leads us to prove the quasi-isometric rigidity of one-ended Kleinian groups: 

\begin{theorem}		\label{thm:one ended kleinian}
Let $G$ be a one-ended group quasi-isometric to a  Kleinian group. 
Then $G$ has a finite index subgroup isomorphic to a one-ended  Kleinian group.
\end{theorem}

Combining this theorem with Theorem \ref{prop:qi-ds} and Proposition \ref{prop:gluing discs}, 
we obtain the quasi-isometric rigidity of  Kleinian groups:

\begin{theorem}		\label{thm:kleinian}
Let $G$ be a group quasi-isometric to a Kleinian group. Then $G$ has a finite index subgroup isomorphic to a Kleinian group.
\end{theorem}

\Proof
{By Remark \ref{torsion free etc}, $G$ is quasi isometric to a torsion free minimally parabolic geometrically finite Kleinian group $K$. Assume that $K$ is not one-ended, in particular $\Omega_K\neq\emptyset$.
As explained in \textsection \ref{splittings finite}, when we decompose an irreducible {three}-manifold along a maximal union of essential discs, we get a graph of groups with trivial edge groups and one-ended or finite vertex groups. If we start with the Kleinian manifold $M_K$, we get a graph of groups whose vertex groups are either trivial or one-ended Kleinian groups. Let $G$ be a group quasi-isometric to $K$.
From Theorems \ref{prop:qi-ds} and \ref{thm:one ended kleinian} we get a graph of groups structure for $G$ with finite edge groups such that every vertex group is either finite or has a finite index subgroup isomorphic to a one-ended geometrically finite Kleinian $K_v$. Since $\Omega_K\neq\emptyset$, $\Omega_{K_v}\neq\emptyset$ hence $\chi(\partial M_{K_v})<0$. From Proposition \ref{prop:gluing discs} we get that $G$ has a finite index subgroup isomorphic to the fundamental group of an atoroidal {three}-manifold. {Since $G$ is quasi-isometric to $K$ which is not one-ended, this three-manifold has non-empty boundary and hence is Haken.} The conclusion follows from Theorem \ref{thm:thurston}.
\endp
}

Finally, Theorem \ref{thm:main2} is a consequence of the latter together with Proposition \ref{prop:minimal}. Now it remains to prove Theorem \ref{thm:one ended kleinian}.

\demode{Theorem \ref{thm:one ended kleinian}}
Let $G$ be a one-ended group quasi-isometric to a 
Kleinian group. By Theorems \ref{thm:dec JSJ} and \ref{thm:paredmain} and Lemma \ref{lma:ibundle}, $G$ has a graph of groups structure $(\G,\{G_v\},\{G_e\}, G_e\hookrightarrow G_{t(e)})$ with the following properties:

\ben
\item vertex groups are of essentially two types:
\ben\item pared groups with paring containing the adjacent edge {groups} and with a pared finite index subgroup isomorphic to a pared Kleinian group;
\item {virtually cyclic or virtually $\Z^2$;}
\een
\item virtually Abelian vertex groups are not adjacent;
\item edge groups are 
virtually cyclic, and edge groups incident to an Abelian vertex group are commensurable (by Lemma \ref{lma:paraz2}).
\een

By Theorem \ref{thm:rfdf} type (a) vertices here are also type A for the definition of Proposition \ref{separability}. Since Abelian groups are LERF (their subgroups are normal), type (b) vertices are also type B for the definition of Proposition \ref{separability}.
For each type (a) vertex $v$ we pick a finite index normal pared {torsion free} Kleinian subgroup $H_v<G_v$ and for each type (b) vertex $w$ 
we pick a torsion free Abelian subgroup $H_w<G_w$. By Proposition \ref{separability}, $G$ has a finite index subgroup $G'$ 
which is the fundamental group of a finite graph of groups $\GGG'=(\G',\{G_v'\},\{G_e'\}, G_e'\hookrightarrow G_{t(e)}')$ 
such that any vertex group is conjugate to a finite index subgroup of $H_v$. 
In particular $\GGG'=(\G',\{G_v'\},\{G_e'\}, G_e'\hookrightarrow G_{t(e)}')$ has the following properties:

\ben
\item Vertex groups are of essentially two types:
\ben
\item pared {torsion free} geometrically finite Kleinian groups with paring containing the adjacent edges;
\item Abelian groups isomorphic  to $\Z$ or $\Z^2$.
\een
\item Edge groups are cyclic, and edge groups incident to an Abelian vertex group all coincide.
\een

{Since any cyclic subgroup of $\Z$ or $\Z^2$ is contained in a unique maximal cyclic subgroup, 
 we may assume that edge groups  are primitive {in type B vertex groups} by Proposition \ref{separability}.}

This graph of groups structure for $G'$ enables us to build a Kleinian group. To each vertex $v$ we associate a compact {three}-manifold $M_v$ with fundamental group $G'_v$.

If $G'_v$ is isomorphic to $\Z$, $M_v$ is a solid torus and for each adjacent edge $e$, we pick an annulus $A_e$ on the boundary of $M_v$ 
such that the map $i_*:\pi_1(A_e)\to\pi_1(M_v)$ induced by the inclusion is also the map $G'_e\hookrightarrow G'_v$ defined by the graph of groups $\GGG'$. 
Since all those edge groups coincide and are primitive, we can choose the annuli corresponding to different edges to be embedded, disjoint and parallel.

If $G'_v$ is isomorphic to $\Z^2$, $M_v$ is a thickened torus $\TT\times I$. For each adjacent edge $e$, we pick an annulus $A_e\subset\TT\times\{0\}$ such that 
$i_*:\pi_1(A_e)\to\pi_1(M_v)$ corresponds to $G'_e\hookrightarrow G'_v$. Again we can choose the annuli corresponding to different edges to be embedded disjoint and parallel.

Otherwise, $G'_v$ is isomorphic to a {torsion free} geometrically finite Kleinian group $K_v$, $M_v$ is the compact manifold whose interior is homeomorphic to $\HH^3/K_v$ and the incident edge groups {and parabolic subgroups of $G$} define a paring on $\partial M_v$ corresponding to the parabolic subgroups of $K_v$. In particular to each adjacent {edge} 
$e$ is associated an annulus $A_v\subset \partial M_v$.

Given an edge $e=(v,v')$, we glue $M_v$ and $M_{v'}$ together along the annuli $A_{\bar e}\subset M_v$ and $A_e\subset M_{v'}$. The manifold thus produced does not depend on the map chosen to identify the annuli (up to homeomorphism). Doing this gluing for each edge, we get a compact {three}-manifold $M$ whose fundamental group is $G'$ \cite{scott:wall}.

By construction, $M$ is irreducible. By Lemma \ref{lma:natural pairing}, we just need to show that $M$ is atoroidal to conclude with the hyperbolization theorem.

{By property (6) of Theorem \ref{thm:dec JSJ} any rank $2$ virtually Abelian subgroup of $G$ lies in a conjugate of a vertex group of ${\mathcal G}$, hence the same is true for $G'$ and ${\mathcal G}'$. Let $A$ be a rank $2$ Abelian subgroup of $G'$, up to conjugation, we have a vertex $v$ such that $A<G'_v$. If $G'_v$ is isomorphic to $\Z^2$ then $M_v$ is a thickened torus $\TT\times I$ and $A$ is a subgroup of $\pi_1(\TT\times\{1\})$. Since we have picked all the gluing annuli on $\TT\times\{0\}$, $\TT\times\{1\}$ is a boundary component of the final manifold $M$. If $G'_v$ is isomorphic to a Kleinian group, $A$ is conjugate to a subgroup of the fundamental group of a torus $T\subset\partial M_v$.}  
{By property (7) of Theorem \ref{thm:dec JSJ} no gluing annulus can be homotoped in $T$. It follows again that $T$ is a boundary component of the final manifold $M$. We have thus proved that $M$ is atoroidal and the  conclusion follows from the hyperbolization theorem.
}
\endp

\subsection{Quasi-isometric rigidity of {three}-manifold groups}    \label{sec:qi general}
We may now combine all our previous results together to deduce Theorem \ref{thm:main1}.

\begin{theorem}		\label{thm:3-man}
Let $G$ be a group quasi-isometric to the fundamental group of a compact {three}-manifold $M$. 
Then $G$ has a finite index subgroup isomorphic to the fundamental group of a compact {three}-manifold.
\end{theorem}

\Proof
{Up to passing to a double cover, we may assume that $M$ is orientable. We may also fill any sphere in $\partial M$ with a ball. In particular $M$ has non-positive Euler characteristic.

Notice that the only geometric {three}-manifolds with negative Euler characteristic are hyperbolic. If $M$ is geometric, the conclusion follows from Theorems \ref{thm:list geom} and \ref{thm:kleinian}. Assume now that $M$ is not geometric.} 

Let us first assume that $G$ is one-ended. It follows that the fundamental group of $M$ is also one-ended so
that $M$ is irreducible and $\partial$-irreducible, cf. \S\,\ref{3 mfds}. 
{In particular, by the Geometrization theorem, Theorem \ref{thm:geometrization}, $M$ is Haken.}
{To make full use of the results of Kapovich and Leeb, let us consider the following coarser torus decomposition. The {\em Euler characteristic decomposition $T_{Eu}\subset T_B$} is the union of the components of {$T_B$} which {either} bound a component of $M\setminus T_B$ with negative Euler characteristic {or are parallel to a component of $\partial M$.} {It follows from Theorem \ref{thm:torus} that a component of {$M\setminus T_{Eu}$} either has zero Euler characteristic or is atoroidal (hence hyperbolic by Theorem \ref{thm:thurston} and Lemma \ref{lma:natural pairing}).} {By Lemma \ref{quasi-invariant} any quasi-isometry $\tilde M\to \tilde M$ preserves the Euler characteristic decomposition so Theorem \ref{prop:qi torus} still holds when we replace $T_B$ with $T_{Eu}$. It follows that}  $G$ has a graph of groups structure $\GGG=(\G,\{G_v\},\{G_e\}, G_e\hookrightarrow G_{t(e)})$ 
such that each vertex group $G_v$ is quasi-isometric to the fundamental group $Q_v$ of a compact {three}-manifold 
which is either hyperbolic or has zero Euler characteristic. Furthermore the quasi-isometry can be chosen to map the 
incident edge groups to conjugates of fundamental groups of boundary components.  In particular, edge groups are virtually 
rank $2$ Abelian by Theorem \ref{thm:abelian}.
 
By Theorems 
\ref{thm:one ended kleinian}, \ref{thm:non geom}  and \ref{thm:quotientman}, for each vertex $v$, $G_v$ has a finite index subgroup $H_v$ 
which is the fundamental group of a compact {three}-manifold {$N_v$} which is either hyperbolic or has zero Euler characteristic. {Let $\PP_v$ be the collection of edge groups associated to the edges adjacent to $v$, let $(H_v,\Q_v)$ be the {adorned subgroup} induced from $(G_v,\PP_v)$ and assume that $G_v$ is not Abelian.  
By Fact \ref{fact:inducedparing}, $(H_v,\Q_v)$ is quasi-isometric to $(G_v,\PP_v)$, which, by Theorem \ref{prop:qi torus}, is quasi-isometric to $(\pi_1(W(v)),P(v))$ where $W(v)$ is a component of {$M\setminus T_{Eu}$} and $P(v)$ is the collection of the fundamental groups of the tori of $\partial \bar W(v)$. 
We equip $N_v$ and $W(v)$ with non-positively curved metrics, cf. Proposition \ref{nonpositive}. From the \v{S}varc--Milnor lemma we get a quasi-isometry between the universal covers $\tilde W(v)$ of $\bar W(v)$ and $\tilde N_v$ of $N_v$. By \cite[Corollary 4.8]{kapovich:leeb:qi3man}, such a quasi-isometry maps boundary flats at bounded distance from boundary flats. It follows that each subgroup $Q_i$ in $\Q_v$ stabilizes, in $\tilde N_v$, a quasi-flat  at bounded distance from a boundary flat $F_i$ of $\tilde N_v$. An element of $H_v$ acting as a covering transformation of $\tilde N_v$ maps boundary flats to boundary flats. Since the Hausdorff distance between two different boundary flats is infinite (see \textsection \ref{scn:torus decomposition}) and $Q_i$ stabilizes a quasi-flat  at bounded distance from  $F_i$, $Q_i$ also stabilizes $F_i$. Notice that $Q_i$ might not be the whole stabilizer $\mathrm{Stab}(F_i)< H_v$, 
but, since $Q_i$ is quasi-isometric to an element of $\PP_v$, it is a rank $2$ Abelian group and hence has finite index in  $\mathrm{Stab}(F_i)$ which is also a rank $2$ Abelian group. On the other hand each boundary flat of $\tilde W(v)$ is stabilized by a subgroup conjugated to an element of $P(v)$. We have thus established that the stabilizer of each boundary flat of $\tilde N_v$ has a finite index subgroup that is conjugated to an element of $\Q_v$. By \cite{hamilton:assh} up to replacing $H_v$ with a finite index subgroup, and replacing $N_v$ with the corresponding finite cover, we may assume that $\Q_v$ is the collection of the fundamental groups of the tori of $\partial N_v$.}

By Theorem \ref{thm:rfdf}, vertices $v$ with {$W(v)$} hyperbolic are type A for Proposition \ref{separability}. Since edge groups are Abelian, it follows from \cite{hamilton:assh} that vertices $v$ such that {$W(v)$} has zero Euler characteristic are type B. 
{It also follows from the definition of the Euler characteristic decomposition that if two vertices $v$ and $w$ are adjacent, then for at least one, say $w$,  $W(w)$ is either hyperbolic or is a thickened torus that contains a component of $\partial M$. In the latter case, $G_w$ embeds into $G_v$ and removing the vertex $w$ (and its adjacent edge) would not change the group $G$. The purpose of those redundant vertices is to guarantee that the collection $P(v)$ of the fundamental groups of the tori of $\partial \bar W(v)$ is also, up to conjugacy, the collection of adjacent edge groups in the graph of groups structure of $\pi_1(M)$ associated to $T_{Eu}$  (this convention simplifies previous statements). 
Thus we may assume that two type B vertices are not adjacent  in the realm of Proposition \ref{separability}.}

Hence $\GGG$ satisfies the hypothesis of Proposition \ref{separability}, so it follows that $G$ has a finite index subgroup $G'$ such that for any vertex $v$, $G'_v{=G'}\cap G_v$ is a normal finite index subgroup of $H_v$. The subgroup $G'_v$ is the fundamental group of a covering {$N'_v$ of $N_v$}. 
The subgroup $G'$ inherits from $\GGG$ a graph of groups structure $\GGG'=(\G',\{G'_v\},\{G'_e\}, G'_e\hookrightarrow G'_{t(e)})$ 
such that any vertex group $G'_v$ is the fundamental group of a compact {three}-manifold {$N'_v$} and incident edge groups are conjugate to fundamental groups of boundary components.

For each edge $e=(v,w)$ of $\G'$ we glue {$N'_v$} to {$N'_w$} along the components of their boundaries corresponding to $e$. This produces a compact {three}-manifold whose fundamental group is $G'$ and the conclusion follows.

When $G$ is two-ended, then it is virtually cyclic and hence has a subgroup isomorphic to the fundamental group of a solid torus.

When $G$ has infinitely many ends, the conclusion follows from the one-ended case together with Proposition \ref{prop:gluing discs} and Theorem \ref{prop:qi-ds}.
\endp

{Let $M$ be a compact orientable three-manifold. Let $D\subset M$ be a maximal union of non isotopic essential spheres and discs. Given the compactification $N$ of a component of $M\setminus D$, we glue a ball along each component of $\partial N$ that is a sphere. We call the resulting three-manifold a {\em piece} of $M\setminus D$. Let $T$ be the union of the  {balanced} characteristic torus decompositions of the pieces of $M\setminus D$ and let $A$ be the union of the balanced annulus decompositions of the pieces of $M\setminus (D\cup T)$. We call {$D\cup T\cup A$} the {\em full decomposition of $M$}. Let $G$ be a finitely generated group quasi-isometric to $\pi_1(M)$. Assuming that no piece of $M\setminus D$ is geometric, we can apply Theorems \ref{prop:qi-ds}, \ref{prop:qi torus} and \ref{thm:dec JSJ}, to get a graph of groups structure for $G$ such that each vertex group is quasi-isometric to the fundamental group of a piece of $M\setminus (D\cup T\cup A)$. Since no piece of $M\setminus D$ is geometric, each piece of {$M\setminus (D\cup T\cup A)$} is either a pared acylindrical hyperbolic three-manifold, a Seifert manifold with boundary or an $I$-bundle over a compact surface with boundary. It follows then from Theorems \ref{thm:paredmain}, \cite[Theorem A]{neumann:commensurability:graph} and more classical results (see for example \cite[Theorem 4.2]{bowditch:course}) that each vertex group is commensurable to the fundamental group of a piece of $M\setminus (D\cup T\cup A)$. Thus we have:

\begin{theorem}     \label{thm:commensurability}
Let $M$ be a compact orientable $3$-manifold and let $G$ be a finitely generated group quasi-isometric to $\pi_1(M)$. Let $F=D\cup T\cup A$ be the full decomposition of $M$ and assume that $M\setminus D$ has no geometric piece. Then $G$ has a graph of groups structure such that each vertex group is commensurable to the fundamental group of a piece of $M\setminus F$.
\end{theorem}

 If a piece $N$ of $M\setminus D$ is geometric, we can still apply Theorem \ref{prop:qi-ds} and the characteristic torus (and annulus) decomposition of $N$ is empty unless $N$ is a Sol manifold. We deduce then from \cite[Theorem A]{neumann:commensurability:graph} and \cite{schwartz:qirank1} that Theorem \ref{thm:commensurability} is true when $M\setminus D$ has some geometric pieces as long as those pieces are not Sol manifolds nor closed hyperbolic manifolds. 

}

\bibliographystyle{math}
\bibliography{refs}

\end{document}